\documentclass[oneside]{amsart}
\numberwithin{equation}{section}
\usepackage[margin=1in]{geometry}
\usepackage{microtype}
\usepackage[light,onlyrm,notext,nosf,sfmathbb,frenchstyle]{kpfonts}
\DeclareMathAlphabet\mathfrak{U}{euf}{m}{n}
\SetMathAlphabet\mathfrak{bold}{U}{euf}{b}{n}
\newcommand{\myName}{John Calabrese}
\newcommand{\myTitle}{Donaldson-Thomas Invariants and Flops}

\newcommand{\correct}[1]{\underline{#1}}
\usepackage[british]{babel}
\usepackage[utf8]{inputenc}
\usepackage{hyperref}
\hypersetup{colorlinks=false, pdfauthor={John Calabrese}, pdfcreator={John Calabrese}, pdfsubject={math.ag}, pdftitle={Donaldson-Thomas Invariants and Flops}, linkcolor=blue, citecolor=brown}
\usepackage{tikz}
\usetikzlibrary{arrows,matrix,shapes,decorations}
\usepackage[backgroundcolor=yellow,textsize=tiny]{todonotes}
	\DeclareMathOperator{\Ad}{Ad}

	\DeclareMathOperator{\Aut}{Aut}
	\DeclareMathOperator{\ch}{ch}
	
	\DeclareMathOperator{\Coh}{Coh}
	
	\DeclareMathOperator{\coker}{coker}
	
	\DeclareMathOperator{\cone}{cone}
	\DeclareMathOperator{\DT}{DT}

	\DeclareMathOperator{\exc}{exc}
	\DeclareMathOperator{\Ext}{Ext}
	\DeclareMathOperator{\lExt}{\underline{Ext}}

	\DeclareMathOperator{\Hilb}{Hilb}
	\DeclareMathOperator{\Hom}{Hom}
	\DeclareMathOperator{\lHom}{\underline{Hom}}

	\DeclareMathOperator{\id}{id}

	\DeclareMathOperator{\PT}{PT}
	
	\DeclareMathOperator{\Sch}{Sch}
	
	\DeclareMathOperator{\Spec}{Spec}
	
	\DeclareMathOperator{\St}{St}
	\DeclareMathOperator{\supp}{supp}
	
	\DeclareMathOperator{\topp}{top}
	\DeclareMathOperator{\Tor}{Tor}

	\newcommand{\st}{\,\middle\vert\,}
	
		\newcommand{\longiff}{\Longleftrightarrow}
		\newcommand{\into}{\hookrightarrow}
		\newcommand{\ot}{\leftarrow}
		\newcommand{\from}{\leftarrow}

		\newcommand{\longto}{\longrightarrow}
		\newcommand{\onto}{\twoheadrightarrow}
		
		\newcommand{\longimplica}{\Longrightarrow}
		\newcommand{\dvert}{\vert^{L}} 
		
		\newcommand{\A}{\mathbb{A}}
		\newcommand{\C}{\mathbb{C}}

		\renewcommand{\D}{\mathbb{D}}

		\renewcommand{\L}{\mathbb{L}}

		\newcommand{\curlyO}{\mathcal{O}}
		\renewcommand{\O}{\mathcal{O}}

		\newcommand{\PP}{\mathbb{P}}
		
		\newcommand{\Q}{\mathbb{Q}}

\newcommand{\R}{R}
		\newcommand{\Z}{\mathbb{Z}}
		\newcommand{\cmplx}[1]{\csh{#1}}
		\newcommand{\scmplx}[1]{\csh{#1}}
		\newcommand{\curly}[1]{\mathscr{#1}}
		
		\newcommand{\curvy}[1]{\mathcal{#1}}
		
		\newcommand{\sh}[1]{{#1}}
		\newcommand{\csh}[1]{{#1}}
		\newcommand{\stack}[1]{\mathfrak{#1}} 		
		\newcommand{\astack}[1]{\mathfrak{#1}}
		\newcommand{\acat}[1]{\mathcal{#1}}
\usepackage{thmtools}
\declaretheoremstyle[
	spaceabove=6pt, spacebelow=6pt,
	headfont=\normalfont\bfseries,
	notefont=\mdseries,
	notebraces={(}{)},
	bodyfont=\normalfont,
	postheadspace=.5em,
	headformat=swapnumber
]{tstyle}
\declaretheorem[numberlike=equation,name=Theorem,style=tstyle]{teo}
\declaretheorem[numberlike=equation,name=Lemma,style=tstyle]{lem}
\declaretheorem[numberlike=equation,name=Definition,style=tstyle]{defn}
\declaretheorem[numberlike=equation,name=Corollary,style=tstyle]{cor}
\declaretheorem[numberlike=equation,name=Proposition,style=tstyle]{prop}
\declaretheorem[numberlike=equation,name=Situation,thmbox=M]{situation}
\declaretheoremstyle[
	spaceabove=6pt, spacebelow=6pt,
	headfont=\normalfont\itshape,
	notebraces={(}{)},
	bodyfont=\normalfont,
	postheadspace=.5em,
	headformat=swapnumber
]{rstyle}
\declaretheorem[numberlike=equation,name=Remark,style=rstyle]{rmk}
\declaretheorem[numbered=no,name=Remark,style=rstyle]{rmk*}
\declaretheorem[numbered=no,name=Notation,style=rstyle]{notation}
\declaretheoremstyle[
	spaceabove=0pt, spacebelow=6pt,
	headfont=\small\normalfont\itshape,
	notefont=\small\normalfont\itshape,
	notebraces={(}{)},
	bodyfont=\small\normalfont,
	postheadspace=.5em,
	headpunct=:,
	qed=$\blacksquare$
]{pstyle}
\declaretheorem[numbered=no,name=Proof,style=pstyle]{prf}
\newcommand{\perv}[2]{{^{\scriptscriptstyle #1}{#2}}}
\DeclareMathOperator{\Per}{Per}
\newcommand{\mum}{\stack{Mum}}
\DeclareMathOperator{\Pilb}{Pilb}
\DeclareMathOperator{\Filb}{Filb}
\newcounter{crucialcounter}
\newcommand{\derp}{\leq 1}

\newcommand{\schnitzel}{\leq 1}

\usepackage{tikz-cd}

\begin{document}
	\author{{\tiny \text{\myName}}}
	\address{Mathematical Institute\\ University of Oxford\\ 24-29 St.~Giles'\\ OX1 3LB Oxford\\ UK}
	\email{\texttt{calabrese@maths.ox.ac.uk.com}}
	\title{{\Large \normalfont \scshape \myTitle}}
	\begin{abstract}
		We prove a comparison formula for the Donaldson-Thomas curve-counting invariants of two smooth and projective Calabi-Yau threefolds related by a flop.
		By results of Bridgeland any two such varieties are derived equivalent.
		 Furthermore there exist pairs of categories of perverse coherent sheaves on both sides which are swapped by this equivalence.
		Using the theory developed by Joyce we construct the motivic Hall algebras of these categories.
		These algebras provide a bridge relating the invariants on both sides of the flop.
	\end{abstract}
	\maketitle
	\tableofcontents
	\section*{Introduction} 
\label{sec:introduction}
An interesting question in Donaldson-Thomas (DT) theory is whether there exists a relationship between the DT numbers of two birational Calabi-Yau\footnote{For us, a \emph{Calabi-Yau threefold} $Y$ will be a smooth and projective complex variety of dimension three with trivial canonical bundle $\omega_Y \cong \O_Y$ and torsion fundamental group $H^1(Y,\O_Y) = 0$.} threefolds.
As Calabi-Yau varieties are minimal models, any birational map between them can be broken down into a sequence of flops \cite{kawamata}.
Hence (in principle) it suffices to understand what happens in the case of a single flop.

A flop is a birational morphism fitting in a diagram
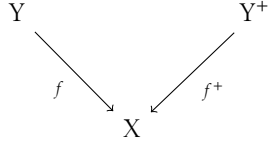
\begin{center}
	\begin{tikzpicture}
		\matrix (m) [matrix of math nodes, row sep=3em, column sep=3em, text height=1.5ex, text depth=0.25ex]
		{
		Y && Y^+\\
		&X&\\
		};
		\path[->,font=\scriptsize]
		(m-1-1) edge node[auto,swap]{$f$} (m-2-2)
		(m-1-3) edge node[auto]{$f^+$} (m-2-2)
		;
	\end{tikzpicture}
\end{center}
where $f$ (respectively $f^+$) is birational and contracts trees of rational curves to points.
In this setting it is possible to write down an explicit formula relating the generating series for the DT invariants of both $Y$ and $Y^+$, as we now explain.

We recall some notation.\footnote{The reader interested in more background on DT theory and, more generally, in curve-counting might turn to \cite{13}.}
For a curve class $\beta \in N_1(Y)$ and an integer $n \in \Z$, we can define $\DT_Y(\beta,n)$, the \emph{$\DT$ number of class $(\beta,n)$ of $Y$}, as the weighted topological Euler characteristic (where the weight is given by Behrend's microlocal function \cite{behrend}) of the Hilbert scheme $\Hilb_Y(\beta,n)$ parameterising quotients of $\O_Y$ of Chern character $(0,0,\beta,n)$.
We formally gather all the DT numbers into a series
\begin{align*}
	\DT(Y) := \sum_{\beta,n} \DT_Y(\beta,n)q^{(\beta,n)}
\end{align*}
where $q$ is a formal variable.

With a flopping contraction $f\colon Y \to X$ one can also associate the DT series of curves \emph{contracted} by $f$, that is
\begin{align*}
	\DT_{\exc}(Y) := \sum_{\substack{\beta, n \\ f_*\beta = 0}} \DT_Y(\beta,n)q^{(\beta,n)}
\end{align*}
where the subscript $_{\exc}$ stands for \emph{exceptional}.
For a pair of Calabi-Yau threefolds $Y, Y^+$ related by a flop we prove the following result.
\begin{teo}[continues=maino]
	If we define the series
	\begin{align*}
		\DT_{\exc}^\vee(Y) := \sum_{\substack{\beta,n \\ f_*\beta=0}}  \DT_Y(-\beta,n)q^{(\beta,n)}
	\end{align*}
	then the following formula\footnote{The rigorous meaning of which is explained in Remark \ref{techfinal}.} holds:
	\begin{align}\label{superid}\tag{$\star$}
		\DT_{\exc}^\vee(Y) \cdot \DT(Y) = \DT_{\exc}^\vee(Y^+) \cdot \DT(Y^+)
	\end{align}
	where we identify the formal variables $q^{(\beta,n)}$ via the flop.\footnote{As the groups of divisor classes on $Y$ and $Y^+$ are isomorphic one identifies curve-classes by dualising. This is spelled out just below \eqref{identification}.}
\end{teo}

The key ingredient here is Bridgeland's derived equivalence between $Y$ and $Y^+$ \cite{tom-flops}, which we denote by $\Phi$.
Inside the derived category $D(Y)$ of $Y$ there is a t-structure whose heart $\Per(Y/X)$ is called the category of \emph{perverse coherent sheaves}.\footnote{As a matter of fact there are different versions of this category, depending on an integer called \emph{perversity}. We shall momentarily ignore this. It will all be made clear in the following section.}
This category is intimately related to the geometry of the flop.
In fact, one can construct $Y^+$ as a moduli space of \emph{point-like objects} in $\Per(Y/X)$.
If one defines $\Per(Y^+/X)$ to be the category of perverse coherent sheaves for $Y^+$, then $\Phi$ restricts to an equivalence of abelian categories $\Per(Y^+/X) \cong \Per(Y/X)$.
This fact can be exploited to compare DT invariants on both sides of the flop and we now explain how.

It turns out that the structure sheaf of $Y$ is a perverse coherent sheaf, $\O_Y \in \Per(Y/X)$.
One can then construct a moduli space $\perv{p}{\Hilb(Y/X)}$ of quotients (in the abelian category $\Per(Y/X)$) of $\O_Y$.
If again we fix a curve class $\beta$ and an integer $n$, it is legitimate to define a \emph{perverse} DT number $\perv{p}{\DT}_{Y/X}(\beta,n)$ as the weighted topological Euler characteristic\footnote{To be precise, our conventions introduce a sign, which will be explained at the beginning of Subsection \ref{sub:conclusion}.} of the moduli space $\perv{p}{\Hilb}_{Y/X}(\beta,n)$ parameterising perverse quotients of $\O_Y$ of class $(\beta,n)$.
We formally write down a generating series for these perverse DT numbers.
\begin{align*}
	\perv{p}{\DT(Y/X)} := \sum_{\beta,n} \perv{p}{\DT_{Y/X}(\beta,n)}q^{(\beta,n)}
\end{align*}
From the discussion so far it's not clear how $\perv{p}{\DT(Y/X)}$ is related to ordinary DT numbers.
However, if we define analogously $\perv{q}{\DT(Y^+/X)}$ on $Y^+$ (and once we know that $\Phi(\O_{Y^+}) = \O_Y$) it follows immediately that $\perv{q}{\DT(Y^+/X)}$ matches up with $\perv{p}{\DT(Y/X)}$ under the equivalence $\Phi$.

To complete the picture, we will prove that $\perv{p}{\DT(Y/X)}$ is (almost) equal to the left hand side of equation \eqref{superid}.
To do this, we will use the incarnation of motivic Hall algebras found in \cite{tom-hall}, but adapted to the category $\Per(Y/X)$ (see also \cite{koso, joycobooko}).
Perverse coherent sheaves are complexes $E \in \Per(Y/X)$ concentrated in degrees $[-1,0]$.
Moreover, $H^{-1}(E)[1]$ and $H^0(E)$ are also perverse coherent, so any $E$ sits in a canonical exact sequence
\begin{align*}
	H^{-1}(E)[1] \into E \onto H^0(E)
\end{align*}
of perverse coherent sheaves.
The Hall algebra is designed to encode precisely this kind of information.
Given an epimorphism of perverse coherent sheaves $\O_Y \onto E$ we obtain a surjection of sheaves $\O_Y \to E \to H^0(E)$.
In a nutshell, these latter surjections know about the ordinary DT invariants of $Y$, while $H^{-1}(E)$ is relevant for $\DT_{\exc}^\vee(Y)$.

The case of a contraction of a disjoint union of $(-1,-1)$-curves was originally dealt with by Hu and Li in \cite{huli}.
Our formula owes a lot to work of Toda, who gave a different approach in \cite{toda-dt2}, using Van den Bergh's non-commutative resolution of $X$ \cite{vdb} and wall-crossing techniques.\footnote{We should also mention \cite{toda-gen}, where one can already find the idea of exploiting the dualising functor to study related problems.}
Our identity \eqref{superid} is related to \cite[Theorem 5.8]{toda-dt2} via \cite[Theorem 5.6]{toda-dt2} (with slightly different notation).
Strictly speaking, Toda's result works only for the naive counting invariants (defined using the ordinary, unweighted, Euler characteristic).
His proof relies on a yet unproved (but widely believed to be true) result regarding the local structure of the moduli stack of the objects of the derived category \cite[Remark 2.32]{toda-dt2}.
Finally, we are delighted to mention that the flop formula has proved to be useful in the remarkable work of Maulik \cite{davesh}.


\subsubsection*{Outline} In the first section we recall what we need about flops and construct the moduli stack of Bridgeland's perverse coherent sheaves.
The second section is devoted to checking that the theory of motivic Hall algebras applies to perverse coherent sheaves.
The third section contains the main result and its proof.
For the reader interested not only in flops but in the more general setup where the base $X$ is allowed one-dimensional singularities, we included a final section to show how $DT_{\exc}(Y/X)$ compares to $DT_{\exc}(Y)$. 
We relegated to the appendix a few basic facts about Lieblich's moduli stack of objects of the derived category \cite{lieblich}.

\subsubsection*{Acknowledgements} I would like to thank my advisor Tom Bridgeland for suggesting the problem and for providing invaluable help overcoming many technical difficulties.
Conversations with Dominic Joyce, Richard Thomas, Max Lieblich, Ed Segal, Mattia Talpo, Fabio Tonini and Arend Bayer were quite helpful and 
I am very grateful to Roland Abuaf for being up-to-date in editorial matters.
I would also like to thank the referee for helpful comments.
Among the things I would like to say to J{\o}rgen Rennemo is ``thank you'' for spotting a mistake in the proof of Proposition \ref{sucaiurgen}.

\subsubsection*{Conventions} In what follows $\C$ will denote the field of complex numbers and all stacks and morphisms will be over $\C$.
Given a scheme $(X,\curlyO_X)$ we denote by $\text{D}(\curlyO_X)$ the derived category of $\curlyO_X$-modules, by $\text{D}(X) = \text{D}_\text{coh}(\curlyO_X)$ the derived category of $\curlyO_X$-modules with coherent cohomology. 
By $D^{[a,b]}(X)$ we shall mean the subcategory of $D(X)$ consisting of complexes with cohomology concentrated between $a$ and $b$.
By $D^{\leq n}(X)$ we mean $D^{[-\infty,n]}(X)$, and similarly $D^{\geq n}(X) = D^{[n,\infty]}(X)$.
We shall sometimes write $D^{[n]}(X)$ or $\Coh(X)[-n]$ for $D^{[n,n]}(X)$.
Given a complex $\cmplx{E} \in \text{D}(\O_X)$ we denote by $H^i(E) \in \O_X\text{-Mod}$ the $i$-th cohomology sheaf and by $H^i(X,\cmplx{E}) = R^i\Gamma(X,\cmplx{E})$ the $i$-th (hyper)cohomology group.
Whenever we have a diagram of schemes $T \stackrel{u}{\to} S \stackrel{\pi}{\from}X$ we often denote a fibre product as $X_T$ together with induced maps $\pi_T:X_T \to T,$ $u_X: X_T \to X.$
The derived pullback $Lu_X^*\cmplx{E}$ of an object $\cmplx{E} \in \text{D}(\curlyO_X)$ will simply be denoted by $\cmplx{E}\vert_{X_T}^L$.
All schemes (and all algebraic stacks) will be assumed to be locally of finite type over $\C$.

\section{Flops} 
\label{sec:flops}
In this section we recall a few facts about the categories of perverse coherent sheaves and construct the corresponding moduli spaces.
\subsection{Perverse Coherent Sheaves} 
\label{sub:perverse_coherent_sheaves}
Henceforth we assume to be working within the following setup.
\begin{situation}\label{situation}
	Fix a smooth and projective variety $Y$ of dimension three, over $\C$, with trivial canonical bundle $\omega_Y \cong \O_Y$ and satisfying $H^1(Y,\O_Y) = 0$.
	Fix a map $f: Y \to X$ satisfying the following properties:
	\begin{itemize}
		\item $f$ is birational and its fibres are at most one-dimensional;
		\item $X$ is projective and Gorenstein;
		\item $Rf_* \O_Y = \O_X$.
	\end{itemize}
\end{situation}
We point out that these assumptions allow for $f$ to contract a divisor to a curve, which will be of importance in \cite{cala}.
Notice that it follows that $X$ has rational singularities \cite{kovacs}, that its canonical bundle is trivial, $\omega_X \cong \O_X$ and that $f$ is crepant.
Also, for any sheaf $G$ on $Y$, $R^i f_* G = 0$ for $i \geq 2$.
\begin{rmk}
	For our main result, Corollary \ref{maino}, we will specialise to the case where $X$ has zero-dimensional zero locus.
	The reader purely interested in flops can go ahead and add this assumption throughout.
	However, this more general setup (together with its baggage of complications) is relevant in \cite{cala}: a wall-crossing-type formula between $DT(Y/X)$ and $DT(Y)$ still holds.
	This is why we presently do not impose additional conditions on the singular variety $X$.
\end{rmk}
The main protagonist of this paper is Bridgeland's category of perverse coherent sheaves $\perv{p}{\Per(Y/X)}$ of $Y$ over $X$.
As mentioned in the introduction there are different versions of it, indexed by an integer $p$ called the \emph{perversity}.
We shall only need two of them, corresponding to the $-1$ and $0$ perversity.
One way to define these categories is by using a torsion pair \cite{torsion}, which we now recall (see also \cite[Section 3]{vdb}).
\begin{notation}\label{notationfory}
	For compactness we will often denote $\Coh(Y)$ by $\acat{A}$ and $\perv{p}{\Per(Y/X)}$ by $\perv{p}{\acat{A}}$.
\end{notation}
Let
\begin{align*}
	\curvy{C} = \left\{ \sh{E}\in \Coh(Y) \st Rf_* \sh{E}=0 \right\}
\end{align*}
and consider the following subcategories of $\acat{A}$:
\begin{align*}
	\perv{0}{\curvy{T}} &= \left\{ \sh{T} \in \acat{A} \st R^1f_* \sh{T} = 0 \right\}\\
	\perv{-1}{\curvy{T}} &= \left\{ \sh{T} \in \acat{A} \st R^1f_* \sh{T} = 0, \Hom(\sh{T},\curvy{C})=0 \right\}\\
	\perv{-1}{\curvy{F}} &= \left\{ \sh{F} \in \acat{A} \st f_*\sh{F}=0 \right\}\\
	\perv{0}{\curvy{F}} &= \left\{ \sh{F} \in \acat{A} \st f_*\sh{F}=0, \Hom(\curvy{C},\sh{F})=0 \right\}.
\end{align*}
The pair $(\perv{p}{\acat{T}},\perv{p}{\acat{F}})$ is a torsion pair on $\acat{A}$, for $p = -1 ,0$, and the tilt of $\acat{A}$ with respect to it is the category of perverse coherent sheaves $\perv{p}{\acat{A}}$.
Notice that we picked the convention where
\begin{align*}
	\perv{p}{\curvy{F}[1]} \subset \perv{p}{\acat{A}} \subset \text{D}^{[-1,0]}(Y).
\end{align*}
We mention in passing that the structure sheaf is perverse coherent, $\O_Y \in \perv{p}{\acat{T}} \subset \perv{p}{\acat{A}}$.
\begin{notation}
	For convenience (and unless otherwise stated) we shall adopt the convention where $p$ stands for either $-1$ or $0$ and $q = -(p+1)$.
	In other words, if $p$ stands for one perversity, $q$ will stand for the other.
\end{notation}

Before moving on, three lemmas. 
\begin{lem}\label{cohomology}
	For all $\sh{T} \in \perv{p}{\acat{T}}$ we have $H^i(Y,T) = H^i(X,f_*T)$, for all $i$.
	For all $\sh{F} \in \perv{p}{\acat{F}}$ we have $H^i(Y,F) = H^{i-1}(X,R^1f_*F)$.
\end{lem}
\begin{prf}
	This follows easily from the Leray spectral sequence $H^i(X,R^jf_*E) \Rightarrow H^{i+j}(Y,E)$ once we recall that $R^{>1}f_*$ of any sheaf vanishes.
\end{prf}
\begin{lem}\label{previous-lemma}
	Let $\varphi\colon A \to B$ be a morphism in $\perv{p}{\acat{A}}$, with $A \in \Coh(Y)$ a coherent sheaf.
	Denote by $\perv{p}{\coker}\varphi$ the cokernel of $\varphi$ in the abelian category $\perv{p}{\acat{A}}$.
	The following are equivalent:
	\begin{enumerate}
		\item the perverse cokernel $\perv{p}{\coker}\varphi$ lies in $\perv{p}{\acat{F}[1]}$;
		\item the cone of $\varphi$ is a complex concentrated in degrees $\leq -1$ -- i.e.~it belongs to $\text{D}^{\leq -1}(Y)$;
		\item the induced morphism $H^0(A) \to H^0(B)$ of sheaves is surjective.
	\end{enumerate} 
\end{lem}
\begin{prf}
	Let $C$ be the cone of $\varphi$.
	Perhaps the easiest way to prove the lemma is by using the spectral sequence $H^i\left(H^j_{\perv{p}{\acat{A}}}(C)\right) \Rightarrow H^{i+j}(C)$, where by $H^i_{\perv{p}{\acat{A}}}$ we mean the cohomology object with respect to the heart $\perv{p}{\acat{A}}$.
	It follows that $C \in \perv{p}{\acat{F}[1]}$ if and only if $C \in \text{D}^{\leq -1}(Y)$ if and only if $H^0H^0_{\perv{p}{\acat{A}}}(C) = 0$, which is equivalent to $\varphi$ being surjective in $H^0$.
\end{prf}
\begin{lem}\label{closedunderquotients}
	The subcategory $\perv{p}{\acat{F}_{\leq 1}[1]} \subset \perv{p}{\acat{A}_{\leq 1}}$ is closed under quotients and extensions.
\end{lem}
\begin{prf}
	Closure under extensions is immediate by looking at the cohomology sheaves long exact sequence.
	Let $F[1] \to G$ be an epimorphism in $\perv{p}{\acat{A}}$, with kernel $K$.
	First, $H^0(F[1]) = 0 = H^0(G)$, thus $G = H[1]$, for a sheaf $H$.
	By pushing down to $X$, we have a short exact sequence 
	$$0 \to Rf_*K \to R^1f_*F \to R^1f_*H \to 0$$
	where the middle term is a skyscraper sheaf (this follows from $f_*F=0$ and the fact that $\dim \supp F = 1$).
	In turn, $R^1f_*H$ is a skyscraper sheaf and thus $\dim \supp H = 1$ (see Lemma \ref{evvai}).
\end{prf}

\subsection{Moduli} 
\label{sub:moduli}
To define the motivic Hall algebra of $\perv{p}{\acat{A}}$ in the next section we need, first of all, an algebraic stack $\perv{p}{\astack{A}}$ parameterising objects of $\perv{p}{\acat{A}}$.
We build it as a substack of the stack $\mum_Y$, which was constructed by Lieblich \cite{lieblich} and christened \emph{the mother of all moduli of sheaves}.
For its definition and some further properties we refer the reader to the appendix.
We only recall that $\mum_Y$ parameterises objects in the derived category of $Y$ with no negative self-extensions.
This last condition is key to avoid having to enter the realm of higher stacks.
We point out that as $\perv{p}{\acat{A}}$ is the heart of a t-structure its objects satisfy this condition.

Notice that the definition of $\perv{p}{\acat{A}}$ is independent of the ground field and is stable under field extension.
Concretely, take $\csh{E} \in \mum_Y(T)$ a family of complexes over $Y$ parameterised by a scheme $T$ and $t: \Spec k \to T$ a geometric point.
We can consider $\csh{E}\vert^L_{Y_t}$, the derived restriction of $\csh{E}$ to the fibre $Y_t$ of $Y_T$ over $t$, and it makes sense to write $\csh{E}\vert^L_{Y_t} \in \perv{p}{\acat{A}}$ (where the latter category is interpreted relatively to $k$).
\begin{prop}
	Define a prestack\footnote{We use the term \emph{prestack} in analogy with \emph{presheaf}.} by the rule
	\begin{align*}
		\perv{p}{\astack{A}}(T) = \left\{ \csh{E} \in \stack{Mum}_Y(T) \st \forall t \in T, \csh{E}\vert^L_{Y_t} \in \perv{p}{\acat{A}} \right\}
	\end{align*}
	with restriction maps induced by $\mum_Y$ and where by $t \in T$ we mean that $t: \Spec k \to T$ is a geometric point of $T$.
	The prestack $\perv{p}{\astack{A}}$ is an open substack of $\mum_Y$.
\end{prop}
\begin{prf}
	As mentioned earlier, objects of $\perv{p}{\Per(Y/X)}$ have vanishing self-extensions and therefore can be glued.
	In other words, $\perv{p}{\astack{A}}$ satisfies descent.
	To prove that the inclusion $\perv{p}{\astack{A}} \to \mum_Y$ is open, we employ Van den Bergh's projective generators.
	For this, we introduce some auxiliary spaces.
	
	If $U \subset X$ is open, we can consider the restriction $g\colon V = f^{-1}(U) \to U$ of the morphism $f$.
	The category of perverse coherent sheaves $\perv{p}{\Per(V/U)} =: \perv{p}{\acat{A}_U}$ still makes sense and the corresponding stack $\perv{p}{\astack{A}_U}$ satisfies descent.
	Notice that $\perv{p}{\astack{A}_X} = \perv{p}{\astack{A}}$.
		
	When $U$ is affine, there exists a vector bundle $P$ (a projective generator \cite[Proposition 3.2.5]{vdb}) on $V$ such that an object of the derived category $E$ is perverse coherent (relatively to $g\colon V \to U$) if and only if $\Hom_V(P,E[i]) =0$ for $i \neq 0$.
	In other words, $E$ is perverse coherent if and only if the complex $Rg_* R\lHom(P,E)$ is concentrated in degree zero.
	From this we automatically deduce that the morphism $\perv{p}{\astack{A}_U} \to \mum_V$ is open as this last condition is open.
	
	To pass from local to global, we recall that in \cite{vdb} it was also proved that one check whether a complex is a perverse coherent sheaf on an open cover of $X$.
	In other words, if $E \in D(Y)$ and if $U \to X$ is an open affine cover and $V  = f^{-1}(U)$, then $E \in \perv{p}{\Per(Y/X)}$ if and only if $E\vert V \in \perv{p}{\Per(V/U)}$.
	
	We have restriction morphisms $\mum_Y \to \mum_V$ and $\perv{p}{\astack{A}} \to \astack{A}_U$.
	When $U \to X$ is an open affine cover, we can realise $\perv{p}{\astack{A}}$ as the fibre product of $\perv{p}{\astack{A}_U} \to \mum_U \ot \mum_Y$.
	This is enough to conclude that the morphism $\perv{p}{\astack{A}} \to \mum_Y$ is an open immersion.
\end{prf}

It will be important for us to also have moduli spaces for the torsion and torsion-free subcategories $\perv{p}{\curvy{T}}$, $\perv{p}{\curvy{F}}$.
We define them similarly as above.
\begin{align*}
	\perv{p}{\astack{F}}(T) &= \left\{ \sh{E} \in \stack{A}(T) \st \forall t \in T,  \sh{E}\dvert_{Y_t} \in \perv{p}{\acat{F}} \right\} \\
	\perv{p}{\astack{T}}(T) &= \left\{ \sh{E} \in \stack{A}(T) \st \forall t \in T, \sh{E}\dvert_{Y_t} \in \perv{p}{\acat{T}} \right\}
\end{align*}
Notice that $\perv{p}{\astack{T}} = \perv{p}{\astack{A}} \cap \astack{A}$ and $\perv{p}{\astack{F}[1]} = \perv{p}{\astack{A}} \cap \astack{A}[1]$.
One has the expected open inclusions of algebraic stacks
\begin{align*}
	\perv{p}{\astack{T}}, \perv{p}{\astack{F}} &\subset \stack{A} \subset \mum^{[-1,0]}_Y \\
	\perv{p}{\astack{T}}, \perv{p}{\astack{F}[1]} &\subset \perv{p}{\astack{A}} \subset \mum^{[-1,0]}_Y
\end{align*}
where $\mum^{[-1,0]}_Y$ is the substack of $\mum_Y$ parameterising complexes concentrated in degrees $-1$ and $0$.

We conclude this section with a technical result regarding the structure of $\perv{p}{\astack{A}}$.
This will essentially allow us to carry all the proofs to set up the motivic Hall algebra of $\perv{p}{\acat{A}}$ from the case of coherent sheaves.
\begin{prop}\label{something}
	Let $p=-1$.
	There is a collection of open substacks $\perv{p}{\astack{A}}_n \subset \perv{p}{\astack{A}}$ which jointly cover $\perv{p}{\astack{A}}$.
	Each $\perv{p}{\astack{A}}_n$ is isomorphic to an open substack of $\astack{A}$.
\end{prop}
To prove this result we start by remarking that, as a consequence of our assumptions on $Y$, the structure sheaf $\curlyO_Y$ is a spherical object \cite[Definition 8.1]{huybrechts-fm} in $\text{D}^\text{b}(Y)$.
Thus the Seidel-Thomas spherical twist around it is an autoequivalence of $\text{D}^\text{b}(Y).$
This functor can also be described as the Fourier-Mukai transform with kernel the ideal sheaf of the diagonal of $Y$ shifted by one.
We thus get an exact auto-equivalence $\tau$ of $\text{D}^\text{b}(Y)$ and we notice that the subcategory of complexes with no negative self-extensions is invariant under $\tau$.
As Fourier-Mukai transforms behave well in families \cite[Proposition 6.1]{nahm} we also obtain an automorphism (which by abuse of notation we still denote by $\tau$) of the \emph{stack} $\mum_Y$.

Let us now fix an ample line bundle $\sh{L}$ downstairs on $X.$
Tensoring with $f^*\sh{L}^n$ also induces an automorphism of $\mum_Y$.
The automorphism $\tau_n \in \Aut (\mum_Y)$ is then defined by $\tau_n(E) = \tau(E \otimes f^*\sh{L}^n)$.
The following lemma tells us how to use the automorphisms $\tau_n$ to deduce the proposition above.
\begin{lem}\label{tau}
	Let $p=-1$ and let $\scmplx{E} \in \perv{p}{\curvy{A}}$ be a perverse coherent sheaf.
	Then there exists an $n_0$ such that for all $n\geq n_0$
	$$\widetilde{\tau_n}(\scmplx{E}) = \tau_n(\scmplx{E})[-1] \in \curvy{A}.$$
\end{lem}
\begin{prf}
	The two key properties we use of $\tau_n$ are that it is an exact functor and that for a complex $\scmplx{G}$ we have an exact triangle
	\begin{align*}
		H^\bullet(Y,\scmplx{G}(n))\otimes_\C \curlyO_Y \stackrel{\text{ev}}{\longto}
		\scmplx{G}(n) \to
		\tau_n(\scmplx{G}) \nrightarrow
	\end{align*}
	where $\scmplx{G}(n) = \scmplx{G}\otimes_{\curlyO_Y} f^*\sh{L}^n.$
	
	Let now $\scmplx{E} \in \perv{p}{\curvy{A}}$ be a perverse coherent sheaf together with its torsion pair exact sequence (in $\perv{p}{\acat{A}}$)
	\begin{align*}
		\sh{F}[1] \into \scmplx{E} \onto \sh{T}
	\end{align*}
	where $\sh{F} \in \perv{p}{\acat{F}}$, $\sh{T} \in \perv{p}{\acat{T}}.$
	Using Leray's spectral sequence, the projection formula, Lemma \ref{cohomology} and Serre vanishing on $X$ we can pick $n$ big enough so that all hypercohomologies involved, $H^\bullet(Y,\sh{F}[1](n))$, $H^\bullet(Y,\sh{T}(n))$, $H^\bullet(Y,\csh{E}(n))$, are concentrated in degree zero.
	
	From the triangle
	\begin{align*}
		H^\bullet(Y,\scmplx{E}(n))\otimes_\C \curlyO_Y \stackrel{\text{ev}}{\longto}
		\scmplx{E}(n) \to
		\tau_n(\scmplx{E}) \nrightarrow
	\end{align*}
	we have that $\tau_n(\csh{E}) \in \text{D}^{[-1,0]}(Y)$, similarly for $\tau_n(\sh{F}[1])$ and $\tau_n(\sh{T}).$
	From the triangle
	\begin{align*}
		H^\bullet(Y,\sh{F}[1](n))\otimes_\C \curlyO_Y \stackrel{\text{ev}}{\longto}
		\sh{F}[1](n) \to
		\tau_n(\sh{F}[1]) \nrightarrow
	\end{align*}
	we obtain that $H^0(\tau_n(\sh{F}[1]))=0.$
	
	From the triangle
	\begin{align*}
		\tau_n(\sh{F}[1]) \to
		\tau_n(\csh{E}) \to
		\tau_n(\sh{T}) \nrightarrow
	\end{align*}
	arising from exactness of $\tau_n$ we have that $H^0(\tau_n(\sh{T}))\simeq H^0(\tau_n(\scmplx{E})).$
	Thus to prove the lemma it suffices to show that $H^0(\tau_n(\sh{T}))=0.$
	
	Finally, from the triangle
	\begin{align*}
		H^\bullet(Y,\sh{T}(n))\otimes_\C \curlyO_Y \to
		\sh{T}(n) \to
		\tau_n(\sh{T}) \nrightarrow
	\end{align*}
	one obtains the following exact sequence.
	\begin{align*}
		0 \to H^{-1}(\tau_n(\sh{T})) \to H^0(Y,\sh{T}(n))\otimes_\C \curlyO_Y
		\stackrel{\alpha}{\longto} \sh{T}(n)
		\stackrel{\beta}{\longto} H^0(\tau_n(\sh{T})) \to 0
	\end{align*}
	Thus we have
	$$ \tau_n(\scmplx{E})[-1] \in \acat{A} \iff H^0(\tau_n(\scmplx{E})) \simeq H^0(\tau_n(\sh{T}))=0 \iff \beta = 0. $$
	Let $\sh{K} = \ker \beta$.
	We then have two short exact sequences
	\begin{center}
		\begin{tikzpicture}
			\matrix (m) [matrix of math nodes, row sep=.5em, column sep=3em, text height=1.5ex, text depth=0.25ex]
			{
			H^{-1}(\tau_n(\sh{T})) & H^0(Y,\sh{T}(n))\otimes_\C \curlyO_Y & \sh{K}\\
			\sh{K} & \sh{T}(n) & H^0(\tau_n(\sh{T})) \\
			};
			\path[->,font=\scriptsize]
			(m-1-1) edge[right hook->] (m-1-2)
			(m-1-2) edge[->>] node[auto]{$\gamma$} (m-1-3)
			(m-2-1) edge[right hook->] node[auto]{$\delta$} (m-2-2)
			(m-2-2) edge[->>] node[auto]{$\beta$} (m-2-3)
			;
		\end{tikzpicture}
	\end{center}
	and notice that $\delta \gamma = \alpha.$
	By pushing forward the first sequence via $f_*$ we have that $R^1f_* \sh{K}= 0,$ as $R^1f_*\curlyO_Y=0.$
	Pushing forward the second sequence yields the exact sequence
	\begin{align*}
		f_*\sh{K} \into f_* \sh{T}(n) \onto f_* H^0(\tau_n(\sh{T}))
	\end{align*}
	and $R^1f_*H^0(\tau_n(\sh{T}))=0,$ as $R^1f_*\sh{T}(n)=0$ (this last is a consequence of Lemma \ref{cohomology} and the projection formula).
	
	By taking $n$ even bigger we can assume $f_*\sh{T}(n)$ to be generated by global sections and thus we can assume $f_*\alpha$ to be surjective.
	As $\alpha = \delta \gamma$ we obtain that $f_*\delta$ is surjective and thus $f_*H^0(\tau_n(\sh{T}))=0.$
	As a consequence we have that $H^0(\tau_n(\sh{T})) \in \curvy{C}$.
	
	The sheaf $\sh{T}(n)$ is in $\perv{p}{\acat{T}}$ (this is a simple computation, the key fact to notice is that $\curvy{C}(n) = \curvy{C}$).
	Finally, as $\sh{T}(n) \in \perv{p}{\acat{T}}$ and $H^0(\tau_n(T)) \in \curvy{C}$, $\beta = 0$.
\end{prf}
To prove Proposition \ref{something} we define $\perv{p}{\acat{A}}_n$ to be the subcategory of $\perv{p}{\acat{A}}$ consisting of elements $\csh{E}$ such that $\widetilde{\tau}_n(\csh{E}) \in \acat{A}$.
We can produce a moduli stack for $\perv{p}{\astack{A}_n}$ via the following composition of cartesian diagrams. \begin{center}
	\begin{tikzpicture}
		\matrix (m) [matrix of math nodes, row sep=3em, column sep=3em, text height=1.5ex, text depth=0.25ex]
		{
		\perv{p}{\astack{A}_n} & \widetilde{\tau}_n^{-1}(\astack{A}) & \astack{A} \\
		\perv{p}{\astack{A}} & \mum_X & \mum_X \\
		};
		\path[->,font=\scriptsize]
		(m-1-1) edge [right hook->] (m-1-2)
		(m-1-2) edge (m-1-3)
		(m-1-1) edge [right hook->] (m-2-1)
		(m-2-1) edge [right hook->] (m-2-2)
		(m-2-2) edge node [auto]{$\widetilde{\tau}_n$} (m-2-3)
		(m-1-2) edge (m-2-2)
		(m-1-3) edge [right hook->] (m-2-3)
		;
	\end{tikzpicture}
\end{center}
We obtain that $\perv{p}{\astack{A}}_n$ is an open substack of $\perv{p}{\astack{A}}$ and is isomorphic to an open substack of $\astack{A}$ via $\widetilde{\tau}_n$.
From the previous lemma we have that the sum of the inclusions $\coprod_n \perv{p}{\astack{A}_n} \to \perv{p}{\astack{A}}$ is surjective.
\begin{rmk}\label{crucial perversity remark flops}
	The proof we just presented here of Proposition \ref{something} worked for $p = -1$, and we do not know a direct way to extend this result to the zero perversity.
	However, we can work around this issue by making the following additional assumption (which will hold in the cases which are of interest to us, i.e.~for flops and the McKay correspondence \cite{cala}): we assume the existence of a Fourier-Mukai equivalence taking $\perv{0}{\Per(Y/X)}$ to $\perv{-1}{\Per}(W/X)$, with $W$ a variety over $X$ satisfying the same assumptions as $Y$.
	Using this, we obtain a variant of Proposition \ref{something}, namely for $q = 0$ there exists a collection $\{\perv{q}{\astack{A}_n}\}_n$ of open substacks of $\perv{q}{\astack{A}}$ such that, for every $n$, $\perv{q}{\astack{A}_n}$ is isomorphic to an open substack of the stack of coherent sheaves on $W$.
	We highlight three key places where this is used: Propositions \ref{pexactsequences}, \ref{semiclassicalpminusone} and Theorem \ref{mainmanor}.
	Henceforth, we will tacitly assume this extra hypothesis so that this strategy of passing to $W$ can be applied.
\end{rmk}

\section{Hall Algebras} 
\label{sec:hall_algebras}
This section is devoted to constructing the motivic Hall algebra of perverse coherent sheaves.
We start by recalling the general setup and then move on to check that we can port the construction of the Hall algebra of coherent sheaves to the perverse case.
\begin{rmk*}
	As this porting process relies on Proposition \ref{something}, we remind the reader of Remark \ref{crucial perversity remark flops}.
\end{rmk*}
\subsection{Grothendieck Rings and the Hall Algebra of Coherent Sheaves} 
\label{sub:grothendieck_rings}
In this section we construct the Hall algebra $\text{H}(\perv{p}{\acat{A}})$ of our perverse coherent sheaves, which is a module over $K(\St/\C)$, the Grothendieck ring of stacks over $\C.$
We start by recalling the definition of the latter.
All the omitted proofs can be found, for example, in \cite{joyce-book,tom-hall}.
\begin{defn}
	The \emph{Grothendieck ring of schemes} $K(\Sch/\C)$ is defined to be the $\Q$-vector space spanned by isomorphism classes of schemes of finite type over $\C$ modulo the \emph{cut \& paste} relations:
	\begin{align*}
		[X] = [Y] + [X\setminus Y]
	\end{align*}
	for all $Y$ closed in $X$.
	The ring structure is induced by $[X \times Y] = [X]\cdot[Y].$
\end{defn}
Notice that the zero element is given by the empty scheme and the unit for the multiplication is given by $[\Spec \C]$.
Also, the Grothendieck ring disregards any non-reduced structure, as $[X_\text{red}] = [X]- 0$.
This ring can equivalently be described in terms of geometric bijections and Zariski fibrations.
\begin{defn}
	A morphism $f: X \to Y$ of finite type schemes is a \emph{geometric bijection} if it induces a bijection on $\C$-points $f(\C): X(\C) \to Y(\C)$.

	A morphism $p: X \to Y$ is a \emph{Zariski fibration} if there exists a \emph{trivialising} Zariski open cover of $Y.$
	That is, there exists a Zariski open cover $\{Y_i\}_i$ of $Y$ together with schemes $F_i$ such that $p^{-1}(Y_i) \cong Y_i \times F_i,$ as $Y_i$-schemes.
	
	Two Zariski fibrations $p: X \to Y$, $p': X' \to Y$ \emph{have the same fibres} if there exists a trivialising open cover for both fibrations such that the fibres are isomorphic $F_i \cong F_i'.$
\end{defn}
\begin{lem}
	We can describe the ring $K(\Sch/\C)$ as the $\Q$-vector space spanned by isomorphism classes of schemes of finite type over $\C$ modulo the following relations.\footnote{The three relations we present here are actually redundant, cf.~\cite[Lemma 2.9]{tom-hall}, although the same is not true for stacks.}
	\begin{enumerate}
		\item $[X_1 \amalg X_2] = [X_1]+[X_2]$, for every pair of schemes $X_1$, $X_2$.
		\item $[X_1] = [X_2]$ for every geometric bijection $f: X_1 \to X_2$.
		\item $[X_1]= [X_2]$ for every pair of Zariski fibrations $p_i:X_i \to Y$  with same fibres.
	\end{enumerate}
\end{lem}
We now consider the Grothendieck ring of stacks.
\begin{defn}
	A morphism of finite type algebraic stacks $f: X_1 \to X_2$ is a \emph{geometric bijection} if it induces an equivalence of groupoids on $\C$-points $f(\C): X_1(\C) \to X_2(\C)$.\footnote{We point out that geometric bijections are relative algebraic spaces \cite[Lemma 2.3.9]{abramovich}.}
	
	A morphism of algebraic stacks $p: X \to Y$ is a \emph{Zariski fibration} if given any morphism from a scheme $T \to Y$ the induced map $X \times_Y T \to T$ is a Zariski fibration of schemes.
	In particular a Zariski fibration is a schematic morphism.
	
	Two Zariski fibrations between algebraic stacks $p_i: X_i\to Y$ \emph{have the same fibres} if the two maps $X_i \times_Y T \to T$ induced by a morphism from a scheme $T \to Y$ are two Zariski fibrations with the same fibres.
\end{defn}

\begin{defn}\label{grothendieck_ring_of_stacks}
	The \emph{Grothendieck ring of stacks} $K(\St/\C)$ is defined to be the $\Q$-vector space spanned by isomorphism classes of Artin stacks of finite type over $\C$ with affine geometric stabilisers, modulo the following relations.
	\begin{enumerate}
		\item $[X_1 \amalg X_2] = [X_1] + [X_2]$ for every pair of stacks $X_1, X_2$.
		\item $[X_1] = [X_2]$ for every geometric bijection $f: X_1 \to X_2.$
		\item $[X_1] = [X_2]$ for every pair of Zariski fibrations $p_i : X_i \to Y$ with the same fibres.
	\end{enumerate}
\end{defn}
Let us call $\L = [\A^1]$ the element represented by the affine line.
The obvious ring homomorphism $K(\Sch/\C) \to K(\St/\C)$ becomes an isomorphism after inverting the elements $\L$ and $(\L^k - 1)$, for $k\geq 1$ \cite[Lemma 3.9]{tom-hall}.
Thus the ring homomorphism factors as follows.
\begin{align*}
	K(\Sch/\C) \to K(\Sch/\C)[\L^{-1}] \to K(\St/\C)
\end{align*}
We also mention that through the lens of the Grothendieck ring one cannot tell apart varieties from schemes or even algebraic spaces \cite[Lemma 2.12]{tom-hall}.

It also makes sense to speak of a \emph{relative} Grothendieck group $K(\St/S)$, where $S$ is a fixed base stack which we assume to be Artin, \emph{locally} of finite type over $\C$ and with affine geometric stabilisers.
We define $K(\St/S)$ to be spanned by isomorphism classes of morphisms $[W \to S]$ where $W$ is an Artin stack of finite type over $\C$ with affine geometric stabilisers, modulo the following relations.
\begin{enumerate}
	\item $[f_1\amalg f_2: X_1 \amalg X_2 \longto S] = [X_1 \stackrel{f_1}{\to} S] + [X_2 \stackrel{f_2}{\to} S]$, for every pair of stacks $X_i$.
	\item For a morphism $f: X_1 \to X_2$ over $S$, with $f$ a geometric bijection, $[X_1 \to S] = [X_2 \to S]$.
	\item For every pair of Zariski fibrations with the same fibres $X_1 \to Y \ot X_2$ and every morphism $Y \to S$
	\begin{align*}
		[X_1 \to Y \to S] = [X_2 \to Y \to S].
	\end{align*}
\end{enumerate}
Given a morphism $a:S \to T$ we have a pushforward map
\begin{align*}
	a_*: K(\St/S) &\longto K(\St/T)\\
	[X \to S] &\longmapsto [X \to S \stackrel{a}{\to} T]
\end{align*}
and given a morphism of finite type $b: S \to T$ we have a pullback map 
\begin{align*}
	b^*: K(\St/T) &\longto K(\St/S)\\
	[X \to T] &\longmapsto [X\times_T S \to S].
\end{align*}
The pushforward and pullback are functorial and satisfy base-change.
Furthermore, given a pair of stacks $S_1,S_2$ there is a K\"unneth map 
\begin{align*}
	\kappa: K(\St/S_1)\otimes K(\St/S_2) &\longto K(\St/S_1 \times S_2) \\
	[X_1 \to S_1] \otimes [X_2 \to S_2] &\longmapsto [X_1 \times X_2 \to S_1 \times S_2].
\end{align*}

Take now $\astack{A}$ to be the stack of coherent sheaves on $X,$ where $X$ is smooth and projective over $\C$, and denote by $\text{H}(\acat{A})$ the Grothendieck ring $K(\St/\astack{A})$ (where $\acat{A}$ stands for $\Coh X$).
We can endow $\text{H}(\acat{A})$ with a \emph{convolution product}, coming from the abelian structure of $\acat{A}$.
The product is defined as follows.
Let $\astack{A}^{(2)}$ be the stack of exact sequences in $\acat{A}.$
There are three natural morphisms $a_1,b,a_2:\astack{A}^{(2)} \to \astack{A}$ which take an exact sequence
\begin{align*}
	\sh{A}_1 \into \sh{B} \onto \sh{A}_2
\end{align*}
to $\sh{A}_1,\sh{B},\sh{A}_2$ respectively.
Consider the following diagram.
\begin{center}
	\begin{tikzpicture}
		\matrix (m) [matrix of math nodes, row sep=3em, column sep=3em, text height=1.5ex, text depth=0.25ex]
		{
		\astack{A}^{(2)} & \astack{A}\\
		\astack{A} \times \astack{A} & \\
		};
		\path[->,font=\scriptsize]
		(m-1-1) edge node[auto]{$b$} (m-1-2)
		(m-1-1) edge node[auto,swap] {$(a_1,a_2)$}(m-2-1)
		;
	\end{tikzpicture}
\end{center}
We remark that $(a_1,a_2)$ is of finite type \cite[Lemma 4.2]{tom-hall}.
A \emph{convolution product} can be then defined as follows:
\begin{align*}
	m: \text{H}(\acat{A}) \otimes \text{H}(\acat{A}) \longto \text{H}(\acat{A})\\
	m= b_* (a_1,a_2)^* \kappa.
\end{align*}
Explicitly, given two elements $[X_1 \stackrel{f_1}{\to} \astack{A}]$, $[X_2 \stackrel{f_2}{\to} \astack{A}]$ we write $f_1 * f_2 = m(f_1 \otimes f_2)$ for their product which is given by the top row of the following diagram.
\begin{center}
	\begin{tikzpicture}
		\matrix (m) [matrix of math nodes, row sep=3em, column sep=3em, text height=1.5ex, text depth=0.25ex]
		{
		Z&\astack{A}^{(2)} & \astack{A}\\
		X_1 \times X_2 & \astack{A}\times \astack{A} &\\
		};
		\path[->,font=\scriptsize]
		(m-1-1) edge[color=white] node[color=black]{$\square$} (m-2-2)
		(m-1-1) edge (m-2-1)
				edge (m-1-2)
		(m-1-2) edge node[auto]{$b$} (m-1-3)
				edge node[auto]{$(a_1,a_2)$} (m-2-2)
		(m-2-1) edge node[auto]{$f_1 \times f_2$} (m-2-2)
		(m-1-1) edge[bend left=30,densely dashed] node[auto]{$f_1 * f_2$} (m-1-3)
		;
	\end{tikzpicture}
\end{center}
The convolution product endows $\text{H}(\acat{A})$ with an associative $K(\St/\C)$-algebra structure with unit element given by $[\Spec \C = \astack{A}_0 \subset \astack{A}]$, the inclusion of the zero object.

\subsection{The Hall Algebra of Perverse Coherent Sheaves} 
\label{sub:hall_algebra_of_perverse_coherent_sheaves}
We now assume to be working in Situation \ref{situation}.
We want to replace $\acat{A}$ by $\perv{p}{\acat{A}}$ and construct the analogous algebra $\text{H}(\perv{p}{\acat{A}})$.
We first need the moduli stack $\perv{p}{\astack{A}}^{(2)}$, which parameterises short exact sequences in $\perv{p}{\acat{A}}$.
Define a prestack $\perv{p}{\astack{A}}^{(2)}$ as follows.
To each scheme $T$ we assign a groupoid $\perv{p}{\astack{A}}^{(2)}(T)$, whose objects are exact triangles
\begin{align*}
	\csh{E}_1 \to \csh{E} \to \csh{E}_2 \nrightarrow
\end{align*}
with vertices belonging to $\perv{p}{\astack{A}}(T)$ and whose morphisms are isomorphisms of triangles.
The restriction functors are exact as they are given by derived pullback.
\begin{prop}\label{pexactsequences}
The prestack $\perv{p}{\astack{A}}^{(2)}$ is an Artin stack locally of finite type over $\C$ with affine stabilisers.
\end{prop}
\begin{prf}
	This prestack is well-defined and satisfies descent.
	In fact, given the existence of the stack of objects of $\perv{p}{\acat{A}}$, the only issue arises in gluing automorphisms.
	This is taken care of by noticing that $\Ext^{<0}_{\perv{}{\acat{A}}}(\csh{A},\csh{B})$ vanishes for any two objects $\csh{A},\csh{B} \in \perv{p}{\acat{A}}$ or, in other words, by appealing to \cite[Lemma 2.1.10]{apst}.
	Take now $p=-1$.
	We want to use the functors $\widetilde{\tau}_n$ of Lemma \ref{tau}.
	Notice that the subcategory $\perv{p}{\acat{A}}_n \subset \perv{p}{\acat{A}}$, of objects which become coherent after a twist by $\widetilde{\tau}_n$, is extension-closed.
	Hence, we have a well-defined stack of exact sequences $\perv{p}{\astack{A}_n^{(2)}}$.
	As $\perv{p}{\astack{A}}_n$ is an open substack of $\perv{p}{\astack{A}}$, by Proposition \ref{something}, we deduce that $\perv{p}{\astack{A}_n}^{(2)}$ an open substack of $\perv{p}{\astack{A}}^{(2)}$.
	Using once more Proposition \ref{something} and the fact that $\widetilde{\tau}_n$ is an exact functor we can embed $\perv{p}{\astack{A}_n^{(2)}}$ inside $\astack{A}^{(2)}$, thus proving that $\perv{p}{\astack{A}_n^{(2)}}$ is algebraic.
	
	The sum $\coprod_n \perv{p}{\astack{A}_n^{(2)}} \to \perv{p}{\astack{A}}^{(2)}$ is surjective, by Lemma \ref{tau}, and thus the stack $\perv{p}{\astack{A}^{(2)}}$ is algebraic.
	All other properties are deduced by the fact that $\perv{p}{\astack{A}_n^{(2)}}$ is an open substack of $\astack{A}^{(2)}$.
	For $p= 0$, one appeals to Remark \ref{crucial perversity remark flops}.
\end{prf}
\noindent The proof actually produces more: it gives an analogue of Proposition \ref{something}.

As for coherent sheaves, the stack $\perv{p}{\astack{A}}^{(2)}$ comes equipped with three morphisms $a_1,b,a_2$, sending a triangle of perverse coherent sheaves
\begin{align*}
	\csh{E}_1 \to \csh{E} \to \csh{E}_2 \nrightarrow
\end{align*}
to $\csh{E}_1, \csh{E}, \csh{E}_2$ respectively.
The exact functor $\widetilde{\tau}_n$ yields a commutative diagram
\begin{center}
	\begin{tikzpicture}
		\matrix (m) [matrix of math nodes, row sep=3em, column sep=3em, text height=1.5ex, text depth=0.25ex]
		{
		\perv{p}{\astack{A}_n^{(2)}} & \astack{A}^{(2)}\\
		\perv{p}{\astack{A}_n} \times \perv{p}{\astack{A}_n} & \astack{A} \times \astack{A}\\
		};
		\path[->,font=\scriptsize]
		(m-1-1) edge node [auto,swap]{$(a_1,a_2)$} (m-2-1)
		(m-1-1) edge[right hook->] (m-1-2)
		(m-1-2) edge (m-2-2)
		(m-2-1) edge [right hook->] (m-2-2)
		;
	\end{tikzpicture}
\end{center}
where the vertical arrow on the right is the corresponding morphism for coherent sheaves, which is of finite type.
From this last observation and the fact that being of finite type is local on the target, we automatically have that the (global) morphism $(a_1,a_2): \perv{p}{\astack{A}}^{(2)} \to \perv{p}{\astack{A}}^2$ is of finite type.
To define the convolution product on $K(\St/\perv{p}{\astack{A}})$ (or equivalently the algebra structure of $\text{H}(\perv{p}{\acat{A}})$) we may proceed analogously as for coherent sheaves.
As usual, this discussion is valid for $p= -1$, but an entirely parallel one can be carried out for $p= 0$ using Remark \ref{crucial perversity remark flops}.

\subsection{More Structure on Hall Algebras} 
\label{sub:more_structure_on_hall_algebras}
There is a natural way to bestow a grading upon our Hall algebras.
Recall that for a triangulated category $\acat{T}$ and the heart $\acat{H}$ of a bounded t-structure on $\acat{T}$, the Grothendieck groups $K(\acat{T})$ and $K(\acat{H})$ coincide (by taking alternating sums of cohomology objects).
In particular, $K(\text{D}^\text{b}(Y))$ can be viewed as both $K(\acat{A})$ or $K(\perv{p}{\acat{A}})$.
The \emph{Euler form} $\chi$ is defined as
\begin{align*}
	\chi(\sh{E},\sh{F}) = \sum_j (-1)^{j}\dim_\C \Ext_Y^j(\sh{E},\sh{F})
\end{align*}
on coherent sheaves $\sh{E}, \sh{F}$ and then extended to the whole of $K(Y)$.
By Serre duality the left and right radicals of $\chi$ are equal and we define the \emph{numerical} Grothendieck group of $Y$ as $N(Y) = K(Y)/K(Y)^\perp$.
As the numerical class of a complex stays constant in families, we have a decomposition
\begin{align*}
	\mum_Y = \coprod_{\alpha \in N(Y)} \mum_{Y,\alpha}
\end{align*}
where $\mum_{Y,\alpha}$ parameterises complexes of class $\alpha$.
Let $\Gamma$ denote the \emph{positive cone} of coherent sheaves, i.e.~the image of objects of $\acat{A}$ inside $N(Y)$.
It is a submonoid of $N(Y)$ and for $\astack{A}$ the previous decomposition can be refined to
\begin{align*}
	\astack{A} = \coprod_{\alpha \in \Gamma} \astack{A}_{\alpha}.
\end{align*}
We can also define sub-modules $\text{H}(\acat{A})_\alpha \subset \text{H}(\acat{A})$, where $\text{H}(\acat{A})_\alpha$ denotes $K(\St/\astack{A}_\alpha)$ (which can be thought as spanned by classes of morphisms $[W \to \astack{A}]$ factoring through $\astack{A}_\alpha$).
We then get a $\Gamma$-grading
\begin{align*}
	\text{H}(\acat{A}) = \bigoplus_{\alpha \in \Gamma} \text{H}(\acat{A})_\alpha.
\end{align*}
Analogously, we have a positive cone $\perv{p}{\Gamma} \subset N(X)$ of perverse coherent sheaves.
The Hall algebra thus decomposes as follows
\begin{align*}
	\text{H}(\perv{p}{\acat{A}}) = \bigoplus_{\alpha \in \perv{p}{\Gamma}} \text{H}(\perv{p}{\acat{A}})_\alpha.
\end{align*}

We mentioned earlier that the morphism from the Grothendieck ring of varieties to the Grothendieck ring of stacks factors as follows
\begin{align*}
	K(\Sch/\C) \to K(\Sch/\C)[\L^{-1}] \to K(\St/\C).
\end{align*}
Let $R = K(\Sch/\C)[\L^{-1}]$.
One can define a subalgebra \cite[Theorem 5.1]{tom-hall} $\text{H}_\text{reg}(\acat{A})$ of \emph{regular elements} as the $R$-module spanned by classes $[W \to \astack{A}]$ with $W$ a scheme.
We have an analogous setup for perverse coherent sheaves.
\begin{prop}\label{semiclassicalpminusone}
	Let $\text{H}_\text{reg}(\perv{p}{\acat{A}})$ to be the sub-$R$-module spanned by classes $[W \to \perv{p}{\astack{A}}]$ with $W$ a scheme.
	Then $\text{H}_\text{reg}(\perv{p}{\acat{A}})$ is closed under the convolution product and the quotient
	\begin{align*}
		\text{H}_\text{sc}(\perv{p}{\acat{A}}) = \text{H}_\text{reg}(\perv{p}{\acat{A}})/(\L -1)\text{H}_\text{reg}(\perv{p}{\acat{A}})
	\end{align*}
	is a commutative $K(\Sch/\C)$-algebra.
\end{prop}
\begin{prf}
	Once again, we may appeal to the case of coherent sheaves by using the functors $\widetilde{\tau}_n$.
	Let $p = -1$.
	Let $[f_1:S_1 \to \perv{p}{\astack{A}}]$, $[f_2:S_2 \to \perv{p}{\astack{A}}]$ be two elements of $\text{H}(\perv{p}{\acat{A}})$ such that the $S_i$ are schemes.
	Consider the two morphisms
	\begin{align*}
		f_1\times f_2:& S_1 \times S_2 \to \perv{p}{\astack{A}} \times \perv{p}{\astack{A}}\\
		(a_1,a_2):& \perv{p}{\astack{A}^{(2)}} \to \perv{p}{\astack{A}} \times \perv{p}{\astack{A}}
	\end{align*}
	used to define the product $f_1 * f_2$ in $\text{H}(\perv{p}{\acat{A}})$.
	It suffices to show that the fibre product
	\begin{align*}
		T = (S_1 \times S_2) \times_{\perv{p}{\astack{A}} \times \perv{p}{\astack{A}}} \perv{p}{\astack{A}^{(2)}}
	\end{align*}
	is a regular element.
	Consider the open cover $\{\perv{p}{\astack{A}_n}\}_n$  of $\perv{p}{\astack{A}}$ given in Proposition \ref{something}.
	The first thing we notice is that the collection $\{\perv{p}{\astack{A}_n} \times \perv{p}{\astack{A}_n}\}_n$ is an open cover of $\perv{p}{\astack{A}}\times \perv{p}{\astack{A}}$ (it covers the whole product via Lemma \ref{tau}).
	Pulling it back via $f_1 \times f_2$ yields open covers $\{S_{i,n}\}_n$ for each of the $S_i$ and an open cover $\{S_{1,n}\times S_{2,n}\}_n$ of $S_1 \times S_2$.
	
	On the other hand, by the proof of Proposition \ref{pexactsequences} we have an open cover $\{\perv{p}{\astack{A}^{(2)}_n}\}_n$ of $\perv{p}{\astack{A}^{(2)}}$.
	By pulling back we obtain an open cover $\{T_n\}_n$ of $T$.
	By chasing around base-changes one can see that
	\begin{align*}
		T_n = (S_{1,n} \times S_{2,n}) \times_{\perv{p}{\astack{A}_n} \times \perv{p}{\astack{A}_n}} \perv{p}{\astack{A}_n^{(2)}}.
	\end{align*}
	The functor $\widetilde{\tau}_n$ induces morphisms $\perv{p}{\astack{A}_n} \times \perv{p}{\astack{A}_n} \to \astack{A} \times \astack{A}$, $\perv{p}{\astack{A}_n^{(2)}} \to \astack{A}^{(2)}$ and it is easy to check that
	\begin{align*}
		\perv{p}{\astack{A}_n^{(2)}} = (\perv{p}{\astack{A}_n} \times \perv{p}{\astack{A}_n}) \times_{\astack{A}\times \astack{A}} \astack{A}^{(2)}
	\end{align*}
	thus $T_n = (S_{1,n} \times S_{2,n}) \times_{\astack{A} \times \astack{A}} \astack{A}^{(2)}$ and by \cite[Theorem 5.1]{tom-hall} it is a regular element.
	We conclude that $T$ is also a regular element.
	
	When $p= 0$ one may use Remark \ref{crucial perversity remark flops}.
\end{prf}

We now briefly turn back to the case of coherent sheaves.
The \emph{semi-classical} Hall algebra of coherent sheaves $\text{H}_\text{sc}(\acat{A})$, defined as $\text{H}_\text{reg}(\acat{A})/(\mathbb{L}-1)\text{H}_\text{reg}(\acat{A})$, can be equipped with a Poisson bracket given by
\begin{align*}
	\{f,g\} = \frac{f*g - g* f}{\mathbb{L}- 1}.
\end{align*}
There is another Poisson algebra $\Q_\sigma[\Gamma]$, which depends on a choice $\sigma \in \{-1,1\}$,
defined as the $\Q$-vector space spanned by symbols $q^\alpha$, with $\alpha \in \perv{p}{\Gamma}$, together with a product
\begin{align*}
	q^{\alpha_1} * q^{\alpha_2} = \sigma^{\chi(\alpha_1,\alpha_2)} q^{\alpha_1 + \alpha_2}.
\end{align*}
and a Poisson bracket
\begin{align*}
	\{ q^{\alpha_1}, q^{\alpha_2}\} = \sigma^{\chi(\alpha_1,\alpha_2)} \chi(\alpha_1,\alpha_2)q^{\alpha_1 + \alpha_2} =
	\chi(\alpha_1,\alpha_2)(q^{\alpha_1}*q^{\alpha_2}).
\end{align*}
Given a locally constructible function \cite[Chapter 2]{joyce-book} $\lambda: \astack{A}(\C) \to \Z$, there exists a so-called (at least when $\lambda$ satisfies some properties) \emph{integration morphism}
\begin{align*}
	I\colon \text{H}_\text{sc}(\acat{A})\to \Q_\sigma[\Gamma].
\end{align*}
For convenience of the reader we compactly recall its properties  \cite[Theorem 5.2]{tom-hall}.
The map $I$ is the unique homomorphism of rational vector spaces such that if $V$ is a variety and $f\colon V \to \astack{A}$ factors through $\astack{A}_\alpha$, for $\alpha \in \Gamma$, then 
\begin{align*}
	I \left( [f] \right) = \chi\left( V, f^*\lambda \right) q^\alpha
\end{align*}
where
\begin{align*}
	\chi_\text{top}(V,f^*\lambda) = \sum_{n \in \Z}n\chi_\text{top}((\lambda \circ f)^{-1}(n))
\end{align*}
and where, for a variety $V$, $\chi_\text{top}(V)$ denotes the topological Euler characteristic.
Moreover, $I$ is a homomorphism of commutative algebras if, for all $F,G \in \acat(A)$,
\begin{align*}
	\lambda(F \oplus G) = \sigma^{\chi(F,G)}\lambda(F)\lambda(G)
\end{align*}
and is a homomorphism of Poisson algebras if the expression
\begin{align*}
	\chi\left( \PP\Ext_{\acat{A}}^1(F,G), \lambda({E_\theta}) - \lambda({E_0}) \right)
\end{align*}
is symmetric in $F$ and $G$.
The notation $E_\theta$ stands for the extension
\begin{align*}
	0 \to G \to E_\theta \to F \to 0
\end{align*}
corresponding to a class $\theta \in \Ext_\acat{A}^1(F,G)$.

For $\sigma = 1$ one can choose $\lambda$ to be identically equal to $1$. 
This gives a well-defined integration morphism which in turn leads to \emph{naive} curve counting invariants.
We are more interested in the case $\sigma = -1$ (although what follows certainly holds for the naive invariants as well) where one takes Behrend's microlocal function $\nu$.
For $\text{H}_\text{sc}(\acat{A})$ we know \cite[Theorem 5.5]{joyce-book} that the Behrend function satisfies the necessary hypotheses and thus yields an integration morphism.

To define an integration morphism in the context of perverse coherent sheaves we first define $\Q_\sigma[\perv{p}{\Gamma}]$ analogously as $\Q_\sigma[\Gamma]$, but using the cone of perverse coherent sheaves.
In this context, we may still use Behrend's function.
More precisely, every Artin stack $\astack{M}$ locally of finite type over $\C$ comes equipped with a Behrend function $\nu_{\astack{M}}$ and given any smooth morphism $f: \astack{M}' \to \astack{M}$ of relative dimension $d$ we have $f^*\nu_{\astack{M}} = (-1)^d \nu_{\astack{M}'}$.
To obtain an integration morphism on $\text{H}(\perv{p}{\acat{A}})$ the Behrend function must satisfy the assumptions of \cite[Theorem 5.2]{tom-hall}.
As the latter can be checked on neighbourhoods of points of $\perv{p}{\astack{A}}$, and we know that $\perv{p}{\astack{A}}$ is locally isomorphic to $\astack{A}$, the assumptions are satisfied and we have a well-defined integration morphism
\begin{align*}
	I : \text{H}(\perv{p}{\acat{A}}) \to \Q_\sigma[\perv{p}{\Gamma}].
\end{align*}

	\section{Identities} 
\label{sec:identities}
As hinted at in the introduction, the proof of our main result can be roughly divided into two blocks: the first is concerned with proving a formula relating `perverse' DT invariants with ordinary ones; the second uses this formula to compare the DT invariants over a flop.
We will start by focusing on the former.

Recall that we denote by $\acat{A}$ the category of coherent sheaves of $Y$.
In the previous sections we reminded ourselves of the category of perverse coherent sheaves $\perv{p}{\acat{A}}$ and of the subcategories $\perv{p}{\acat{T}}, \perv{p}{\acat{F}}$.
We also reminded ourselves of the motivic Hall algebra of coherent sheaves $\text{H}(\acat{A})$, defined as the Grothendieck ring $K(\St/\astack{A})$ of stacks over the stack of coherent sheaves $\astack{A}$ equipped with the convolution product.
We also constructed a moduli stack $\perv{p}{\astack{A}}$ parameterising objects in $\perv{p}{\acat{A}}$ and the Hall algebra $\text{H}(\perv{p}{\acat{A}})$ of perverse coherent sheaves, together with the subalgebra of regular elements $\text{H}_\text{reg}(\perv{p}{\acat{A}})$, its semi-classical limit $\text{H}_\text{sc}(\perv{p}{\acat{A}})$ and the integration morphism $I: \text{H}_\text{sc}(\perv{p}{\acat{A}}) \to \Q_\sigma[\perv{p}{\Gamma}]$.
Recall that $\perv{p}{\Gamma}$ is the cone of perverse coherent sheaves sitting inside the numerical Grothendieck group $N(Y)$ and we take $\sigma = -1,1$ depending on the choice of a locally constructible function on $\perv{p}{\astack{A}}$ (either the function identically equal to one or the Behrend function).
\subsection{A Route} 
\label{sub:a_route}
Before we start off, we would like to give a moral proof our main result, which will later guide us through the maze of technical details.
As we are interested in counting \emph{curves}, we will restrict to sheaves (and complexes) supported in dimension at most one.
All the constructions and definitions so far restrict to this setting, and we will append a $\leq 1$ subscript to notify this change (e.g.~we deal with the Hall algebra $H(\perv{p}{\acat{A}}_{\leq 1})$ of perverse coherent sheaves supported in dimension at most one).
The two key results are the identities \eqref{identity-respect}, \eqref{identity-two}.
Continuing from the introduction, our goal is to understand the relationship between perverse DT numbers $\perv{p}{\DT}(Y/X)$ and ordinary DT numbers $\DT(Y)$.
For simplicity, let us for the moment assume $X$ to have a singular locus of dimension zero. 

The Hilbert scheme of curves and points $\Hilb_{\leq 1}(Y)$ maps to $\astack{A}_{\leq 1}$ by taking a quotient $\O_Y \onto \sh{E}$ to $\sh{E}$, thus defining an element $\curly{H}_{\leq 1} \in \text{H}(\acat{A}_{\leq 1})$.\footnote{Strictly speaking this is false as $\Hilb(Y)$ is not of finite type. We shall later enlarge our Hall algebra precisely to deal with this issue.}
From the previous section we know that the integration morphism is related to taking weighted Euler characteristics and in fact integrating $\curly{H}_{\leq 1}$ gives the generating series for the DT invariants\footnote{Again, this is slightly imprecise, there is a sign issue to be explained at the beginning of Subsection \eqref{sub:conclusion}.}
\begin{align*}
	I(\curly{H}_{\leq 1}) \text{``=''} \DT(Y) := \sum_{\beta,n} \DT_Y(\beta,n)q^{(\beta,n)}
\end{align*}
where $\beta \in N_1(Y)$ ranges among curve-classes in $Y$ and $n \in \Z$ is a zero-cycle.
The perverse Hilbert scheme $\perv{p}{\Hilb}_{\leq 1}(Y/X)$ produces a corresponding element $\perv{p}{\curly{H}}_{\leq 1}$ of $\text{H}(\perv{p}{\acat{A}}_{\leq 1})$, which upon being integrated produces $\perv{p}{\DT{(Y/X)}}$.

The first thing we remark is that, as quotients (in $\acat{A}$) of $\O_Y$ lie in $\perv{p}{\acat{T}}$ and $\perv{p}{\acat{T}} \subset \perv{p}{\acat{A}}$, we can interpret $\curly{H}_{\leq 1}$ as an element of $\text{H}(\perv{p}{\acat{A}}_{\leq 1})$.
There is an element $1_{\perv{p}{\acat{F}_{\leq 1}[1]}}$ in $\text{H}(\perv{p}{\acat{A}}_{\leq 1})$ represented by the inclusion $\perv{p}{\astack{F}_{\leq 1}[1]} \subset \perv{p}{\astack{A}}_{\leq 1}$.
There is also a stack parameterising objects of $\perv{p}{\acat{F}_{\leq 1}[1]}$ together with a morphism from $\O_Y$. 
This stack maps down to $\perv{p}{\astack{A}}_{\leq 1}$ by forgetting the morphism, yielding an element $1_{\perv{p}{\acat{F}_{\leq 1}[1]}}^\O$.
We will prove that there is an identity
\begin{align}\label{identity-respect}
	\perv{p}{\curly{H}}_{\leq 1} * 1_{\perv{p}{\acat{F}_{\leq 1}[1]}} = 1_{\perv{p}{\acat{F}_{\leq 1}[1]}}^\O * \curly{H}_{\leq 1}
\end{align}
in the Hall algebra of perverse coherent sheaves.
Let us see how one might deduce this.

We extend the notation $1_{\perv{p}{\acat{F}[1]}}, 1^\O_{\perv{p}{\acat{F}[1]}}$ to general subcategories $\acat{B} \subset \perv{p}{\acat{A}}$ (whenever we have an open inclusion of stacks $\astack{B} \subset \perv{p}{\astack{A}}$) producing elements $1_\acat{B}, 1^\O_\acat{B}$ in $\text{H}(\perv{p}{\acat{A}})$, and similarly for $\text{H}(\acat{A})$.
As $(\perv{p}{\acat{T}},\perv{p}{\acat{F}})$ is a torsion pair in ${\acat{A}}$, we have an identity $1_{{\acat{A}}} = 1_{\perv{p}{\acat{T}}} * 1_{\perv{p}{\acat{F}}}$.
This follows from the fact that for any coherent sheaf $\sh{E}$ there is a unique exact sequence $\sh{T} \into \sh{E} \onto \sh{F}$ with $\sh{T} \in \perv{p}{\acat{T}}$, $\sh{F} \in \perv{p}{\acat{F}}$.
Notice that the product $1_{\perv{p}{\acat{T}}} * 1_{\perv{p}{\acat{F}}}$ is given by $[Z \to {\astack{A}}]$ where $Z$ parameterises exact sequences $\sh{T} \into \sh{E} \onto \sh{F}$ and the morphism $Z \to \astack{A}$ sends such an exact sequence to $\sh{E}$.

We also have an identity $1^\O_{{\acat{A}}} = 1^\O_{\perv{p}{\acat{T}}} * 1^\O_{\perv{p}{\acat{F}}}$.
This is a consequence of the previous identity plus the fact that $\Hom(\O_Y,\perv{p}{\acat{F}}) = 0$ (Lemma \ref{cohomology}).
This last fact also tells us that $1^\O_{\perv{p}{\acat{F}}} = 1_{\perv{p}{\acat{F}}}$.
Moreover, the first isomorphism theorem for the abelian category $\acat{A}$ is reflected in the identity $1^\O_\acat{A} = \curly{H} * 1_\acat{A}$ (any morphism $\O_Y \to \sh{E}$ factors through its image).
Combining everything together (and restricting to sheaves supported in dimension at most one) we see that $\curly{H}_{\leq 1} = 1^\O_{\acat{A}_{\leq 1}} * 1^{-1}_{\acat{A}_{\leq 1}} = 1^\O_{\perv{p}{\acat{T}_{\leq 1}}} * 1^{-1}_{\perv{p}{\acat{T}}_{\leq 1}}$.

A parallel argument can be carried out for $\perv{p}{\acat{A}}$ yielding
\begin{align*}
	\perv{p}{\curly{H}}_{\leq 1} = 1^\O_{\perv{p}{\acat{A}}_{\leq 1}} * 1^{-1}_{\perv{p}{\acat{A}}_{\leq 1}} = 1^\O_{\perv{p}{\acat{F}_{\leq 1}[1]}} * ( 1^\O_\perv{p}{\acat{T}_{\leq 1}} * 1^{-1}_{\perv{p}{\acat{T}_{\leq 1}}}) * 1^{-1}_{\perv{p}{\acat{F}_{\leq 1}[1]}} = 1^\O_{\perv{p}{\acat{F}_{\leq 1}[1]}} * {\curly{H}_{\leq 1}} * 1^{-1}_{\perv{p}{\acat{F}_{\leq 1}[1]}}
\end{align*}
from which we extract \eqref{identity-respect}.
Notice that for the identity $1^\O_{\perv{p}{\acat{A}}_{\leq 1}} = 1^\O_{\perv{p}{\acat{F}_{\leq 1}[1]}} * 1^\O_{\perv{p}{\acat{T}_{\leq 1}}}$ one uses $\Hom_Y(\O_Y,\perv{p}{\acat{F}_{\leq 1}}[2]) = 0$.

We now want to understand how to rewrite $1_{\perv{p}{\acat{F}_{\leq 1}[1]}}^\O$ in a more familiar form.
It turns out that duality almost interchanges $\perv{q}{\acat{T}}$ and $\perv{p}{\acat{F}}$, where $q = -(p+1)$.
Precisely, let $\acat{Q}$ be the subcategory of $\acat{A}$ consisting of sheaves with no subsheaves supported in dimension zero.
Let $\acat{Q}_{\exc}$ denote the subcategory of $\acat{Q}$ made up of sheaves $Q$ such that $\R f_* Q$ is supported in dimension zero and let $\perv{q}{\acat{T}}_\bullet = \acat{Q}_{\exc} \cap \perv{q}{\acat{T}}$.
It is a simple computation (Lemma \ref{duality lemma}) to check that the duality functor $\D = R\lHom_Y(-,\O_Y)[2]$ takes $\perv{q}{\acat{T}}_\bullet$ to $\perv{p}{\acat{F}_{\leq 1}}$.
The category $\acat{Q}$ is related to DT invariants in the following way.

There is an identity $1_\acat{Q}^\O = \curly{H}^\# * 1_\acat{Q}$ in $\text{H}(\acat{A})$, where $\curly{H}^\#$ corresponds to (yet another) Hilbert scheme of a tilt $\acat{A}^\#$ of $\acat{A}$, where $\acat{A}^\#$ is the category in which quotients of $\O_Y$ are the so-called \emph{stable pairs} of Pandharipande and Thomas \cite{pandartsp} (see also \cite{tom-cc}).
We can restrict to sheaves with zero-dimensional pushdown, which yields an identity $1_{\acat{Q}_{\exc}}^\O = \curly{H}^\#_{\exc} * 1_{\acat{Q}_{\exc}},$
which can be refined to	$1_{\perv{q}{\acat{T}_\bullet}}^\O = \curly{H}_{\exc}^\# * 1_{\perv{q}{\acat{T}_\bullet}}.$
Integrating $\curly{H}_{\exc}^\#$ gives the generating series for the Pandharipande-Thomas (PT) invariants of $Y$ \cite[Lemma 5.5]{tom-cc}
\begin{align*}
	I(\curly{H}_{\exc}^\#) \text{``=''} \PT_{\exc}(Y) := \sum_{\substack{\beta,n \\ f_*\beta = 0}} \PT_Y(\beta,n)q^{(\beta,n)}
\end{align*}
where $\beta$ ranges over the curve-classes contracted by $f$.
If we let
\begin{align*}
	\DT_0(Y) := \sum_n \DT_Y(0,n)q^n
\end{align*}
we know \cite[Theorem 1.1]{tom-hall} that the \emph{reduced} DT invariants $\DT'(Y) :=  \DT(Y)/\DT_0(Y)$ coincide with the PT invariants $\PT(Y)$.

Now, the (shifted) duality functor $\D' = \D[1]$ induces a anti-homomorphism between Hall algebras\footnote{More precisely it induces a morphism between certain subalgebras to be defined below.} and takes $\perv{q}{\acat{T}_\bullet}$ to $\perv{p}{\acat{F}_{\leq 1}[1]}$, so we have $\D'(1_{\perv{q}{\acat{T}_\bullet}}) = 1_{\perv{p}{\acat{F}_{\leq 1}[1]}}$.
Furthermore, as a consequence of Serre duality, $\D'(1_{\perv{q}{\acat{T}_\bullet}}^\O) = 1_{\perv{p}{\acat{F}_{\leq 1}[1]}}^\O$.
As a result we have
\begin{align}\label{identity-two}
	1_{\perv{p}{\acat{F}_{\leq 1}[1]}}^\O = 1_{\perv{p}{\acat{F}_{\leq 1}[1]}} * \D'\left(\curly{H}_{\exc}^\# \right)
\end{align}
as	$
	1_{\perv{p}{\acat{F}_{\leq 1}[1]}}^\O =
	\D'\left(1_{\perv{q}{\acat{T}_\bullet}}^\O\right) = 
	\D'\left(\curly{H}_{\exc}^\# * 1_{\perv{q}{\acat{T}_\bullet}} \right) = 
	\D'\left( 1_{\perv{q}{\acat{T}_\bullet}} \right) * \D'\left(\curly{H}_{\exc}^\# \right) = 
	1_{\perv{p}{\acat{F}_{\leq 1}[1]}} * \D'\left(\curly{H}_{\exc}^\# \right)
	$
(notice that duality is an \emph{anti}-equivalence and thus swaps extensions).
We can rewrite \eqref{identity-respect} as follows.
\begin{align}\label{identity-rewritten}
	\perv{p}{\curly{H}_{\leq 1}} * 1_{\perv{p}{\acat{F}_{\leq 1}[1]}} = 
	1_{\perv{p}{\acat{F}_{\leq 1}[1]}} * \D'\left( \curly{H}_{\exc}^\# \right) * \curly{H}_{\leq 1}
\end{align}
Duality and integration can be interchanged up to a flip in signs.
Precisely
\begin{align*}
	I\left( \D'\left( \curly{H}_{\exc}^\# \right) \right)
	\text{``=''} \PT_{\exc}^\vee(Y)
	:= \sum_{\substack{\beta,n \\ f_* \beta = 0}} \PT_Y(-\beta,n) q^{(\beta,n)}.
\end{align*}
Upon integrating the two sides of \eqref{identity-rewritten} the two $1_{\perv{p}{\acat{F}_{\leq 1}[1]}}$ cancel out\footnote{This is the content of Proposition \ref{adexp}, a consequence of an important result of Joyce.} and we are left with the identity 
\begin{align*} 
	\perv{p}{\DT{(Y/X)}} = \PT_{\exc}^\vee(Y) \cdot \DT(Y).
\end{align*}
\subsection{The Perverse Hilbert Scheme} 
\label{sub:first_subsection}
We now proceed along the route traced in the previous subsection, but taking care of technical details.
We revert to Situation \ref{situation}.
Let us start by working in \emph{infinite-type} versions $\text{H}_\infty(\acat{A}_{\schnitzel})$, $\text{H}_\infty(\perv{p}{\acat{A}_{\schnitzel}})$ of our Hall algebras.
The advantage of $\text{H}_\infty$ is that we include stacks \emph{locally} of finite type over $\C$ (e.g.~$\perv{p}{\astack{A}}_{\schnitzel}$), the disadvantage is that we do not have an integration morphism at our disposal.
To define this algebra we proceed exactly as in the previous section: the only differences being that we allow our stacks to be \emph{locally} of finite type over $\C$, we insist that geometric bijections be finite type morphisms and we disregard the disjoint union relation.\footnote{If we allowed both the disjoint union relation and spaces of infinite type then we would be left with the zero ring.
Indeed the standard trick would apply: by removing a point from an infinite disjoint union of points we would conclude that one is equal to zero.
The finite type assumption for geometric bijections is there to avoid pathologies such as an infinite disjoint union of points representing the same class as a line.
}

\begin{rmk}
	It seems worthwhile to point out the following.
	Of course, here we make the (often unjustified) assumption that not only the reader has survived this far, but that he/she is also paying attention to all the details.
	Rather than the infinite-type Hall algebra we just defined, we should really be working in the \emph{Laurent} Hall algebra $\text{H}_\Lambda$ of Section \ref{sub:laurent_elements}.
	We decided not to burden the reader with yet another definition and to temporarily work with $\text{H}_\infty$ instead.
	The identities we prove, starting from \eqref{hilb-f=fo-philb} and eventually leading up to \eqref{identity-leq1}, make sense and are true (with identical proofs) in the Laurent Hall algebra.
\end{rmk}

The first element we consider is $\curly{H}_{\schnitzel} \in \text{H}_\infty(\acat{A}_{\schnitzel})$ corresponding to the Hilbert scheme of $Y$, which parameterises quotients of $\curlyO_Y$ in $\acat{A}_{\schnitzel}$.
To be precise, $\curly{H}_{\schnitzel}$ is represented by the forgetful morphism $\Hilb_{\schnitzel}(Y) \to \astack{A}_{\schnitzel}$, which takes a quotient $\curlyO_Y \onto \sh{E}$ to $\sh{E}$.
For us, the important thing to notice is that if $\curlyO_Y \onto \sh{E}$ is a quotient in $\acat{A}_{\schnitzel}$, then $\sh{E} \in \perv{p}{\acat{T}}$.
This is a consequence of $\curlyO_Y \in \perv{p}{\acat{T}}$ and of the fact that the torsion part of a torsion pair is closed under quotients.
Thus the morphism $\Hilb_{\schnitzel}(Y) \to \astack{A}_{\schnitzel}$ factors through $\perv{p}{\astack{T}_{\schnitzel}}$.
As $\perv{p}{\astack{T}_{\schnitzel}} \subset \perv{p}{\astack{A}_{\schnitzel}}$, $\curly{H}_{\schnitzel}$ can be interpreted as an element of $\text{H}_\infty(\perv{p}{\acat{A}}_{\schnitzel})$.

\begin{rmk}\label{addedremark}
	For the reader interested in the general setup of Situation \ref{situation}, we must replace $\Hilb_{\leq 1}(Y)$ with an open subset $\Hilb_{\leq 1}^\partial(Y) \subset \Hilb_{\leq 1}(Y)$ and thus replace $\curly{H}_{\leq 1}$ in the discussion above with $\curly{K}_{\leq 1}$, the element corresponding to this ``partial'' Hilbert scheme.
	The space $\Hilb_{\leq 1}^\partial(Y)$ parameterises quotients $\varphi \colon\O_Y \to E$, such that the perverse cokernel of $\varphi$ lies in $\perv{p}{\acat{F}_{\leq 1}[1]}$ -- as opposed to $\perv{p}{\acat{F}[1]}$ (cfr. Lemma \ref{previous-lemma}).
	For the reader interested solely in flops, i.e.~when the singular locus of $X$ is zero-dimensional, then $\perv{p}{\acat{F}_{\leq 1}} = \perv{p}{\acat{F}}$, so $\Hilb^\partial_{\leq 1}(Y) = \Hilb_{\leq 1}(Y)$ and so $\curly{K}_{\leq 1} = \curly{H}_{\leq 1}$.
	
	We confess here that it is not clear to us whether the ``partial'' Hilbert scheme is strictly contained the full Hilbert scheme.
	Related to this, we also expect $\perv{p}\Hilb_{\leq 1}$ to be strictly contained in the perverse Hilbert scheme parameterising perverse quotients $\O_Y \onto E$, with $\ch_0(E) = \ch_1(E) =0$, although we have not been able to construct a point in the complement.
\end{rmk}

Once and for all we establish some general notation.
For $\acat{B} \subset \acat{A}$ a subcategory we denote $1_\acat{B}$ the element of $\text{H}_\infty(\acat{A})$ represented by the inclusion of stacks $\astack{B} \subset \astack{A}$, when this is an open immersion (analogous notation for $\acat{A}_{\schnitzel}$ and $\perv{p}{\acat{A}_{\schnitzel}}$).
Another important stack is $\astack{A}_{\schnitzel}^\curlyO$, the stack of \emph{framed} coherent sheaves \cite[Section 2.3]{tom-cc}, which parameterises sheaves with a fixed global section $\curlyO_Y \to \sh{E}$.
By considering surjective sections we can realise $\Hilb_{\schnitzel}(Y)$ as an open subscheme of $\astack{A}_{\schnitzel}^\curlyO$.
We have a forgetful map $\astack{A}_{\schnitzel}^\curlyO \to \astack{A}_{\schnitzel}$, which takes a morphism $\O_Y \to \sh{E}$ to $\sh{E}$.
Given an open substack $\astack{B} \subset \astack{A}_{\schnitzel}$, we can consider the fibre product $\astack{B}^\curlyO = \astack{B}\times_{\astack{A}_{\schnitzel}} \astack{A}_{\schnitzel}^\curlyO$, which gives an element $1_\acat{B}^\curlyO \in \text{H}_\infty(\acat{A}_{\schnitzel})$.

We want to emulate this last construction for $\text{H}_\infty(\perv{p}{\acat{A}_{\schnitzel}})$.
We define the stack $\perv{p}{\astack{A}}^\O$ of \emph{framed perverse coherent sheaves} as the prestack taking a base $S$ to a family of perverse coherent sheaves $P$ together with a morphism $\O_{S \times Y} \to P$.
It is useful for us to realise $\perv{p}{\astack{A}^\O}$ as a fibre product as follows.

Note first that we also have a stack $\astack{C}$ parameterising coherent sheaves on $X$.
Pushforward of complexes induces a morphism of stacks $\perv{p}{\astack{A}} \to \astack{C}$.
In fact, for this to be well-defined, we simply need to check that given a family of perverse coherent sheaves $P$ over a base $S$, the pushforward $Rf_{S,*}P$ is a coherent sheaf.
This can be verified on fibres.
If $s \in S$ is a point, then $Ls^* Rf_{S,*}P = Rf_{s,*}P\vert^L_{Y_s}$\footnote{For a proof of this non-flat base-change, we refer the reader to \cite[Proposition 6.3]{heinrich-cusps}.}, which is a coherent sheaf as $P\vert^L_{Y_s}$ is a perverse coherent sheaf.

Moreover, there is a corresponding stack of framed sheaves $\astack{C}^\curlyO$ \cite[2.3]{tom-cc}.
For $\csh{P} \in \perv{p}{\acat{A}}$, morphisms $\curlyO_Y \to \csh{P}$ correspond (by adjunction) to morphisms $\curlyO_X \to Rf_*\csh{P}$.
We know that $Rf_*\csh{P}$ is a sheaf, so morphisms $\curlyO_Y \to \csh{P}$ correspond to points of $\astack{C}^\curlyO$.

To make the argument work in families, we notice that over a base $S$ we still have $R f_{S,*}\O_{S \times Y} = \O_{S \times X}$ (this follows from flatness of $S \to \Spec \C$ and base change).
Hence, the considerations made above still apply and $\perv{p}{\astack{A}^\O}$ sits in the cartesian diagram below.
%
\begin{center}
	\begin{tikzpicture}
		\matrix (m) [matrix of math nodes, row sep=3em, column sep=3em, text height=1.5ex, text depth=0.25ex]
		{
		\perv{p}{\astack{A}^\curlyO} & \astack{C}^\curlyO\\
		\perv{p}{\astack{A}} & \astack{C}\\
		};
		\path[->,font=\scriptsize]
		(m-1-1) edge (m-1-2)
				edge (m-2-1)
		(m-1-2) edge (m-2-2)
		(m-2-1) edge (m-2-2)
		;
	\end{tikzpicture}
\end{center}
Once again, we have an obvious substack $\perv{p}{\astack{A}^\O_{\schnitzel}}$, which can also be described as the preimage of $\perv{p}{\astack{A}_{\schnitzel}}$.

We have elements $1_{\perv{p}{\acat{F}_{\schnitzel}[1]}}, 1_{\perv{p}{\acat{T}_{\schnitzel}}} \in \text{H}_\infty(\perv{p}{\acat{A}_{\schnitzel}})$ corresponding to the subcategories $\perv{p}{\acat{F}_{\schnitzel}[1]}, \perv{p}{\acat{T}_{\schnitzel}}$ of $\perv{p}{\acat{A}_{\schnitzel}}$.
By taking fibre products with $\perv{p}{\acat{A}_{\schnitzel}}^\curlyO \to \perv{p}{\acat{A}_{\schnitzel}}$ we produce elements $1^\curlyO_{\perv{p}{\acat{F}_{\schnitzel}[1]}}, 1^\curlyO_{\perv{p}{\acat{T}_{\schnitzel}}} \in \text{H}_\infty(\perv{p}{\acat{A}_{\schnitzel}})$.

We also want a \emph{perverse Hilbert scheme} $\perv{p}{\Hilb_{\schnitzel}(Y/X)}$ of $Y$ over $X$ parameterising quotients of $\curlyO_Y$ in $\perv{p}{\acat{A}_{\schnitzel}}$.
One can realise it as an open substack of $\perv{p}{\astack{A}^\O_{\schnitzel}}$.
Indeed, for $\alpha:\curlyO_Y \to \csh{P}$ with $\csh{P} \in \perv{p}{\acat{A}_{\schnitzel}}$, being surjective is equivalent to the cone of $\alpha$ lying in $\perv{p}{\acat{A}}[1]$, which we know to be an open condition on $\perv{p}{\astack{A}^\O_{\schnitzel}}$.
Thus we have an element $\perv{p}{\curly{H}_{\schnitzel}} \in \text{H}_\infty(\perv{p}{\acat{A}_{\schnitzel}})$.

\subsection{A First Identity}
We want to prove the identity
\begin{align}\label{hilb-f=fo-philb}
	\perv{p}{\curly{H}_{\schnitzel}} * 1_{\perv{p}{\acat{F}}_{\schnitzel}[1]} = 1^\curlyO_{\perv{p}{\acat{F}_{\schnitzel}}[1]} * \curly{K}_{\schnitzel}
\end{align}
which we motivated in the beginning of this section (see also Remark \ref{addedremark}). 
The left hand side is represented by a stack $\astack{M}_\text{L}$, parameterising diagrams
\begin{center}
	\begin{tikzpicture}
		\matrix (m) [matrix of math nodes, row sep=3em, column sep=3em, text height=1.5ex, text depth=0.25ex]
		{
		\curlyO_Y &&\\
		\cmplx{P}_1 & \cmplx{E} & \cmplx{P}_2\\
		};
		\path[->,font=\scriptsize]
		(m-1-1) edge[->>] (m-2-1)
		(m-2-1) edge[right hook->] (m-2-2)
		(m-2-2) edge[->>] (m-2-3)
		;
	\end{tikzpicture}
\end{center}
where all objects are in $\perv{p}{\acat{A}_{\schnitzel}}$, the sequence $\csh{P}_1 \into \csh{E} \onto \csh{P}_2$ is exact in $\perv{p}{\acat{A}_{\schnitzel}}$, $\curlyO_Y \onto \csh{P}_1$ is surjective in $\perv{p}{\acat{A}_{\schnitzel}}$ and $\csh{P}_2 \in \perv{p}{\acat{F}_{\schnitzel}[1]}$.\footnote{To be precise, over a base $U$, the groupoid $\astack{M}_\text{L}(U)$ consists of diagrams as above which, upon restricting to fibres of points of $U$, satisfy the required properties. Similar remarks will be implicit for the other stacks we define below.}

The right hand side is represented by a stack $\astack{M}_\text{R}$ parameterising diagrams
\begin{center}
	\begin{tikzpicture}
		\matrix (m) [matrix of math nodes, row sep=3em, column sep=3em, text height=1.5ex, text depth=0.25ex]
		{
		\curlyO_Y && \curlyO_Y\\
		\sh{F}[1] & \csh{E} & \sh{T}\\
		};
		\path[->,font=\scriptsize]
		(m-1-1) edge (m-2-1)
		(m-1-3) edge node[auto]{\tiny sur} (m-2-3)
		(m-2-1) edge[right hook->] (m-2-2)
		(m-2-2) edge[->>] (m-2-3)
		;
	\end{tikzpicture}
\end{center}
where the horizontal maps form a short exact sequence in $\perv{p}{\acat{A}_{\schnitzel}},$ $\sh{F} \in \perv{p}{\acat{F}_{\schnitzel}}, \sh{T} \in \perv{p}{\acat{T}_{\schnitzel}}$ and the map $\curlyO_Y \to T$ is both surjective as a morphism in $\acat{A}_{\schnitzel}$ and has perverse cokernel lying in $\perv{p}{\acat{F}_{\leq 1}}[1]$. 
We remind ourselves that, given a perverse coherent sheaf $\csh{E} \in \perv{p}{\acat{A}_{\leq 1}}$, there is a unique exact sequence $\sh{F}[1] \into \csh{E} \onto \sh{T}$, with $\sh{F}\in \perv{p}{\acat{F}_{\schnitzel}}, \sh{T}\in \perv{p}{\acat{T}_{\schnitzel}}$. 

As the proof of the required identity goes through a chain of geometric bijections and Zariski fibrations, we draw a diagram for future reference.
\begin{center}
	\begin{tikzpicture}
		\matrix (m) [matrix of math nodes, row sep=3em, column sep=3em, text height=1.5ex, text depth=0.25ex]
		{
		\astack{M}_\text{L} & & \astack{M} && \astack{M}_\text{R} \\
		& \astack{M}' && \astack{N} & \\
		};
		\path[->,font=\scriptsize]
		(m-1-1) edge (m-2-2)
		(m-1-3) edge (m-2-2)
		(m-1-3) edge (m-2-4)
		(m-1-5) edge (m-2-4)
		;
	\end{tikzpicture}
\end{center}

We now define a stack $\astack{M}'$ parameterising diagrams of the form
\begin{center}
	\begin{tikzpicture}
		\matrix (m) [matrix of math nodes, row sep=3em, column sep=3em, text height=1.5ex, text depth=0.25ex]
		{\curlyO_Y \\
		\csh{E}\\
		};
		\path[->,font=\scriptsize]
		(m-1-1) edge node[auto]{$\varphi$} (m-2-1)
		;
	\end{tikzpicture}
\end{center}
where $\perv{p}{\coker} \, \varphi \in \perv{p}{\curvy{F}_{\schnitzel}}[1].$
By Lemma \ref{previous-lemma}, this last condition is equivalent to $\cone(\varphi) \in \text{D}^{\leq -1}(Y)$, which is open.
Thus $\astack{M}'$ is an open substack of the stack of framed perverse sheaves $\perv{p}{\astack{A}}_{\schnitzel}^\curlyO$.

\begin{prop}\label{backtothefuture}
	There is a map $\astack{M}_\text{L} \to \astack{M}'$ induced by the composition $\curlyO_Y \onto \csh{P}_1 \into \csh{E}.$
	This map is a geometric bijection.
\end{prop}
\begin{prf}
	By taking the composition $\O_Y \to \csh{P}_1 \to \csh{E}$ in the diagram defining $\astack{M}_\text{L}$ (and using the previous lemma) we see that $\astack{M}_\text{L} \to \astack{M}'$ is an equivalence on $\C$-points.
	To prove finite typeness of the morphism we use a fact that shall be proved later: $\perv{p}{\curly{H}_{\leq 1}}$, $1_{\perv{p}{\acat{F}_{\leq 1}[1]}}$ are \emph{Laurent} elements of our Hall algebra (Propositions \ref{prop-1flaurent}, \ref{prop-laurent-laurent}), in the sense of Definition \ref{laurent-definition}.
	As the stack $\astack{M}_\text{L}$ is the product of these two elements, it is also Laurent.
	Thus, for any numerical class $\alpha$, we have a morphism $\astack{M}_{\text{L},\alpha} \to \astack{M}'_\alpha$.
	As $\astack{M}_{\text{L},\alpha}$ is of finite type, we are done.
\end{prf}
We define another stack $\astack{M}$ parameterising diagrams of the form
\begin{center}
	\begin{tikzpicture}
		\matrix (m) [matrix of math nodes, row sep=3em, column sep=3em, text height=1.5ex, text depth=0.25ex]
		{
		&\curlyO_Y&\\
		\sh{F}[1] & \csh{E} & \sh{T}\\
		};
		\path[->,font=\scriptsize]
		(m-1-2) edge node[auto]{$\varphi$} (m-2-2)
		(m-2-1) edge[right hook->] (m-2-2)
		(m-2-2) edge[->>] (m-2-3)
		;
	\end{tikzpicture}
\end{center}
where the horizontal maps form a short exact sequence of perverse sheaves, $\sh{F} \in \perv{p}{\acat{F}_{\schnitzel}},$ $\sh{T} \in \perv{p}{\acat{T}_{\schnitzel}}$ and $\perv{p}{\coker}\, \varphi \in \perv{p}{\acat{F}_{\schnitzel}[1]}$.
This stack can be obtained as a fibre product as follows.
The element $1_{\perv{p}{\acat{F}_{\schnitzel}}[1]}* 1_{\perv{p}{\acat{T}_{\schnitzel}}}$ is represented by a morphism $Z \to \perv{p}{\astack{A}_{\schnitzel}}$ and $\astack{M}$ is the top left corner of the following cartesian diagram.
\begin{center}
	\begin{tikzpicture}
		\matrix (m) [matrix of math nodes, row sep=3em, column sep=3em, text height=1.5ex, text depth=0.25ex]
		{
		\astack{M} & \astack{M}'\\
		Z & \perv{p}{\astack{A}_{\schnitzel}}\\
		};
		\path[->,font=\scriptsize]
		(m-1-2) edge (m-2-2)
		(m-2-1) edge (m-2-2)
		(m-1-1) edge (m-2-1)
		(m-1-1) edge (m-1-2)
		;
	\end{tikzpicture}
\end{center}
\begin{prop}
	The morphism $\astack{M} \to \astack{M}'$ defined by forgetting the exact sequence is a geometric bijection.
\end{prop}
\begin{prf}
	The morphism in question is precisely the top row of the previous diagram.
	The bottom row is obtained by composing the top arrows of the following diagram.
	\begin{center}
		\begin{tikzpicture}
			\matrix (m) [matrix of math nodes, row sep=3em, column sep=3em, text height=1.5ex, text depth=0.25ex]
			{
			Z & \perv{p}{\astack{A}_{\schnitzel}^{(2)}} & \perv{p}{\astack{A}_{\schnitzel}} \\
			\perv{p}{\astack{F}_{\schnitzel}[1]} \times \perv{p}{\astack{T}_{\schnitzel}} & \perv{p}{\astack{A}_{\schnitzel}} \times \perv{p}{\astack{A}_{\schnitzel}} & \\
			};
			\path[->,font=\scriptsize]
			(m-1-1) edge (m-1-2)
					edge (m-2-1)
			(m-2-1) edge (m-2-2)
			(m-1-2) edge node[auto]{$b$} (m-1-3)
					edge (m-2-2)
			;
		\end{tikzpicture}
	\end{center}
	where the bottom row is an open immersion (and thus of finite type) and the morphism $b$ is of finite type (this follows from the fact that $b$ locally is isomorphic to the analogous morphism for coherent sheaves).
	The morphism $Z \to \perv{p}{\astack{A}_{\schnitzel}}$ induces an equivalence on $\C$-points because $(\perv{p}{\acat{F}[1]},\perv{p}{\acat{T}_{}})$ is a torsion pair in $\perv{p}{\acat{A}_{}}$ (and thus any perverse coherent sheaf has a unique short exact sequence with torsion kernel and torsion-free cokernel) and because an automorphism of a short exact sequence which is the identity on the middle term is trivial.
	As $\astack{M} \to \astack{M}'$ is a base change of $Z \to \perv{p}{\astack{A}_{\schnitzel}}$ we are done.
\end{prf}
Thus the identity \eqref{hilb-f=fo-philb} boils down to proving that $\astack{M}$ and $\astack{M}_\text{R}$ represent the same element in $\text{H}_\infty(\perv{p}{\acat{A}})$.
To do this we use one last stack $\astack{N}$ and build a pair of Zariski fibrations with same fibres.
We define the stack $\astack{N}$ to be the moduli of the following diagrams
\begin{equation}\label{flops zar fib diagram}
	\begin{tikzpicture}[baseline=(current  bounding  box.center)]
		\matrix (m) [matrix of math nodes, row sep=3em, column sep=3em, text height=1.5ex, text depth=0.25ex]
		{
		&& \curlyO_Y\\
		\sh{F}[1] & \csh{E} & \sh{T}\\
		};
		\path[->,font=\scriptsize]
		(m-1-3) edge node[auto]{\tiny sur} (m-2-3)
		(m-2-1) edge[right hook->] (m-2-2)
		(m-2-2) edge[->>] (m-2-3)
		;
	\end{tikzpicture}
\end{equation}
where the horizontal maps form a short exact sequence of perverse sheaves, $\sh{F} \in \perv{p}{\acat{F}_{\schnitzel}},$ $\sh{T} \in \perv{p}{\acat{T}_{\schnitzel}}$ and the map $\curlyO_Y \to \sh{T}$ is both surjective as a morphism of coherent sheaves and has perverse cokernel lying in $\perv{p}{\acat{F}_{\leq 1}}[1]$. 
This stack is also a fibre product of known stacks (compare with the element $1_{\perv{p}{\acat{F}_{\schnitzel}}[1]}* \curly{K}_{\schnitzel}$). 
Notice that there are two maps $\astack{M}\to \astack{N} \ot \astack{M}_\text{R}.$
The map $\astack{M}_\text{R} \to \astack{N}$ is given by forgetting the morphism $\curlyO_Y \to \sh{F}[1]$.
The map $\astack{M}\to \astack{N}$ is given by composition $\curlyO_Y \to \csh{E} \to \sh{T}$ (which is a surjective morphism thanks to Lemma \ref{previous-lemma}).
\begin{prop}\label{sucaiurgen}
	The maps $\astack{M} \to \astack{N} \ot \astack{M}_{\text{R}}$ are two Zariski fibrations with the same fibres.
\end{prop}
\begin{prf}
	Keeping in mind diagram \eqref{flops zar fib diagram}, the idea is that over a perverse coherent sheaf $\csh{E}$ the morphism $\astack{M}_{\text{R}} \to \astack{N}$ has fibres $\Hom_Y(\curlyO_Y,\sh{F}[1])$ while $\astack{M} \to \astack{N}$ has fibres lifts $\curlyO_Y \to \csh{E}$ such that the perverse cokernel lies in $\perv{p}{\acat{F}_{\schnitzel}[1]}$.\footnote{
		This is where we thank J{\o}rgen Rennemo for spotting a tricky subtlety.
		When the singular locus of $X$ is not zero-dimensional, the condition $\perv{p}\coker \in \perv{p}{\acat{F}_{\schnitzel}[1]}$ is not automatic and the argument ceases to work unless we either replace $\curly{H}$ by $\curly{K}$ or we restrict ourselves to the smaller category $\perv{p}{\acat{A}_{\exc}}$ (see the next section).
	} 
	The long exact sequence 
	\begin{align*}
		0 \to \Hom_Y(\curlyO_Y,\sh{F}[1]) \to \Hom_Y(\curlyO_Y,\csh{E}) \to \Hom_Y(\curlyO_Y,\sh{T}) \to 0
	\end{align*}
	tells us that given a choice of a lift of $\curlyO_Y \to \sh{T}$ all lifts are in bijection with $\Hom_Y(\curlyO_Y,\sh{F}[1])$. 
	We will now show that in fact \emph{any} lift of $\O_Y \to T$ is such that the corresponding perverse cokernel lies in $\perv{p}{\acat{F}_{\leq 1}}[1]$.
	
	Let $\psi\colon \O_Y \to T$ be a morphism with $T \in \perv{p}{\acat{T}_{\leq 1}}$ and $\perv{p}{\coker\psi} \in \perv{p}{\acat{F}_{\leq 1}[1]}$.
	Let $\varphi\colon \O_Y \to E$ be any lift.
	Consider the following diagram
	\begin{center}
		\begin{tikzcd}
			0 \ar{d}\ar{r} &
			\O_Y \ar[-,double equal sign distance]{r} \ar{d}{\varphi} & \O_Y \ar{d}{\psi} \\
			F[1] \ar[hook]{r} & E \ar[->>]{r} & T 
		\end{tikzcd}
	\end{center}
	and notice that the rows are short exact sequences.
	The snake lemma implies an exact sequence
	$$ F[1] \to \perv{p}{\coker}\varphi \to \perv{p}{\coker \psi} \to 0 $$
	and using Lemma \ref{closedunderquotients} it follows that $\perv{p}{\coker}\varphi \in \perv{p}{\acat{F}_{\leq 1}[1]}$.
	
	Let's see how to make this argument work in families.
	Let $S$ be an affine and connected scheme and let $S \to \astack{N}$ correspond to a diagram 
	\begin{center}
		\begin{tikzpicture}
			\matrix (m) [matrix of math nodes, row sep=3em, column sep=3em, text height=1.5ex, text depth=0.25ex]
			{
			&& \curlyO_Y\\
			\sh{F}[1] & \csh{E} & \sh{T}\\
			};
			\path[->,font=\scriptsize]
			(m-1-3) edge node[auto]{\tiny sur} (m-2-3)
			(m-2-1) edge[right hook->] (m-2-2)
			(m-2-2) edge[->>] (m-2-3)
			;
		\end{tikzpicture}
	\end{center}
	on $Y_S$.
	First of all notice that Lemma \ref{cohomology} (and base change to $S$) tells us that $Rp_{S,*}\sh{F}$ is just $H^1(Y_S,\sh{F})$ shifted by one, where $p_S:Y_S \to S$ is the projection.
	In addition, $H^1(Y_S,\sh{F})$ is flat over $S$, or in other words $\curlyO_{Y_S}$ and $\sh{F}$ have constant $\Ext$ groups in the sense of \cite[Section 6.1]{tom-hall} (all the others vanish).
	
	Let $W$ be the fibre product $\astack{M}_{\text{R}} \times_{\astack{N}} S$.
	This is actually a functor which associates to an affine $S$-scheme $q:T \to S$ the group $H^1(Y_T,q_Y^*\sh{F})$ and we know by loc.~cit.~that it is represented by a vector bundle over $S$ of rank the rank of $H^1(Y_S,\sh{F})$.
	
	Similarly, the fibre product $\astack{M} \times_{\astack{N}} S$ is represented by an affine bundle of rank the rank of $H^1(Y_S,\sh{F})$ (notice that because of the previous arguments the exact sequence at the beginning of the proof still holds over $S$).
	This allows us to conclude that $\astack{M} \to \astack{N}$ and $\astack{M}_\text{R} \to \astack{N}$ are two Zariski fibrations with same fibres.
\end{prf}

\subsection{PT Invariants} 
\label{sub:duality}
We are still left with the task of understanding what we obtain by integrating $\perv{p}{\curly{H}_{\schnitzel}}.$
To achieve this goal we first substitute $1^\curlyO_{\perv{p}{\acat{F}_{\schnitzel}}[1]}$ with something more recognisable (from the point of view of the integration morphism $I$).
Recall \cite[Section 2.2]{tom-cc} that on $\acat{A}$ there is a torsion pair $(\acat{P},\acat{Q})$, where $\acat{P}$ consists of sheaves supported in dimension zero and $\acat{Q}$ is the right orthogonal of $\acat{P}$.
In particular, an element $Q \in \acat{Q}$ which is supported in dimension one is pure.
Notice also that $\curlyO_Y \in \acat{Q}$.
We denote by $\acat{A}^\#$ the tilt with respect to $(\acat{P},\acat{Q})$, but with the convention
\begin{align*}
	\acat{P}[-1] \subset \acat{A}^\# \subset \text{D}^{[0,1]}(Y).
\end{align*}

There exists a scheme $\Hilb_{\schnitzel}^\#(Y)$ parameterising quotients of $\curlyO_Y$ in $\acat{A}^\#$ supported in dimension at most one.
Using \cite[Lemma 2.3]{tom-cc} one constructs an element $\curly{H}_{\schnitzel}^\# \in \text{H}_\infty(\acat{A}_{\schnitzel})$ which eventually leads to the PT invariants of $Y$.
We recall that quotients of $\curlyO_Y$ in $\acat{A}^\#$ are exactly morphisms $\curlyO_Y \to \sh{Q}$, with cokernel in $\acat{P}$ and $\sh{Q} \in \acat{Q}$.

In $\text{H}_\infty(\acat{A}_{\schnitzel})$ we have an element $1_{\acat{Q}_{\schnitzel}}$ given by the inclusion of the stack parameterising objects in $\acat{Q}_{\schnitzel}$ inside $\astack{A}_{\schnitzel}$ and its framed version $1^\curlyO_{\acat{Q}_{\schnitzel}}$.
There is also an identity \cite[Section 4.5]{tom-cc}
\begin{align*}
	1^\curlyO_{\acat{Q}_{\schnitzel}} = \curly{H}_{\schnitzel}^\# * 1_{\acat{Q}_{\schnitzel}}.
\end{align*}
We want to restrict the element $\curly{H}_{\schnitzel}^\#$ further by considering only quotients whose derived pushforward $\R f_*$ is supported in dimension zero.
We thus define the following subcategories.
\begin{align*}
	\acat{Q}_{\exc} &= \left\{ \sh{Q}\in\curvy{Q}_{\leq 1} \st \dim \supp \R f_*\sh{Q} = 0 \right\}\\
	\perv{p}{\acat{A}_{\exc}} &= \left\{ \csh{E} \in \perv{p}{\acat{A}_{\leq 1}} \st \dim \supp \R f_* \csh{E} = 0 \right\} \\
	\perv{p}{\acat{T}_{\exc}} &= \perv{p}{\acat{T}} \cap \perv{p}{\acat{A}_{\exc}} \\
	\perv{p}{\acat{T}_\bullet} &= \perv{p}{\acat{T}_{\exc}} \cap \acat{Q}_{\exc}
\end{align*}
We can also consider the scheme $\Hilb_{\exc}^\#(Y)$ parameterising quotients of $\curlyO_Y$ in $\acat{A}_{\schnitzel}^\#$ with target having zero-dimensional pushdown (it is indeed an open  subscheme of $\Hilb_{\schnitzel}^\#(Y)$ as we are imposing a restriction on the numerical class of the quotients).
From it we obtain an element $\curly{H}_{\exc}^\# \in \text{H}_\infty(\acat{A})$.
Before we move on to the following result, we point out that $\perv{p}{\acat{T}}_\bullet \subset {\acat{A}}_{\schnitzel}^\#$
\begin{prop}
	The following identity in $\text{H}_\infty(\acat{A})$ is true.
	\begin{align}\label{1ot=hs*1t}
		1_{\perv{p}{\acat{T}_\bullet}}^\curlyO = \curly{H}_{\exc}^\# * 1_{\perv{p}{\acat{T}_\bullet}}
	\end{align}
\end{prop}
\begin{prf}
	We start with a remark.
	If we have a morphism $\curlyO_Y \to \sh{T}$ in $\acat{A}^\#$, with $\sh{T} \in \perv{p}{\acat{T}_\bullet}$, we can factor it through its image (in $\acat{A}^\#$) $\curlyO_Y \to \sh{I} \to \sh{T}$ and we denote by $\sh{Q}$ the quotient, again in $\acat{A}^\#$.
	We already know \cite[Lemma 2.3]{tom-cc} that $\sh{I}$ is a sheaf and that the morphism $\curlyO_Y \to \sh{I}$, as a morphism in $\acat{A}$, has cokernel $\sh{P}$ supported in dimension zero.
	
	Glancing at the cohomology sheaves long exact sequence of $I \to T \to Q$, reveals that $\sh{Q}$ is also a sheaf, thus the sequence $\sh{I} \into \sh{T} \onto \sh{Q}$ is actually a short exact sequence of sheaves.
	The sheaf $\sh{Q}$ is in $\perv{p}{\acat{T}}$, as it is a quotient of $\sh{T}$, and it lies in $\acat{Q}$ as it is an object of $\acat{A}^\#$.
	Also, $\R f_*\sh{Q}$ is supported on points as $\R f_*\sh{T}$ is, thus $\sh{Q} \in \perv{p}{\acat{T}_\bullet}$.
	
On the other hand, given a morphism of sheaves $\curlyO_Y \to \sh{I}$, which is an epimorphism in $\acat{A}^\#$, and given a short exact sequence of coherent sheaves $\sh{I} \into \sh{T} \onto \sh{Q}$, with $\sh{I} \in \acat{Q}_{\exc}$ and $\sh{Q} \in \perv{p}{\acat{T}_\bullet}$, we claim that $\sh{T} \in \perv{p}{\acat{T}_\bullet}$.
The fact that $\sh{T} \in \acat{Q}_{\exc}$ is clear, if we prove that $\sh{I} \in \perv{p}{\acat{T}}$ then we are done.

We know there is an exact sequence $\curlyO_Y \to \sh{I} \onto \sh{P}$, with $\sh{P}$ supported in dimension zero, viz.~a skyscraper sheaf.
Let $\sh{I} \onto \sh{F}$ be the projection to the torsion-free part of $\sh{I}$ (for the $(\perv{p}{\acat{T}}, \perv{p}{\acat{F}})$ torsion pair).
The morphism $\curlyO_Y \to \sh{I} \onto \sh{F}$ is zero, as objects of $\perv{p}{\acat{F}}$ have no sections.
Thus there is a morphism $\sh{P} \to \sh{F}$ such that $\sh{I} \onto \sh{P} \to \sh{F}$ is equal to $\sh{I} \onto \sh{F}$.
As $\sh{P}$ is a skyscraper sheaf, the morphisms from it are determined on global sections, thus $\sh{P} \to \sh{F}$ is zero, which in turn implies that $\sh{I} \onto \sh{F}$ is zero.
Thus $\sh{F} = 0$ and $\sh{I} \in \perv{p}{\acat{T}}$.

Using the remark above we can see that there exists a morphism from the stack parameterising diagrams
\begin{center}
	\begin{tikzpicture}
		\matrix (m) [matrix of math nodes, row sep=3em, column sep=3em, text height=1.5ex, text depth=0.25ex]
		{
		\curlyO_Y & & \\
		\sh{I} & \sh{T} & \sh{Q}\\
		};
		\path[->,font=\scriptsize]
		(m-1-1) edge (m-2-1)
		(m-2-1) edge[right hook->] (m-2-2)
		(m-2-2) edge [->>] (m-2-3)
		;
	\end{tikzpicture}
\end{center}
with $\curlyO_Y \to \sh{I}$ an epimorphism in $\acat{A}^\#$, $\sh{I} \in \acat{Q}_{\exc}$, $\sh{Q} \in \perv{p}{\acat{T}_\bullet}$, to the stack parameterising morphisms $\curlyO_Y \to \sh{T}$, with $\sh{T} \in \perv{p}{\acat{T}_\bullet}$.
This morphism induces an equivalence on $\C$-points and the fact that it is of finite type will follow from Proposition \ref{prop-1flaurent} and Proposition \ref{dualityprop}.
\end{prf}

\subsection{Duality}
We will see now how to link everything together via the duality functor.
\begin{lem}\label{duality lemma}
	Let $\D : \text{D}(Y) \to \text{D}(Y)$ be the anti-equivalence defined by
	\begin{align*}
		\csh{E} \longmapsto \D\left(\csh{E}\right) = R\lHom_Y(\csh{E},\curlyO_Y)[2].
	\end{align*}
	Then
	\begin{align*}
		\D\left(\perv{q}{\acat{T}_\bullet}\right) &= \perv{p}{\acat{F}}_{\leq 1} 
	\end{align*}
	for $q=-(p+1)$.
\end{lem}
The shift $[2]$ in the definition of $\D$ is due to the fact we are dealing with pure sheaves supported in codimension two.
Indeed, if $\acat{Q}_{\leq 1}$ is the category of pure sheaves supported in dimension one, then $\D(\acat{Q}_{\leq 1}) = \acat{Q}_{\leq 1}$ \cite[Lemma 5.6]{tom-cc}.
Notice that any sheaf $F \in \perv{p}{\acat{F}_{\leq 1}}$ is automatically pure, as the existence of a zero-dimensional subsheaf would contradict the condition $f_* F = 0.$
\begin{prf}
	We will prove the two inclusions $\D(\perv{q}{\acat{T}_\bullet}) \subset \perv{p}{\acat{F}_{\leq 1}}$, $\D(\perv{p}{\acat{F}_{\leq 1}}) \subset \perv{q}{\acat{T}_\bullet}$, but first let us make a consideration about the category $\acat{C}_{\leq 1}$ of coherent sheaves supported in dimension at most one with vanishing derived pushforward.
	We have $\D(\acat{C}_{\leq 1}) = \acat{C}_{\leq 1}$.
	In fact, as $\acat{C}_{\leq 1} \subset \acat{Q}_1$, one has $\D(\acat{C}_{\leq 1}) \subset \acat{Q}_1$, thus one needs only check $Rf_*\D(C) = 0$, for all $C \in \acat{C}_{\leq 1}$.
	\begin{align*} 
		Rf_*\D(\sh{C})
		&= Rf_*R\lHom_Y(\sh{C},\curlyO_Y)[2] \\
		&= Rf_*R\lHom_Y(\sh{C},f^!\curlyO_X)[2] \\
		&= R\lHom_X(Rf_*\sh{C},\curlyO_X)[2]=0
	\end{align*}

	Let $F \in \perv{p}{\acat{F}_{\leq 1}}$.
	We first check that $R^1f_*\D(\sh{F})=0$.
	\begin{align*}
		R^1f_* \D(\sh{F})
		&= H^1\left( Rf_*R\lHom_Y(\sh{F},\curlyO_Y))[2] \right)\\
		&= H^3\left( R\lHom_X\left( Rf_*\sh{F},\curlyO_X \right) \right)\\
		&= H^3\left(R\lHom_Y(R^1f_*\sh{F}[-1],\curlyO_X)\right)\\
		&= \lExt_X^4\left( R^1f_*\sh{F},\curlyO_X \right)\\
		&= \Ext_X^4\left( R^1f_*\sh{F},\curlyO_X \right) = 0
	\end{align*}
	where the last equality follows from Serre duality and the second to last is a consequence of the local-to-global spectral sequence and the fact that $R^1f_*\sh{F}$ (and thus $\lExt_Y^4(R^1f_*\sh{F},\curlyO_Y)$) is supported in dimension zero.
	When $p = -1$ this is enough to show that $\D(\perv{p}{\acat{F}_{\leq 1}}) \subset \perv{q}{\acat{T}_\bullet}$.
	When $p = 0$ we are still left to check that $\Hom_Y(\D(\perv{p}{\acat{F}_{\leq 1}}), \acat{C}) = 0$.
	If $F \in \perv{p}{\acat{F}_{\leq 1}}$, then (using the fact that $\D$ is an antiequivalence of $D(Y)$)
	\begin{align*}
		\Hom_Y(\D(F), \curvy{C}_{\leq 1}) = \Hom_Y(\curvy{C}_{\leq 1}, F) \subset \Hom_Y(\acat{C}, F) = 0
	\end{align*}
	where the last equality is by definition of $\perv{0}{\acat{F}}$.
	To complete the proof, we show that if $T \in \acat{A}_{\leq 1}$ is such that $R^1f_* T = 0$ and $\Hom_Y(T, \acat{C}_{\leq 1}) = 0$, then $\Hom_Y(T, \acat{C}) = 0$.
	In fact, let $T \to C$ be a morphism with $C \in \acat{C}$.
	The image $I$ satisfies $R^1f_* I = 0$ as it is a quotient of $T$ and $f_* I = 0$ as it is a subobject of $C$.
	Observing that $T \onto I$ is surjective implies that $I \in \acat{C}_{\leq 1}$ and that $T \to C$ is the zero morphism.
	
	Let now $T \in \perv{q}{\acat{T}_\bullet}$, we check that $f_*\D(\sh{T}) = 0$.
	\begin{align*}
		f_*\D(\sh{T})
		&= H^0\left(Rf_*R\lHom_Y\left(\sh{T},\curlyO_Y\right)[2]\right)\\
		&= H^2\left(R\lHom_X\left(f_*\sh{T},\curlyO_X\right)\right)\\
		&= \lExt^2_X\left(f_*\sh{T},\curlyO_X\right)\\
		&= \Ext^2_X\left(f_*\sh{T},\curlyO_X\right) = 0
	\end{align*}
	where the last two equalities again follow from Serre duality and the dimension of the support of $f_*\sh{T}$.
	Analogously as above, this is enough for $p = -1$, and for $p = 0$ we see that $\Hom_Y(\acat{C}_{\leq 1}, \D(\perv{q}{\acat{T}_\bullet})) = 0$.
	
	Let now $F \in \acat{A}_{\leq 1}$ be such that $f_* F = 0$ and $\Hom_Y(\acat{C}_{\leq 1}, F) = 0$.
	It follows that $\Hom_Y(\acat{C},F) = 0$.
	In fact, if $C \to F$ is a morphism with $C \in \acat{C}$, then the image $I$ satisfies $f_* I = 0$ as it is a subobject of $F$ and satisfies $R^1f_* I = 0$ as it is a quotient of $C$.
	As $I \into F$ is injective, $I \in \acat{C}_{\leq 1}$ which implies that $C \to F$ is the zero morphism.
\end{prf}
We now want to apply the duality functor, or better $\D' = \D[1]$, to our Hall algebras.
As the category $\perv{p}{\acat{F}_{\leq 1}[1]}$  (respectively $\perv{q}{\acat{T}}_\bullet$) is closed by extensions we have an algebra $\text{H}_\infty(\perv{p}{\acat{F}_{\leq 1}[1]})$ (respectively $\text{H}_\infty(\perv{q}{\acat{T}_\bullet})$) spanned by morphisms $[W \to \perv{p}{\astack{F}_{\leq 1}[1]}]$ (respectively $[W \to \perv{q}{\astack{T}_\bullet}]$).
Notice that while the first is a subalgebra of $\text{H}_\infty(\perv{p}{\acat{A}_{\leq 1}})$, the second can be viewed as a subalgebra of both $\text{H}_\infty(\perv{q}{\acat{A}_{\leq 1}})$ and $\text{H}_\infty(\acat{A}_{\leq 1})$, as a distinguished triangle with vertices lying in $\perv{q}{\acat{T}_\bullet}$ is an exact sequence in both $\perv{q}{\acat{A}}$ and $\acat{A}$.
\begin{prop}\label{dualityprop}
The functor $\D'$ induces an anti-isomorphism between $\text{H}_\infty(\perv{q}{\acat{T}}_\bullet)$ and $\text{H}_\infty(\perv{p}{\acat{F}_{\leq 1}[1]})$.
	Furthermore the following identities hold.
	\begin{align*}
		\D'\left( 1_{\perv{q}{\acat{T}_\bullet}} \right) &= 1_{\perv{p}{\acat{F}_{\leq 1}[1]}} \\
		\D'\left( 1_{\perv{q}{\acat{T}_\bullet}}^\curlyO \right) &= 1_{\perv{p}{\acat{F}_{\leq 1}[1]}}^\curlyO
	\end{align*}
\end{prop}
\begin{prf}
	Duality $\D'$ induces an isomorphism between stacks $\perv{q}{\astack{T}}_\bullet$ and $\perv{p}{\astack{F}_{\leq 1}[1]}$.
	The anti-isomorphism between the Hall algebras is then defined by taking a class $[W \to \perv{q}{\astack{T}_\bullet}]$ to $[W \to \perv{q}{\astack{T}_\bullet} \to \perv{p}{\astack{F}_{\leq 1}[1]}]$ and noticing that duality flips extensions \cite[Section 5.4]{tom-cc}.
	Clearly this takes the element $1_{\perv{q}{\acat{T}_\bullet}}$ to $1_{\perv{p}{\acat{F}_{\leq 1}[1]}}$, while the second identity requires a bit of work.
	
	Two remarks are in order.
	The first is that given any $\sh{T} \in \perv{q}{\acat{T}_\bullet}$,
	\begin{align*}
	\Hom_Y(\curlyO_Y,\sh{T})=\Hom_Y(\D'(\sh{T}),\curlyO_Y[3])=\Hom_Y(\curlyO_Y,\D'(T))^\vee.
	\end{align*}
	The second is that, if $\sh{T} \in \perv{q}{\acat{T}}_\bullet$ and $\sh{F} \in \perv{p}{\acat{F}_{\leq 1}}$, then $\dim_\C H^0(Y,\sh{T}) = \chi(\sh{T})$ and similarly $\dim_\C H^1(Y,\sh{F}) = - \chi(\sh{F})$.
	This is useful since, for a family of coherent sheaves, the Euler characteristic is locally constant on the base.
	Thus we can decompose the stack $\perv{q}{\astack{T}_\bullet}$ as a disjoint union according to the value of the Euler characteristic.
	We have a corresponding decomposition of $\perv{q}{\astack{T}_\bullet}^\O$ and we write $\perv{q}{\astack{T}_{\bullet,n}^\curlyO}$ for the nth component of this disjoint union.
	This space maps down to $\perv{q}{\astack{T}_{\bullet,n}}$ by forgetting the section.
	Similarly, the space $\A^n \times \perv{q}{\astack{T}_{\bullet,n}}$ projects onto $\perv{q}{\astack{T}_{\bullet,n}}$.
	As these two maps are Zariski fibrations with same fibres the stacks $\perv{q}{\astack{T}_{\bullet,n}^\curlyO}$ and $\A^n \times \perv{q}{\astack{T}_{\bullet,n}}$ represent the same element in the Grothendieck ring.
	This argument is then extended to the whole $\perv{q}{\astack{T}_{\bullet}^\curlyO}$ proving that
	\begin{align*}
		\left[\perv{q}{\astack{T}_{\bullet}^\curlyO}\right] =
		\left[ \coprod_n \A^n \times \perv{q}{\astack{T}_{\bullet,n}} \right].
	\end{align*}
	We can proceed analogously for $\perv{p}{\astack{F}_{\leq 1}[1]}$.
	The component $\perv{p}{\astack{F}_{\leq 1}[1]_n^\curlyO}$ represents the same element as $\A^n \times \perv{p}{\astack{F}_{\leq 1}[1]_n}$.
	The first remark above implies that duality $\D'$ takes $\perv{q}{\astack{T}_{\bullet,n}}$ to $\perv{p}{\astack{F}_{\leq 1}[1]_n}$, which lets us conclude.
\end{prf}
Thus in our infinite-type Hall algebra we deduce that $1^\curlyO_{\perv{p}{\acat{F}_{\leq 1}[1]}} = \mathbb{D}'(1^\curlyO_{\perv{q}{\acat{T}_\bullet}}) = \mathbb{D}'(\curly{H}^\#_{\exc} * 1_{\perv{q}{\acat{T}_\bullet}}) = 1_{\perv{p}{\acat{F}_{\leq 1}}[1]} * \mathbb{D}'(\curly{H}^\#_{\exc})$.
Accordingly, we have
\begin{align*}
	1^\curlyO_{\perv{p}{\acat{F}_{\leq 1}[1]}} = 
	1_{\perv{p}{\acat{F}_{\leq 1}}[1]} * \mathbb{D}'(\curly{H}^\#_{\exc})
\end{align*}
which, combined with \eqref{hilb-f=fo-philb}, yields
\begin{align*}
	\perv{p}{\curly{H}_{\leq 1}} * 1_{\perv{p}{\acat{F}_{\leq 1}}[1]} = 1_{\perv{p}{\acat{F}_{\leq 1}}[1]} * \mathbb{D}'(\curly{H}^\#_{\exc}) * \curly{K}_{\leq 1}.
\end{align*}
 

\subsection{Laurent Elements} 
\label{sub:laurent_elements}
Our objective is to get rid of the spurious $1_{\perv{p}{\acat{F}_{\leq 1}[1]}}$'s in the identity above.
This is achieved by constructing a (weak) stability condition (in the sense of \cite[Definition 3.5]{joyce-book}) with values in the ordered set $\{1,2\}$, such that $\perv{p}{\acat{F}_{\leq 1}[1]}$ manifests as the class of semi-stable objects of $\mu=2$.
Before we do that, however, we want to define a sort of completed Hall algebra $\text{H}(\perv{p}{\acat{A}})_{\Lambda}$ (parallel to the one in \cite[Section 5.2]{tom-cc}) which morally sits in between $\text{H}(\perv{p}{\acat{A}}_{\leq 1})$ and $\text{H}_\infty(\perv{p}{\acat{A}_{\leq 1}})$.
The reason we need to do so is simple.
On one hand the Hall algebra constructed in the previous section only includes spaces that are of finite type, on the other the infinite type Hall algebra is much too big to support an integration morphism.
To deal with objects such as the Hilbert scheme of curves and points of $Y$ we allow our spaces to be \emph{locally} of finite type while imposing a \emph{Laurent} condition.

We previously mentioned that $\text{H}(\perv{p}{\acat{A}})$ is graded by the numerical Grothendieck group $N(Y)$.
There is a subgroup $N_{\leq 1}(Y)$ generated by sheaves supported in dimension at most one and $\text{H}(\perv{p}{\acat{A}}_{\leq 1})$ is graded by it.
We also notice \cite[Lemma 2.2]{tom-cc} that the Chern character induces an isomorphism
\begin{align*}
	N_{\leq 1}(Y) \ni [\sh{E}] \longmapsto (\ch_2 \sh{E}, \ch_3  \sh{E}) \in N_1(Y) \oplus N_0(Y)
\end{align*}
where by $N_1(Y)$ we mean the group of curve-classes modulo numerical equivalence, and $N_0(Y) \simeq \Z$.
Henceforth we tacitly identify $N_{\leq 1}(Y)$ with $N_1(Y) \oplus \Z$.

We have a pushforward morphism $f_*\colon N_1(Y) \to N_1(X)$.
This morphism is surjective and we denote its kernel by $N_1(Y/X)$.
The short exact sequence
\begin{align*}
	N_1(Y/X) \into N_1(Y) \stackrel{f_*}{\onto} N_1(X)
\end{align*} is of free abelian groups (of finite rank) therefore it splits (non-canonically), $N_1(Y) \cong N_1(X) \oplus N_1(Y/X)$.
Elements of $N_{\leq 1}(Y)$ can then be described by triples $(\gamma,\delta,n) \in N_1(X) \oplus N_1(Y/X) \oplus \Z$.
We denote the image of $\perv{p}{\acat{A}_{\leq 1}}$ (via the Chern character) in $N_{\leq 1}(Y)$ by $\perv{p}{\Delta}$ (this is the cone of perverse coherent sheaves supported in dimensions $\leq 1$).
The algebra $\text{H}(\perv{p}{\acat{A}_{\leq 1}})$ is graded by $\perv{p}{\Delta}$.
Finally, by $\mathscr{E} \subset N_1(Y/X)$ we denote the \emph{effective} curve classes in $Y$ which are contracted by $f$.
\begin{defn}\label{laurent-definition}
	We define a subset $L \subset \perv{p}{\Delta}$ to be \emph{Laurent} if the following conditions hold:
	\begin{itemize}
		\item for all $\gamma$ there exists an $n(\gamma,L)$ such that for all $\delta,n$ with $(\gamma,\delta,n) \in L$ one has that $n \geq n(\gamma,L)$;
		\item for all $\gamma, n$ there exists a $\delta(\gamma,n,L) \in \mathscr{E}$ such that for all $\delta$ with $(\gamma,\delta,n) \in L$ one has that $\delta \leq \delta(\gamma,n,L)$.\footnote{For $\delta,\delta' \in N_1(Y/X)$, by the notation $\delta \leq \delta'$ we mean $\delta' - \delta \in \mathscr{E}$ or equivalently $\delta - \delta' \in -\mathscr{E}$. In general we will write $\alpha \geq 0$ to denote that a certain class is effective.}
	\end{itemize}
	We denote by $\Lambda$ the set of all Laurent subsets of $\perv{p}{\Delta}$.
\end{defn}
Notice that $\Lambda$ does not depend on the choice of the above splitting.
We have the following lemma.
\begin{lem}
	The set $\Lambda$ of Laurent subsets of $\perv{p}{\Delta}$ satisfies the two following properties.
	\begin{enumerate}
		\item If $L_1, L_2 \in \Lambda$ then $L_1 + L_2 \in \Lambda$.
		\item If $\alpha \in \perv{p}{\Delta}$ and $L_1, L_2 \in \Lambda$ then there exist only finitely many decompositions $\alpha = \alpha_1 + \alpha_2$ with $\alpha_j \in L_j$.
	\end{enumerate}
\end{lem}
\begin{prf}
	We start by proving (1).
	Fix a $\gamma$ and let $(\gamma,\delta,n) \in L_1 + L_2$.
	By \cite[Corollary 1.19]{kollar} there are only finitely many decompositions $\gamma = \gamma_1 + \gamma_2$ with $\gamma_i \geq 0$ (i.e.~with $\gamma_i$ effective).
	Given a decomposition $(\gamma,\delta,n) = (\gamma_1 + \gamma_2, \delta_1 + \delta_2, n_1 + n_2)$, with $(\gamma_i,\delta_i,n_i) \in L_i$, we know that $n_i \geq n(\gamma_i,L_i)$ so $n = n_1 + n_2 \geq n(\gamma_1,L_1) + n(\gamma_2,L_2)$.
	By letting the $\gamma_i$'s vary we obtain the desired lower bound for $n$.
	
	Fix now $\gamma, n$.
	We want to find an upper bound for the possible $\delta$'s such that $(\gamma,\delta,n) \in L_1 + L_2$.
	By the argument above we know that for decompositions $(\gamma, \delta,n) = (\gamma_1 + \gamma_2, \delta_1 + \delta_2, n_1 + n_2)$ with $(\gamma_i,\delta_i,n_i) \in L_i$ the possible combinations of $\gamma_i$ and $n_i$ are finite.
	Fix such a decomposition $(\gamma_1+\gamma_2,\delta_1+\delta_2,n_1+n_2)$.
	We know that $\delta_i \leq \delta(\gamma_i,n_i,L_i)$.
	Thus $\delta = \delta_1 + \delta_2 \leq \delta(\gamma_1,n_1,L_1) + \delta(\gamma_2,n_2,L_2)$.
	Take now another decomposition $(\gamma_1' + \gamma_2', \delta_1' + \delta_2', n_1' + n_2')$.
	Running the same argument we have that $\delta \leq \delta(\gamma_1',n_1',L_1) + \delta(\gamma_2',n_2',L_2)$.
	Finally, as $\delta(\gamma_i,n_i,L_i), \delta(\gamma_i',n_i',L_i) \geq 0$, we conclude $\delta \leq \sum_i \delta(\gamma_i,n_i,L_i) + \delta(\gamma_i',n_i',L_i)$.
	By taking the sum for all possible decompositions we have our upper bound for $\delta$.
	
	Let us now prove (2).
	Fix a class $\alpha = (\gamma,\delta,n) \in \perv{p}{\Delta}$ and two Laurent subsets $L_1, L_2$.
	Again by \cite[Corollary 1.19]{kollar} we know that there are only finitely many possible decompositions $\gamma = \gamma_1 + \gamma_2$.
	Thus we may fix $\gamma_1$ and $\gamma_2$.
	Given a decomposition $(\gamma,\delta,n) = (\gamma_1 + \gamma_2,\delta_1 + \delta_2,n_1 + n_2)$, there are again finitely many possible values occurring for $n_1,n_2$, as $n_i \geq n(\gamma_i,L_i)$.
	Thus we may take $n_1,n_2$ also to be fixed.
	Finally, the combinations $(\gamma,\delta,n) = (\gamma_1+\gamma_2,\delta_1+\delta_2,n_1+n_2)$ are again a finite number, as $\delta = \delta_1 + \delta_2$ lives in $\delta(\gamma_1,n_1,L_1) + \delta(\gamma_2,n_2,L_2) -\mathscr{E}$ (thus we can apply \cite[Corollary 1.19]{kollar} again).
\end{prf}

We now have all the ingredients to define a $\Lambda$-completion $\text{H}(\perv{p}{\acat{A}_{\leq 1}})_\Lambda$ of $\text{H}(\perv{p}{\acat{A}}_{\leq 1})$.
Let us give a general definition.
\begin{defn}
	Let $R$ be a $\perv{p}{\Delta}$-graded associative $\Q$-algebra.
	We define $R_\Lambda$ to be the vector space of formal series
	\begin{align*}
		\sum_{(\gamma,\delta,n)} x_{(\gamma,\delta,n)}
	\end{align*}
	with $x_{(\gamma,\delta,n)} \in R_{(\gamma,\delta,n)}$ and $x_{(\gamma,\delta,n)} = 0$ outside a Laurent subset.
	We equip this vector space with a product
	\begin{align*}
		x \cdot y = \sum_{\alpha \in \perv{p}{\Delta}} \sum_{\alpha_1 + \alpha_2 = \alpha} x_{\alpha_1} \cdot y_{\alpha_2}.
	\end{align*}
	The algebra $R$ is included in $R_\Lambda$ as any finite set is Laurent.
	To a morphism $R \to S$ of $\perv{p}{\Delta}$-graded algebras corresponds an obvious morphism $R_\Lambda \to S_\Lambda$.
\end{defn}

There is a subalgebra 
\begin{align*}
	\Q_\sigma[\perv{p}{\Delta}] \subset \Q_\sigma[\perv{p}{\Gamma}]
\end{align*}
spanned by symbols $q^\alpha$ with $\alpha \in \perv{p}{\Delta}$.
Notice that the Poisson structure on $\Q_\sigma[\perv{p}{\Delta}]$ is trivial as the Euler form on $N_{\leq 1}(Y)$ is identically zero.
The integration morphism restricts to $I : \text{H}_\text{sc}(\perv{p}{\acat{A}_{\leq 1}}) \to \Q_\sigma[\perv{p}{\Delta}]$ and so, by taking $\Lambda$-completions, we have a morphism
\begin{align*}
	I_\Lambda: \text{H}_\text{sc}(\perv{p}{\acat{A}_{\leq 1}})_\Lambda \longto \Q_\sigma[\perv{p}{\Delta}]_\Lambda.
\end{align*}
\begin{rmk}
	Notice that given an algebra $R$ as above and an element $r \in R$ with $r_{(0,0,0)} = 0$, the element $1- r$ is invertible in $R_\Lambda$.
	This is due to the fact that the series
	\begin{align*}
		\sum_{k\geq 0} r^k
	\end{align*}
	makes sense in $R_\Lambda$.
\end{rmk}

Now it's time to have a look at what the elements of $\text{H}(\perv{p}{\acat{A}}_{\leq 1})_\Lambda$ look like.
Let $\astack{M}$ be an algebraic stack locally of finite type over $\C$ mapping down to $\perv{p}{\astack{A}}_{\leq 1}$ and denote by $\astack{M}_\alpha$ the preimage under $\perv{p}{\astack{A}}_\alpha$, for $\alpha \in \perv{p}{\Delta}$.
We say that
	\begin{align*}
		\left[\astack{M}\to \perv{p}{\astack{A}}_{\leq 1}\right] \in \text{H}_\infty(\perv{p}{\acat{A}}_{\leq 1})
	\end{align*}
is \emph{Laurent} if $\astack{M}_\alpha$ is a stack of finite type for all $\alpha \in \perv{p}{\Delta}$ and if $\astack{M}_\alpha$ is empty for $\alpha$ outside a Laurent subset.
Such a Laurent element gives an element of $\text{H}(\perv{p}{\acat{A}_{\leq 1}})_\Lambda$ by considering $\sum_\alpha \astack{M}_\alpha$.
\begin{prop}\label{prop-1flaurent}
	The elements $1_{\perv{p}{\acat{F}_{\leq 1}}[1]}$, $1^\O_{\perv{p}{\acat{F}_{\leq 1}}[1]}$ are Laurent.
\end{prop}
\begin{prf}
	Let $\sh{F} \in \perv{p}{\acat{F}_{\leq 1}}$ and let $(\gamma,\delta,n)$ be the class in $N_{\leq 1}(Y)$ corresponding to $[\sh{F}[1]] = -[\sh{F}]$.
	By \cite[Proposition X-1.1.2]{sga6} we know that in the Grothendieck group $\sh{F}$ decomposes as
	\begin{align*}
		\sh{F} = \sum_i l_i [\curlyO_{C_i}] + \tau
	\end{align*}
	where the $C_i$ are the curves comprising the irreducible components of the support of $\sh{F}$ (which is contained in the exceptional locus of $f$), where $l_i \geq 0$ and where $\tau$ is is supported in dimension zero.
	From this decomposition we infer that $\gamma = 0$ and $\delta \leq 0$.
	Finally, Riemann-Roch tells us that $n$ is minus the Euler characteristic of $\sh{F}$ and Lemma \ref{cohomology} gives us that $n \geq 0$.
	To conclude, the finite type axiom is deduced using Lemma \ref{finitetype}, combined with Lemma \ref{duality lemma}.
	
	For $1^\O_{\perv{p}{\acat{F}_{\leq 1}}[1]}$, it is enough to notice that for $F \in \perv{p}{\acat{F}}$, $H^1(Y,F)$ is finite-dimensional.
\end{prf}
Notice also that by the remark above both $1_{\perv{p}{\acat{F}_{\leq 1}}[1]}$ and $1^\O_{\perv{p}{\acat{F}_{\leq 1}}[1]}$ are invertible in $\text{H}(\perv{p}{\acat{A}_{\leq 1}})_\Lambda$.
\begin{prop}\label{prop-laurent-laurent}
	The element $\perv{p}{\curly{H}}_{\leq 1}$ is Laurent.
\end{prop}
\begin{prf}
	By \cite[Theorem 7.3]{tom-flops} if we fix a numerical class $\alpha \in N_{\leq 1}(Y)$ then the space $\perv{p}{\Hilb}_{Y/X}(\alpha)$ is of finite type (it is in fact a projective scheme).
	Thus we are left with checking the second half of the Laurent property.
	Fix then a class $\gamma \in N_1(X)$ and consider a possible quotient $\O_Y \onto \csh{P}$ in $\perv{p}{\acat{A}}$, with $\dim \supp \csh{P} \leq 1$ and with $\csh{P}$ of class $(\gamma,\delta,n)$.
	We need to show that there exists a lower bound on the possible values of $n$.
	By pushing down to $X$ we obtain a quotient (in $\Coh (X)$) $\O_X \onto Rf_* \csh{P}$, and we note that the sheaf $Rf_* \csh{P}$ is of class $(\gamma,n)$.
	If a class $\gamma$ is fixed, it is known that the possible values of the Euler characteristic of a quotient $\O_X \onto Q$ are bounded below (this follows from boundedness of the Hilbert scheme), hence we have the required bound.
	
	To proceed, we let $\gamma$ and $n$ both be fixed and notice that we only really need to focus on exact sequences of both coherent and perverse sheaves, that is on points of $\Hilb_Y \cap \perv{p}{\Hilb}_{\leq 1}(Y/X)$ (which we temporarily denote by $\Pilb(Y)$).
	This is a consequence of the fact that given an epimorphism $\O_Y \onto \csh{P}$, with $\csh{P}$ in $\perv{p}{\acat{A}_{\leq 1}}$ of class $(\gamma,\delta,n)$, we can consider the torsion torsion-free exact sequence
	\begin{align*}
		\sh{F}[1] \into \csh{P} \onto \sh{T}.
	\end{align*}
	In fact, $\sh{F}[1]$ does not contribute towards $\gamma$, contributes negatively towards $\delta$ and positively towards $n$, as seen in the previous proposition.
	Thus we just need to study the possible classes of $\sh{T}$.
	Finally, $\O_Y \onto \csh{P} \onto \sh{T}$ is a quotient in $\perv{p}{\acat{A}}$ but glancing at the cohomology sheaves long exact sequence tells us that it is indeed a quotient in $\acat{A}$ as well.
	Thus we only need to check that, having chosen a $\gamma$ and an $n$, there exists an upper bound $\delta_0$ such that $\Pilb_{Y}(\gamma,\delta,n)$ is empty for $\delta \geq \delta_0$.
	
	Notice that the pushforward induces a morphism from $\perv{p}{\Hilb}({Y/X})$ to $\Hilb(X)$.
	We consider its restriction to $\Pilb(Y)$.
	We would like for the pullback functor to induce a morphism going in the opposite direction.
	A flat family of sheaves on $X$ might, however, cease to be flat once pulled back on $Y$.
	To remedy we impose this condition by hand.
	We define a subfunctor $\Filb_X$ of $\Hilb_X$ by the rule
	\begin{align*}
		\Filb_X(S) = \left\{\curlyO_{X_S} \onto \sh{G} \st \sh{G}, f_S^* \sh{G} \text{ flat over } S \right\}.
	\end{align*}
	If $\sh{U}$ is the structure sheaf of the universal subscheme for $\Hilb(X)$ on $X \times \Hilb(X)$ then one can see that $\Filb(X)$ is represented by the flattening stratification of $\Hilb(X)$ with respect to $f_{\Hilb(X)}^* \sh{U}$.
	From this we deduce that if we fix a numerical class $(\gamma,n)$ on $X$ then $\Filb_{X}(\gamma,n)$ is of finite type.
	
	We claim that the composition of pushing forward and pulling up as just described, $\Pilb(Y) \to \Filb(X) \to \Pilb(Y)$, is the identity.
	Let us see first why this is true on geometric points.
	Take an exact sequence of both coherent and perverse sheaves
	\begin{align*}
		\sh{I} \into \curlyO_Y \onto \sh{E}.
	\end{align*}
	Applying the counit of the adjunction $f^* \dashv f_*$ (and using the fact that the objects above are both sheaves and perverse sheaves) we obtain a commutative diagram
	\begin{center}
		\begin{tikzpicture}
			\matrix (m) [matrix of math nodes, row sep=3em, column sep=3em, text height=1.5ex, text depth=0.25ex]
			{
			 & f^*f_*\sh{I} & \curlyO_Y & f^*f_*\sh{E} & 0 \\
			0 & \sh{I} & \curlyO_Y & \sh{E} & 0 \\
			};
			\path[->,font=\scriptsize]
			(m-1-2) edge (m-1-3)
					edge (m-2-2)
			(m-1-3) edge (m-1-4)
					edge node[auto]{$\id$} (m-2-3)
			(m-1-4) edge (m-1-5)
					edge (m-2-4)
			(m-2-1) edge (m-2-2)
			(m-2-2) edge (m-2-3)
			(m-2-3) edge (m-2-4)
			(m-2-4) edge (m-2-5)
			;
		\end{tikzpicture}
	\end{center}
	with exact rows.
	By \cite[Proposition 5.1]{tom-flops} we have that $f^*f_*\sh{I} \to \sh{I}$ is surjective and so, by a simple diagram chase, $f^*f_*\sh{E} \to \sh{E}$ is an isomorphism.
	This argument indeed works in families, as surjectivity can be checked fibrewise.
	
	Finally, let us fix a $\gamma$ and an $n$ and let $\Pilb_Y(\gamma,n)$ be the subspace of $\Pilb(Y)$ where we've fixed $\gamma$ and $n$ but we let $\delta$ vary.
	By the previous arguments we know that $\Pilb_Y(\gamma,n) \to \Filb_X(\gamma,n) \to \Pilb_Y(\gamma,n)$ composes to the identity.
	As the retract of a quasi-compact space is quasi-compact\footnote{If $A \to B \to A$ composes to the identity, one can start with an open cover $\{A_i\}$ and pull it back to a cover $\{B_i\}$ of $B$. Pick a finite subcover $\{B_j\}$ and pull it back to $A$. This is a finite subcover of $\{A_i\}$.} we obtain that $\Pilb_Y(\gamma,n)$ is of finite type, which is enough to conclude.
\end{prf}

\begin{prop}\label{laurenthilb}
	The element $\curly{K}_{\leq 1}$ is Laurent.
\end{prop} 
\begin{prf}
	It is a known fact that for a fixed numerical class $\alpha \in N_{\leq 1}(Y)$ the scheme $\Hilb_Y(\alpha)$ is of finite type (it is in fact a projective scheme).
	To prove the second half of the Laurent property we start by fixing a class $\gamma \in N_1(X)$.
	If $\O_Y \onto \sh{T}$ is a quotient in $\acat{A}$ with kernel $\sh{I}$, we have an exact sequence
	\begin{align*}
		0 \to f_* \sh{I} \to \O_X \to f_*\sh{T} \to R^1f_* \sh{I} \to 0.
	\end{align*}
	If $\sh{T}$ is of class $(\gamma,\delta,n)$ then $f_* \sh{T}$ is of class $(\gamma,n)$ and $R^1 f_* \sh{I}$ is supported in dimension zero.
	The image $\sh{Q}$ of $\O_X \to f_* \sh{T}$ is of class $(\gamma,m)$ with $m \leq n$.
	As $\gamma$ is fixed we have a lower bound on the possible values of $m$ and a fortiori on the values of $n$.
	
	Let now $\gamma$ and $n$ be fixed.
	We start off with the identity
	\begin{align*}
		\perv{p}{\curly{H}}_{\leq 1} * 1_{\perv{p}{\acat{F}_{\leq 1}[1]}} = 1^\curlyO_{\perv{p}{\acat{F}_{\leq 1}[1]}} * \curly{K}_{\leq 1}
	\end{align*}
	in $\text{H}_\infty(\perv{p}{\acat{A}_{\leq 1}})$.
	By directly applying our definition of $*$ we see that the right hand side is represented by a morphism $[W \to \perv{p}{\astack{A}_{\leq 1}}]$, given by the top row of the following diagram.
	\begin{center}
		\begin{tikzpicture}
			\matrix (m) [matrix of math nodes, row sep=3em, column sep=3em, text height=1.5ex, text depth=0.25ex]
			{
			W& \perv{p}{\astack{A}_{\leq 1}^{(2)}} & \perv{p}{\astack{A}_{\leq 1}}\\
			\perv{p}{\astack{F}_{\leq 1}[1]^\O} \times \Hilb_{Y,\leq 1}^\partial & \perv{p}{\astack{A}_{\leq 1}}\times  \perv{p}{\astack{A}_{\leq 1}} &\\
			};
			\path[->,font=\scriptsize]
			(m-1-1) edge (m-2-1)
					edge (m-1-2)
			(m-1-2) edge node[auto]{$b$} (m-1-3)
					edge node[auto]{$(a_1,a_2)$} (m-2-2)
			(m-2-1) edge (m-2-2)
			;
		\end{tikzpicture}
	\end{center}
	Recall Remark \ref{addedremark} for the definition of $\Hilb^\partial$.
	Similarly, the left hand side is represented by a morphism $[Z \to \perv{p}{\astack{A}_{\leq 1}}]$.
	The main tool we use for the proof is the cover $\{\perv{p}{\astack{A}_\alpha}\}_\alpha$ of $\perv{p}{\astack{A}_{\leq 1}}$, with $\alpha \in \perv{p}{\Delta}$ ranging inside the cone of perverse coherent sheaves.
	
	By taking preimages through $b$ we obtain an open cover $\{U_\alpha\}_\alpha$ of $\perv{p}{\astack{A}^{(2)}_{\leq 1}}$.
	Concretely, $U_\alpha$ parameterises exact sequences $\csh{P}_1 \into \csh{P} \onto \csh{P}_2$ in $\perv{p}{\acat{A}_{\leq 1}}$ with $\csh{P}$ of class $\alpha$.
	
	On the other hand, we can cover $\perv{p}{\astack{A}_{\leq 1}} \times \perv{p}{\astack{A}_{\leq 1}}$ by taking products $\perv{p}{\astack{A}_{\alpha_1}} \times \perv{p}{\astack{A}_{\alpha_2}}$.
	By pulling back via $(a_1,a_2)$ we produce an open cover $\{U_{\alpha_1,\alpha_2}\}_{\alpha_1,\alpha_2}$ of $\perv{p}{\astack{A}_{\leq 1}^{(2)}}$.
	The space $U_{\alpha_1,\alpha_2}$ parameterises exact sequences $\csh{P}_1 \into \csh{P} \onto \csh{P}_2$ in $\perv{p}{\acat{A}_{\leq 1}}$ with $\csh{P}_1$ of class $\alpha_1$ and $\csh{P}_2$ of class $\alpha_2$.
	Notice that the collection $\{U_{\alpha_1,\alpha_2}\}_{\alpha_1 + \alpha_2 = \alpha}$ is an open cover of $U_\alpha$.
	
	By pulling back these covers of $\perv{p}{\astack{A}_{\leq 1}^{(2)}}$ we obtain open covers $\{W_\alpha\}_\alpha$ and $\{W_{\alpha_1,\alpha_2}\}_{\alpha_1,\alpha_2}$ of $W$.
	The same can be done for $Z$.
	
	We remind ourselves that we think of a class $\alpha$ as a triple $(\gamma,\delta,n)$.
	If we fix a $\gamma$ and an $n$, it is a consequence of $\perv{p}{\curly{H}}_{\leq 1} * 1_{\perv{p}{\acat{F}_{\leq 1}[1]}}$ being Laurent that there exists a $\delta'$ such that $Z_{(\gamma,\delta,n)} = \emptyset$ for $\delta \geq \delta'$.
	Because of the identity above, the same holds for $W_{(\gamma,\delta,n)}$.
	
	What we need to prove is that, once we fix $\gamma$ and $n_2$, the space $\Hilb_Y(\gamma,\delta_2,n_2)$ is empty for large $\delta_2$.
	Fix $\delta_1, n_1$ such that $\perv{p}{\astack{F}_{\leq 1}[1]}^\O_{(0,\delta_1,n_1)} \neq \emptyset$.
	The space representing the product 
	$$1^\O_{\perv{p}{\astack{F}_{\leq 1}[1]_{(0,\delta_1,n_1)}}} * \Hilb^\partial_Y(\gamma,\delta_2,n_2)$$ 
	is $W_{(0,\delta_1,n_1),(\gamma,\delta_2,n_2)} \subset W_{(\gamma,\delta_1 + \delta_2, n_1 +n_2)}.$
	We have already remarked that for fixed $\gamma$, $n_1,n_2$ we have an upper bound $\delta'$ such that $W_{(\gamma,\delta_1 + \delta_2,n_1+n_2)} = \emptyset$ for $\delta_1 + \delta_2 \geq \delta'$.
	As $\perv{p}{\astack{F}_{\leq 1}[1]}^\O_{(0,\delta_1,n_1)} \neq \emptyset$, we conclude that $\Hilb^\partial_Y(\gamma,\delta_2,n_2) = \emptyset$ for $\delta_2 \geq \delta' - \delta_1$, in particular the same is true for $\delta_2 \geq 0$.
\end{prf}

\begin{rmk}\label{dualityandpt}
	We need to interpret Proposition \ref{dualityprop} in the Laurent setting.
	Duality $\D'$ acts on $N_{\leq 1}(Y)$ by taking a class $(\gamma,\delta,n)$ to $(-\gamma,-\delta,n)$.
	Even more concretely, an element $\sh{T} \in \perv{q}{\acat{T}_\bullet}$ of class $(0,\delta,n)$ is sent to an element $\D'(\sh{T}) \in \perv{p}{\acat{F}_{\derp}}$ of class $(0,-\delta,n)$.
	This suggests that we should complete the algebra $H(\perv{q}{\acat{T}_\bullet})$ with respect to a sort of dual Laurent subsets.
	
	Let $\perv{q}{\Delta}_{\exc}$ be the subcone of $\perv{q}{\Delta}$ consisting of elements of the form $(0,\delta,n)$, with $n \geq 0$.
	We define $\Lambda'$ as the collection of subsets $L \subset \perv{q}{\Delta}_{\exc}$ such that:
	\begin{itemize}
		\item for all $n \in \Z$, there exists $\delta(n) \in \curly{E}$ such that for all $\delta$, with $(0,\delta,n) \in L$, $\delta \geq \delta(n)$.
	\end{itemize}
	We can complete the algebra $H(\perv{q}{\acat{T}_\bullet})$ with respect to $\Lambda'$, just as we complete the Hall algebra of perverse coherent sheaves with respect to $\Lambda.$
	We denote this completion by $\text{H}(\perv{q}{\acat{T}_\bullet})_{\Lambda'}$.
	
	The elements $1_{\perv{q}{\acat{T}_\bullet}}$ and $1_{\perv{q}{\acat{T}_\bullet}}^\curlyO$ belong $\text{H}(\perv{q}{\acat{T}_\bullet})_{\Lambda'}$ by Proposition \ref{prop-1flaurent} and duality.
	The element $\curly{H}^\#_{\exc}$ also belongs to $\text{H}(\perv{q}{\acat{T}_\bullet})_{\Lambda'}$ by running a similar proof to the one above, using \eqref{1ot=hs*1t}.
	Proposition \ref{dualityprop} now implies $\D'$ defines an isomorphism between $H(\perv{q}{\acat{T}_\bullet})_{\Lambda'}$ and $H(\perv{p}{\acat{F}_{\leq 1}[1]})_\Lambda$, taking $1_{\perv{q}{\acat{T}_\bullet}}$ to $1_{\perv{p}{\acat{F}_{\leq 1}[1]}}$ and $1^\O_{\perv{q}{\acat{T}_\bullet}}$ to $1^\O_{\perv{p}{\acat{F}_{\leq 1}[1]}}$.
	%
	%
	%
	%
	%
	%
\end{rmk}

Going back to $\text{H}(\perv{p}{\acat{A}_{\leq 1}})_\Lambda$, the remark above implies the identity
\begin{align}\label{identity-leq1}
	\perv{p}{\curly{H}}_{\leq 1} = 1_{\perv{p}{\acat{F}_{\derp}[1]}} * \D'(\curly{H}^\#_{\exc}) * \curly{K}_{\leq 1} * 1_{\perv{p}{\acat{F}_{\leq 1}[1]}}^{-1}.
\end{align}
What keeps us from simply applying the integration morphism $I_\Lambda$ is that, although $\D'(\curly{H}_{\exc}^\#)$ and $\curly{K}_{\leq 1}$ are regular (in the sense of Proposition \ref{semiclassicalpminusone}), $1_{\perv{p}{\acat{F}_{\derp}[1]}}$ is not.
\subsection{A Stability Condition} 
\label{sub:a_stability_condition}
We want to proceed analogously as in \cite[Section 6.3]{tom-cc}, proving that 
\begin{align*}
	I_\Lambda(\perv{p}{\curly{H}_{\leq 1}}) = I_\Lambda(\D'(\curly{H}_{\exc}^\#)) \cdot I_\Lambda(\curly{K}_{\leq 1})
\end{align*}
holds nevertheless.
We make use of an important result of Joyce, which we can roughly summarise as follows.
Suppose we are working in the Hall algebra of an abelian category and suppose we are given a stability condition.
The slogan we keep in mind is: ``the product $[\C^\times]\cdot \log(1_{\text{SS}(\mu)})$, where $1_{\text{SS}(\mu)}$ is the element corresponding to the inclusion of semi-stable objects of slope $\mu$, is a regular element.''
It will suffice to combine \cite[Corollary 5.10]{joyce-configii} and \cite[Theorem 8.7]{joyce-configiii}.

In our context the key is to show that $(\L - 1) \cdot \log (1_{\perv{p}{\acat{F}_{\derp}[1]}}) \in H_\text{reg}(\perv{p}{\acat{A}})$.
This can be achieved by constructing an appropriate stability condition such that $\perv{p}{\acat{F}_{\derp}[1]}$ manifests as the set of objects of some fixed slope.
For convenience we work within the category $\perv{p}{\acat{A}_{\exc}}$, whose objects are those perverse coherent sheaves $\csh{P} \in \perv{p}{\acat{A}_{\leq 1}}$ whose pushforward to $X$ is supported on points (in other words such a $\csh{P}$ is of class $(0,\delta,n)$, for some $\delta \in N_1(Y/X)$ and $n \in \Z$ -- see the next section).
We define a stability condition $\mu$, taking values in the ordered set $\{1,2\}$ as follows.
\begin{align*}
	(0,\delta,n) \longmapsto 
		\begin{cases} 
			1 \text{ if } \delta \geq 0 \\
			2 \text{ if } \delta < 0. \\
		\end{cases}
\end{align*}
It is immediate that $\mu$ is indeed a weak stability condition (in the sense of \cite[Definition 3.5]{joyce-book}), as the only axiom one needs to check is the (weak) see-saw property.
\begin{lem}\label{sstob}
	The set of $\mu$-semistable objects of slope $\mu = 2$ is $\perv{p}{\acat{F}_{\derp}[1]}$.
	The set of $\mu$-semistable objects with $\mu = 1$ is $\perv{p}{\acat{T}_{\exc}}$.
\end{lem}
Recall that an object $\csh{P}$ is said to be \emph{semistable} if for all proper subobjects $\csh{P'} \subset \csh{P}$ we have $\mu(\csh{P'}) \leq \mu(\csh{P}/\csh{P'})$.
\begin{prf}
	Let $\csh{P}$ be any semistable perverse coherent sheaf.
	Consider the torsion torsion-free exact sequence
	\begin{align*}
		\sh{F}[1] \into \csh{P} \onto \sh{T}.
	\end{align*}
	If $\sh{F}[1] \neq 0$ and $\sh{T} \neq 0$ then, by semistability, $2 = \mu(\sh{F}[1]) \leq \mu(\sh{T}) = 1$ which is impossible.
	Thus a semistable object must be either torsion or torsion-free.
	
	On the other hand, as $\perv{p}{\acat{F}_{\derp}[1]}$ is stable under quotients and $\perv{p}{\acat{T}_{\exc}}$ is stable under subobjects we conclude.
\end{prf}
The last property we need is permissibility, in the sense of \cite[Definition 4.7]{joyce-configiii}.
\begin{prop}
	The stability condition $\mu$ is permissible.
\end{prop}
\begin{prf}
	The first fact we check is that the category $\perv{p}{\acat{A}_{\exc}}$ is noetherian.
	More generally, this follows from Noetherianness of $\perv{p}{\acat{A}}$.
	The latter can be seen as a consequence of \cite{vdb}, as $\perv{p}{\acat{A}}$ is equivalent to the category of finitely generated modules over a noetherian coherent $\O_X$-algebra.
	
	Now we want to check that if $\csh{P} \in \perv{p}{\acat{A}_{\exc}}$ and $[\csh{P}]=0$ in $N_{\leq 1}(Y)$ then $\csh{P}=0$.
	By pushing forward via $f$ we have that $[Rf_*\csh{P}]=0$ and as $Rf_*\csh{P} \in \Coh (X)$ it follows that $Rf_*\csh{P} = 0$.
	Now, from Leray's spectral sequence we obtain that $f_*H^{-1}(\csh{P}) = f_*H^0(\csh{P}) = 0$.
	Thus, if $p=-1$, $\csh{P} = H^{-1}(\csh{P})[1]$ and, if $p=0$, $\csh{P} = H^0(\csh{P})$.
	In either case we reduce to dealing with a coherent sheaf and so $\csh{P}=0$.
	
	Let now $\perv{p}{\acat{A}_\alpha(i)}$ be the subset of $\perv{p}{\astack{A}_{\exc}(\C)}$ consisting of perverse coherent sheaves which are of numerical class $\alpha$ and semistable with $\mu = i$.
	We now check that these subsets are constructible.
	
	In light of Lemma \ref{duality lemma} and Lemma \ref{sstob}, what remains to be proved is that, given a $\delta$ and an $n$, the stack $\perv{p}{\astack{T}_{(0,\delta,n)}}$ is of finite type.
	This is the content of the following lemma.
	
	To finish, we show that $\mu$ is artinian.
	Consider a chain of subobjects
	\begin{align*}
		\cdots \into \csh{P}_2 \into \csh{P}_1
	\end{align*}
	with $\mu(\csh{P}_{n+1}) \geq \mu(\csh{P}_n /\csh{P}_{n+1})$.
	Let $\csh{P'} \into \csh{P}$ be any two consecutive elements in the chain above and let $\csh{Q}$ be the quotient $\csh{P'}/\csh{P}$ so that we have an exact sequence
	\begin{align*}
		\csh{P'} \into \csh{P} \onto \csh{Q}
	\end{align*}
	with $\mu(\csh{P'}) \geq \mu(\csh{Q})$, which corresponds to the relation $(\delta',n') + (\delta_q,n_q) = (\delta,n)$ in $\perv{p}{\Delta}$.
	As the sheaves we are considering have support contracted by $f$ we see that $n,n',n_q \geq 0$ hence we can assume (by going further down the chain if necessary) that $n = n'$, which in turn implies $n_q = 0$.
	
	When $p = -1$, this implies that $Q \in \perv{p}{\acat{F}_{\leq 1}[1]}$ and $\mu(\csh{Q}) = 2$.
	As a consequence, we have $\mu(\csh{P'})=\mu(\csh{P}) = 2$ and so $\delta,\delta' <0$.
	Finally, as $\delta' \geq \delta$, we can assume $\delta' = \delta$ and so $\delta_q = 0$, from which we gather that $\csh{Q} = 0$, which concludes the proof.
	
	When $p= 0$, the condition $n_q = 0$ implies $Q \in \perv{p}{\acat{T}_\text{exc}}$.
	Glancing at the cohomology sheaves long exact sequence we see that $P' \to P$ is an isomorphism on $H^{-1}$ and an injection on $H^0$.
	If we denote $\delta_0 = \ch_2 (H^{-1}(P))$, we see that $\delta_0 \leq \delta' \leq \delta$.
	Thus, again by descending further down the chain if necessary, we can assume $\delta = \delta'$ and we are done.
	%
	%
	%
	%
	%
	%
	%
	%
	%
\end{prf}

\begin{lem}\label{finitetype}
	Let $\delta \in N_1(Y/X)$ and let $n \in \Z$.
	Then, the stack $\perv{p}{\astack{T}_{(0,\delta,n)}}$ is of finite type.
\end{lem}
\begin{prf}
	We can use the criterion found for example in \cite[Lemma 1.7.6]{shaves}.
	Let $S$ be a finite type scheme and let $T \in \Coh(S\times Y)$ be a flat family of sheaves such that for any closed point $s \in S$ the restriction to the fibre $T_s$ lies in $\perv{p}{\acat{T}_{\exc}}$ and is of class $(0,\delta,n)$.
	We show that there exists a sheaf surjecting onto all the $T_s$.
	It is sufficient to prove that $T_s$ is generated by global sections, as then the sheaf $H^0(Y,T_s) \otimes_\C \O_Y$ will surject onto $T_s$ and $\dim H^0(Y,T_s) = n$ is independent of $s$.
	
	Let $I$ and $C$ be respectively the image and the cokernel of the evaluation morphism $H^0(Y,T_s)\otimes_\C \O_Y \to T_s$.
	The sheaf $I$ belongs to $\perv{p}{\acat{T}}$ and together with the exact sequence $I \into T_s \onto C$ we have
	\begin{align*}
		f_*I \into f_* T_s \onto f_*C.
	\end{align*}
	As $f_*T_s$ is supported on points, the morphism $H^0(X,f_*T_s) \otimes_\C \O_X \to f_*T_s$ is surjective, which (by adjunction) implies that $f_*I \to f_*T_s$ is surjective, which by the above exact sequence implies that $f_*C=0$.
	As $R^1f_* C = 0$ as well, by the properties of perverse coherent sheaves it follows that $T_s \to C$ is the zero morphism, which in turn implies $C = 0$.
	Hence the claim.
\end{prf}
\begin{prop}\label{adexp}
	In $\text{H}(\perv{p}{\curvy{A}}_{\leq 1})_\Lambda$ we have $1_{\perv{p}{\acat{F}_{\derp}}[1]} = \exp (\epsilon)$, with $\eta = (\L - 1) \cdot \epsilon \in \text{H}_\text{reg}(\perv{p}{\acat{A}}_{\leq 1})_\Lambda$ a regular element.
	Furthermore the automorphism
	\begin{align*}
		\Ad_{1_{{\perv{p}{\acat{F}}_{\derp}[1]}}}: \text{H}(\perv{p}{\curvy{A}}_{\leq 1})_\Lambda \longto \text{H}(\perv{p}{\acat{A}}_{\leq 1})_\Lambda
	\end{align*}
	preserves regular elements.
	The induced Poisson automorphism of $\text{H}_\text{sc}(\perv{p}{\acat{A}_{\leq 1}})_\Lambda$ is given by
	\begin{align*}
		\Ad_{1_{\perv{p}{\acat{F}}_{\derp}[1]}} = \exp \{ \eta,- \}.
	\end{align*}
\end{prop}
\begin{prf}
	We can draw an argument entirely parallel to the one in Theorem 6.3 and Corollary 6.4 of \cite{tom-cc}.
	The only thing to check here is that the class $[\C^*] \cdot \log ( 1_{\perv{p}{\acat{F}}_{\derp}[1]} )$ is a regular element, which can be done in the subalgebra $H(\perv{p}{\acat{A}_{\text{exc}}})_\Lambda$, exploiting the stability condition we just described.
\end{prf}

\subsection{Main Identity} 
\label{sub:conclusion}
At last, we have all the ingredients to prove our main result.
Before we proceed, we must deal with an issue of signs.

The Hilbert scheme $\Hilb_Y(\beta,n)$ comes with two constructible functions which are of interest to us.
The first ($\nu$) is Behrend's microlocal function.
The second ($\mu$) is the pullback along $\Hilb_Y(\beta,n) \to \astack{A}$ of the Behrend function of the stack $\astack{A}$.
Given a homology class $\beta$ and an integer $n$, the DT number of class $(\beta,n)$ is defined to be
\begin{align*}
	\DT_Y(\beta,n) := \chi_{\nu}\left(\Hilb_Y(\beta,n)\right) := \sum_{k \in \Z} k \chi_{\topp}(\nu^{-1}(k))
\end{align*}
where $\chi_{\topp}$ is the topological Euler characteristic.
We package all these numbers into a generating series
\begin{align*}
	\DT(Y) = \sum_{\beta,n} \DT_Y(\beta,n) q^{(\beta,n)}
\end{align*}
which can be interpreted as a Laurent series according to our definition.
As we work with the integration morphism, it is convenient for us to define a variant of the DT series:
\begin{align*}
	\correct{\DT}(Y) \coloneqq I_\Lambda(\curly{H}_{\leq 1}).
\end{align*}
Unpacking the definitions, we see that, if we write
\begin{align*}
	\correct{\DT}_Y(\beta,n) := \chi_{\mu}\left(\Hilb_Y(\beta,n)\right) := \sum_{k \in \Z} k \chi_{\topp}(\mu^{-1}(k))
\end{align*}
then
\begin{align*}
	I_\Lambda\left(\curly{H}_{\leq 1}\right) = \correct{\DT}(Y) = \sum_{\beta,n} \correct{\DT}_Y(\beta,n) q^{(\beta,n)}.
\end{align*}
\begin{rmk}\label{flops signs remark}
	It is shown in \cite[Theorem 3.1]{tom-cc} that there is a simple relationship between $\DT$ and $\correct{\DT}$, given as follows.
	\begin{align*}
		\correct{\DT}_Y(\beta,n) = (-1)^n \DT_Y(\beta,n)
	\end{align*}
\end{rmk}

We pause a moment to notice that on $\Hilb(Y)$, not only do we have the pullback of the Behrend function of $\astack{A}$, but also the pullback of the Behrend function of $\perv{p}{\astack{A}}$.
However, no ambiguity arises, as the two morphisms factor through $\perv{p}{\astack{T}}$, which is open in both $\astack{A}$ and $\perv{p}{\astack{A}}$.

To deal with the general setup of Situation \ref{situation}, we had to work with the moduli spaces $\Hilb^\partial_Y(\beta,n) \subset \Hilb^\partial_Y(\beta,n)$ (see Remark \ref{addedremark}), corresponding to elements in the Hall algebra $\curly{K}_{(\beta,n)}$ -- whose sum is the Laurent element $\curly{K}_{\leq 1}$.
In turn this leads to ``partial'' DT invariants 
$$ \correct{DT}^\partial_Y(\beta,n) 
	:= I(\curly{K}_{(\beta,n)}) = \chi_\mu \left( \Hilb^\partial_Y(\beta,n) \right) 
	= (-1)^n \chi_\nu \left( \Hilb^\partial_{Y}(\beta,n) \right)
	=: (-1)^n DT^\partial_Y(\beta,n).
$$
Recall that when $X$ has zero-dimensional singular locus then $\Hilb^\partial = \Hilb$ and thus $DT^\partial = DT$.

Analogously as above, we have the perverse Hilbert scheme $\perv{p}{\Hilb_{Y/X}}(\beta,n)$.
Its weighted Euler characteristic produces perverse DT invariants
$$ \perv{p}DT_{Y/X}(\beta,n) := \chi_\mu \left( \perv{p}\Hilb_{Y/X}(\beta,n) \right). $$
However, to deal with the case where $X$ has one-dimensional singular locus, we must replace $\perv{p}{\Hilb_{Y/X}(\beta,n)}$ by $\perv{p}{\Hilb_{Y/X,\leq 1}}(\beta,n) = \perv{p}{\Hilb_{Y/X}^\partial}(\beta,n)$, which is the open subspace corresponding to epimorphisms $\O_Y \to E$, with $E$ supported in dimension one.
We can sum all the invariants together obtaining
$$ \perv{p}{\correct{DT}}^\partial(Y/X):= I_\Lambda (\perv{p}{\curly{H}_{\leq 1}}) = \sum_{\beta,n} \perv{p}{\correct{DT}}^\partial_{Y/X}(\beta,n) q^{(\beta,n)}.$$
As usual, when $X$ has singular locus of dimension zero then the ``partial'' invariants are the same as the ordinary ones.

Fortunately, for the purposes of this paper, we needn't be concerned with comparing $\chi_\mu(\perv{p}{\Hilb_{Y/X}(\beta,n)})$ with $\chi_\nu(\perv{p}{\Hilb_{Y/X}(\beta,n)})$.

\begin{rmk}\label{techfinal}
	We point out that the identities we write down below should be interpreted as taking place in the algebra $\Q_\sigma[\perv{p}{\Delta}]_\Lambda$, defined in Subsection \ref{sub:laurent_elements}.
\end{rmk}

We introduce the following sums,
\begin{align*}
	{\DT}_0(Y) &:= \sum_{n \in \Z}{\DT}_Y(0,n) q^{(0,n)}\\
	{\DT}_{\exc}(Y) &:= \sum_{\substack{\beta, n \\ f_* \beta = 0}} {\DT}_Y(\beta,n) q^{(\beta,n)}\\
	{\DT}^\vee_{\exc}(Y) &:= \sum_{\substack{\beta,n \\ f_* \beta = 0}} {\DT}_Y(-\beta,n) q^{(\beta,n)}
\end{align*}
and their $\correct{\DT}$ analogues.
\begin{teo}\label{mainmanor}
	Assume to be working in Situation \ref{situation}.
	The following identity holds.
	\begin{align}
		\perv{p}{\correct{\DT}^\partial(Y/X)} = \frac{\correct{\DT}_{\exc}^\vee(Y) \cdot \correct{\DT}^\partial(Y)}{\DT_0(Y)}
	\end{align}
\end{teo}
\begin{prf}
	The Poisson bracket on $\Q_\sigma[\perv{p}{\Delta}]$ is trivial, so Proposition \ref{adexp}, together with \eqref{identity-leq1}, yields the identity
	\begin{align*}
		I_\Lambda(\perv{p}{\curly{H}_{\leq 1}}) = I_\Lambda(\D'(\curly{H}_{\exc}^\#)) \cdot I_\Lambda(\curly{K}_{\leq 1}).
	\end{align*}
	The left hand side is equal to $\perv{p}{\correct{\DT}^\partial(Y/X)}$ and $I_\Lambda(\curly{K}_{\leq 1}) = \correct{\DT}^\partial(Y)$.
	As remarked in Subsection \ref{sub:a_route}, \cite[Lemma 5.5 and Theorem 1.1]{tom-cc} tell us how $\curly{H}^\#$ is related to DT invariants.
	In fact, combining these with Remark \ref{dualityandpt} we see that
	\begin{align*}
		I_\Lambda(\D'(\curly{H}_{\exc}^\#)) = \frac{\correct{\DT}_{\exc}^\vee(Y)}{\DT_0(Y)}
	\end{align*}
	and hence the claim.
\end{prf}
Notice that, as we are working with the assumption of Remark \ref{crucial perversity remark flops}, the Theorem above holds for both perversities, hence the series $\perv{p}{\correct{\DT}^\partial(Y/X)}$ is independent of the perversity $p$.
We will therefore drop the superscript $p$.
\subsection{Conclusion}
Now that we understand how the category of perverse coherent sheaves relates to DT invariants we can prove our promised formula for flops.
\begin{situation}\label{situflop2}
	Recall Situation \ref{situation} and assume moreover $X$ to have zero-dimensional singular locus (in other words $f\colon Y \to X$ is an isomorphism in codimension one).
	Let $f^+ \colon Y^+ \to X$ be the flop of $f$.
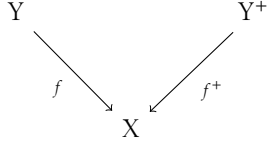
	\begin{center}
		\begin{tikzpicture}
			\matrix (m) [matrix of math nodes, row sep=3em, column sep=3em, text height=1.5ex, text depth=0.25ex]
			{
			Y && Y^+\\
			&X&\\
			};
			\path[->,font=\scriptsize]
			(m-1-1) edge node[auto,swap]{$f$} (m-2-2)
			(m-1-3) edge node[auto]{$f^+$} (m-2-2)
			;
		\end{tikzpicture}
	\end{center}	
\end{situation}
Notice that with these additional assumptions it follows automatically that $\perv{p}{\acat{F}} = \perv{p}{\acat{F}_{\derp}}$ (for $p = -1,0$).
In turn, all ``partial'' superscripts become superfluous and $DT^\partial = DT$. 

Following \cite{tom-flops}, we know that the variety $Y^+$ can be constructed as the moduli space of point-like objects of $\perv{-1}{\Per(Y/X)} = \perv{-1}{\acat{A}}$, the category of perverse coherent sheaves with minus one perversity.
The pair $(Y^+,f^+)$ satisfies the same assumptions as $(Y,f)$, so the categories of perverse coherent sheaves $\perv{q}{\Per(Y^+/X)} = \perv{q}{\acat{A}^+}$ (for $q = -1,0$) make sense as well.
Moreover, Bridgeland proved that there is a derived equivalence $\Phi$ (with inverse $\Psi$) between $Y$ and $Y^+$ restricting to an equivalence
\begin{align*}
	\Phi\colon \perv{-1}{\acat{A}^+} \leftrightarrows \perv{0}{\acat{A}} :\! \Psi
\end{align*}
which is the key to transport DT invariants from one side of the flop to the other.

The following lemma will be useful.
\begin{lem}\label{structuresheaf}
	Assume to be working in Situation \ref{situflop2}, then $\Phi(\O_{Y^+}) = \O_{Y}$.
\end{lem}
\begin{prf}
	First of all, it is shown in \cite[(4.4)]{tom-flops} that the equivalence $\Phi$ commutes with pushing down to $X$.
	The object $\Phi(\O_{Y^+}) =: L$ is a line bundle as, for any closed point $y \in Y$, the complex
	\begin{align*}
		\R\Hom_{Y}(\Phi(\O_{Y^+}),\O_y) = 
		\R\Hom_{Y^+}(\O_{Y^+},\Psi(\O_y)) = 
		\R\Hom_X(\O_X,\R f^+_* \Psi(\O_y)) = 
		\R\Hom_X(\O_X,Rf_*\O_y)
	\end{align*}
	 is concentrated in degree zero and has dimension one.
	The bundle $L$ pushes down to the structure sheaf, $\R f_* L = \O_X$.
	By adjunction, morphisms $\O_X \to \R f_* L$ correspond to morphisms $\O_Y \to L$, so that we deduce the existence of a non-zero section of $L$.
	Using Grothendieck duality for $f$, we see that $\R f_* L^\vee = \R f_* \R\lHom(L, f^!\O_X) = (\R f_* L)^\vee = \O_X$, hence $L^\vee$ has a non-zero section as well.
	As $Y$ is proper and integral, it follows that $L = \Phi(\O_{Y^+})$ must be the structure sheaf $\O_{Y}$.
\end{prf}

Gathering all the results so far, the only task left to accomplish is to compare the generating series for the perverse DT invariants on both sides of the flop: $\correct{\DT}(Y/X)$, $\correct{\DT}(Y^+/X)$.
%

The functor $\Phi$ induces an isomorphism between the numerical $K$-groups of $Y$ and $Y^+$, which restricts to an isomorphism
\begin{align}\label{identification}
	\phi\colon N_{\leq 1}(Y^+) \leftrightarrows N_{\leq 1}(Y) :\! \psi.
\end{align}
We can sharpen this result, by noticing that a class $(\beta,n) \in N_{\leq 1}(Y^+)$ is sent to $(\varphi(\beta),n)$, where $\varphi$ can be described as follows.
The smooth locus $U$ of $X$ is a common open subset of both $Y$ and $Y^+$.
By the Gysin exact sequence, we have an identification between the numerical groups of divisors of $Y$ and $Y^+$, via pulling back to $U$.
The inverse of the transpose of this identification is precisely $\varphi$, as the equivalence $\Phi$ restricts to the identity on $U$.

%
As the Fourier-Mukai equivalence $\Phi$ is an exact functor, and in light of Lemma \ref{structuresheaf}, we deduce an isomorphism of perverse Hilbert schemes $\perv{-1}{\Hilb(Y^+/X)} \simeq \perv{0}{\Hilb(Y/X)}$.
We can sharpen this result by noticing that, for a class $(\beta,n) \in N_1(Y^+) \oplus \Z$, we have
\begin{align*}
	\perv{-1}{\Hilb_{Y^+/X}(\beta,n)} \simeq \perv{0}{\Hilb_{Y/X}(\varphi(\beta),n)}.
\end{align*}
Taking weighted Euler characteristics and summing over all $\beta$'s and $n$'s we obtain
\begin{align*}
	\sum_{\beta,n} \correct{\DT}_{Y^+/X}(\beta,n) q^{(\beta,n)} = \sum_{\beta,n} \correct{\DT}_{Y/X}(\varphi(\beta),n) q^{(\beta,n)}
\end{align*}
which can be rephrased as a theorem.
\begin{teo}\label{pervdtlemma}
	Assume to be working in Situation \ref{situflop2}.
	Then, identifying variables via $\phi$, the following identity holds.
	\begin{align*}
		\correct{\DT}(Y^+/X) = \correct{\DT}(Y/X)
	\end{align*}
\end{teo}
The identity \eqref{superid} promised in the introduction now follows.
\begin{cor}\label{maino}
	The following identity holds.
	\begin{align*}
		\DT^\vee_{\exc}(Y^+)\DT(Y^+) = \DT^\vee_{\exc}(Y)\DT(Y)
	\end{align*}
\end{cor}
Concretely, for a class $\beta = (\gamma,\delta) \in N_1(X) \oplus N_1(Y/X)$ and an integer $n$ we have
\begin{align*}
	\sum_{\substack{\delta_1 + \delta_2 = \delta \\ n_1 + n_2 = n}} \DT_{Y^+}(0,-\delta_1,n_1) \DT_{Y^+}(\gamma,\delta_2,n_2) - \DT_Y(0,-\varphi(\delta_1),n_1) \DT_Y(\gamma,\varphi(\delta_2),n_2) = 0.
\end{align*}
\begin{prf}
	The hard work is done, as we already have Theorem \ref{mainmanor}.
	To prove this last identity we first observe that $\DT_0(Y)$ is an expression depending only on the topological Euler characteristic of $Y$ \cite{behrend-fantechi}.
	A result of Batyrev \cite{batyrev} tells us that $\chi_{\topp}(Y) = \chi_{\topp}(Y^+)$, so that the combination of Theorem \ref{pervdtlemma}, Theorem \ref{mainmanor} and Remark \ref{flops signs remark} imply the desired identity.
\end{prf}

	\section{One-dimensional Singular Locus} 
\label{sec:one_dimensional_singular_locus} 
Let's revert to Situation \ref{situation}.
As we've mentioned earlier, divisors on $Y$ are now allowed to be contracted to curves on $X$.
In particular, there are elements $F \in \perv{p}{\acat{F}}$ with $\dim \supp F = 2$.
From this it follows that there are perverse coherent sheaves $E \in \perv{p}{\acat{A}}$, which are numerically curve classes (namely, $\ch_0(E) = \ch_1(E) = 0$) but with $\dim \supp E = 2$ (recall that the support of a complex is defined to be union of the supports of its cohomology sheaves).
In particular, such an $E$ will have $H^{-1}(E)[1]$ supported on a surface.
Thus, the category of perverse coherent sheaves numerically supported on a curve is not closed under subobjects (nor quotients).

Consider now $\perv{p}{\acat{A}_{\leq 1}}$, the category of $E \in \perv{p}{\acat{A}}$ with $\dim \supp E \leq 1$.
What makes things worse is that $\perv{p}{\acat{A}_{\leq 1}}$ is not closed under subobjects either!
This fact is the core reason why we had to replace $\curly{H}$ with $\curly K$ and is what causes the appearance of the ``partial'' DT numbers.

It is not clear whether there is a tidy formula in the Hall algebra relating $\perv{p}{\curly{H}}$ with $\curly{H}$.
To complicate matters further, the identity (in $H_\infty$) $1_{\perv{p}{\acat{A}}}^\O = 1_{\perv{p}{\acat{F}[1]}}^\O * 1_{\perv{p}{\acat{T}}}^\O$ does not even hold in the full Hall algebra: given $E \in \perv{p}{\acat{A}}$, together with its torsion and torsion free part $F[1], T$, there is an exact sequence
$$ 0 \to \Hom(\O_Y, F[1]) \to \Hom(\O_Y,E) \to \Hom(\O_Y, T) \to \Hom(\O_Y,F[2])$$
where the last group is equal to $H^2(Y,F) = H^1(X,R^1f_*F)$.
This group vanishes when $F \in \perv{p}{\acat{F}_{\leq 1}}$ (this assumption, in conjunction with $f_*F=0$, forces the support of $F$ to be a union of curves contracted to points, hence the support of $R^1f_*$ is zero-dimensional), but in general it may not be true.

\medskip
This being said, let us come to the good news.
If we restrict to the subcategory $\perv{p}{\acat{A}_{\exc}}$ then we can work around these obstacles.
\begin{lem}\label{evvai}
	Let $E \in \perv{p}{\acat{A}}$ be such that $\ch_0(E) = \ch_1(E) = 0$ and $\dim \supp Rf_*E = 0$.
	Then $\dim \supp E \leq 1$, in other words $E \in \perv{p}{\acat{A}_{\leq 1}}$.
	
	In other words, given $E \in \perv{p}{\acat{A}}$ with $Rf_*E$ a skyscraper, $E \in \perv{p}{\acat{A}_{\leq 1}}$ if and only if $\ch_0(E) = \ch_1(E) = 0$ (i.e.~for these complexes, being supported in dimension at most one is a condition on their chern characters).
\end{lem}
Before we prove this lemma, we recall the key technical result of Van den Bergh \cite[Lemma 3.1.3, 3.1.5]{vdb}.
\begin{lem}[Van den Bergh]
	Consider the counit morphism $f^*f_*T \to T$.
	The objects in $\perv{-1}{\acat T}$ are precisely those $T \in \Coh(Y)$ such that $f^*f_* T \to T$ is surjective.
	
	Given $F \in \Coh(Y)$, there is a canonical map $\phi_F\colon F \to H^{-1}(f^!R^1f_*F)$.
	The objects in $\perv{0}{\acat F}$ are precisely those $F \in \Coh(Y)$ such that $\phi_F$ is injective.
\end{lem}
\begin{prf}[of Lemma \ref{evvai}]
	Let us now prove that a perverse coherent sheaf $E$ satisfying our assumptions is actually supported on a curve.
	Let $T = H^0(E)$ and $F = H^{-1}(E)$.
	As $Rf_*E$ is a skyscraper sheaf, it follows that both $f_*T$ and $R^1f_*F$ are skyscraper sheaves as well.
	
	As $[E] = [T] - [F]$ in $K_0(Y)$ and $\ch_0(E) = \ch_1(E) = 0$, it follows that $T$ is supported in dimension at most one if and only if $F$ is.
	When the perversity is $p=-1$, we have that $ f^{-1}(\supp f_*T) = \supp f^*f_*T \supset \supp T$.
	But, as $f_*T$ is a skyscraper sheaf and the fibres of $f$ are at most one-dimensional, we have $\dim \supp f^*f_*T \leq 1$ and we can conclude.
	
	When the perversity is $p=0$, we have $\supp F \subset \supp H^{-1}(f^! R^1f_*F) \subset \supp f^! R^1f_*F \subset f^{-1}(\supp R^1f_*F)$.
	Again, as $R^1f_*F$ is a skyscraper sheaf, we are done.
\end{prf}
What truly makes the category $\perv{p}{\acat{A}_{\exc}}$ robust is the following lemma.
\begin{lem}\label{thick}
	The category $\perv{p}{\acat{A}_{\exc}}$ is closed under extensions, quotients and subobjects.
\end{lem}
\begin{prf}
	Let $A \to B \to C$ be a short exact sequence in $\perv{p}{\acat{A}}$.
	First of all, $Rf_*B$ is a skyscraper if and only if both $Rf_*A$ and $Rf_*C$ are skyscraper sheaves.
	As $[B] = [A] + [C]$ in $K_0(\perv{p}{\acat{A}})$, it follows from Lemma \ref{evvai} that $\perv{p}{\acat{A}_{\exc}}$ is closed under extensions.
	
	Thus we are left to show that if $B \in \perv{p}{\acat{A}_{\exc}}$ then $A,C \in \perv{p}{\acat{A}_{\exc}}$.
	Consider the long exact sequence of cohomology sheaves.
	$$ 0 \to H^{-1}(A) \to H^{-1}(B) \to H^{-1}(C) \to H^0(A) \to H^0(B) \to H^0(C) \to 0$$
	As $B \in \perv{p}{\acat{A}_{\exc}}$ it follows that $H^{-1}(A)$ and $H^0(C)$ are both supported in dimension at most one.
	Thus the only obstruction to concluding is the support of either $H^{-1}(C)$ or $H^0(A)$.
	However, we know that both $R^1f_*H^{-1}(C)$ and $f_*H^0(A)$ are skyscraper sheaves.
	Using Van den Bergh's lemma we conclude that at least one is (and hence both are) supported in dimension at most one.
\end{prf}
Consequently, the previous section can be adapted to the category $\perv{p}{\acat{A}_{\exc}}$ without needing any ``partial'' invariants.
\begin{teo}
	Assume Situation \ref{situation}.
	Then
	$$ \perv{p}{\correct{DT}}_{\exc}(Y/X) = \frac{ \correct{DT}^\vee_{\exc} \cdot \correct{DT}_{\exc}(Y) } { DT_0(Y) } $$
	holds, where
	$$
		\perv{p}{\correct{DT}}_{\exc}(Y/X) = \sum_{\substack{\beta,n \\ f_*\beta = 0}} \correct{DT}_{Y/X}(\beta,n) q^{(\beta,n)} ;
		\quad \quad
		\correct{DT}_{Y/X}(\beta,n) = \chi_\mu \left( \perv{p}{\Hilb}_{Y/X}(\beta,n) \right).$$
\end{teo}

Finally, we conclude by reiterating Remark \ref{addedremark}.
Let $\perv{p}\Hilb_{\ch_{\leq 1}}$ parameterise perverse quotients $\O_Y \onto E$ with $\ch_0(E) = 0 = \ch_1(E)$.
Although we expect it to be true, we do not whether $\perv{p}\Hilb_{\leq 1}$ is strictly contained inside $\perv{p}{\Hilb}_{\ch_{\leq 1}}$.
We feel that this matter should be investigated further, especially in light of the formulae deduced in \cite{cala}.

\appendix
\section{Substacks} 
\label{sec:substacks}
At the core of the construction of the Hall algebra of an abelian category lies the existence of a moduli stack\footnote{The author would like to thank Fabio Tonini for patiently explaining to him many things about stacks.} parameterising its objects (and a moduli of short exact sequences).
In our case this amounts, first of all, to proving the existence of the moduli stack $\perv{p}{\astack{A}}$, parameterising perverse coherent sheaves.
We have mentioned in the first section that as the category $\perv{p}{\acat{A}}$ is the heart of a t-structure, its objects have no negative self-extensions.
This simple remark is actually key, as we construct $\perv{p}{\astack{A}}$ as an open substack of one big moduli stack $\stack{Mum},$ which Lieblich refers to as \emph{the mother of all moduli of sheaves} \cite{lieblich}.
Let us recall its definition.

First, fix a flat and proper morphism of schemes $\pi: X \to S$.
\begin{defn}
	An object $\cmplx{E} \in \text{D}(\curlyO_X)$ is (relatively over $S$) \emph{perfect and universally gluable} if the following conditions hold.
	\begin{itemize}
		\item There exists an open cover $\{U_i\}$ of $X$ such that $\cmplx{E}\vert_{U_i}$ is quasi-isomorphic to a bounded complex of quasi-coherent sheaves flat over $S$.
		\item For any $S$-scheme $u:T \to S$ we have
		\begin{align*}
			R\pi_{T,*}R\lHom_{X_T}(Lu_X^*\cmplx{E},Lu_X^*\cmplx{E}) \in \text{D}^{\geq 0}(\curlyO_T)
		\end{align*}
		where $\pi_T$ and $u_X$ denote the maps induced by $\pi$ and $u$ respectively on the base-change $X_T$.
	\end{itemize}
	We denote the category of perfect and universally gluable sheaves on $X$ (over $S$) as $\text{D}_\text{pug}(\curlyO_X).$
\end{defn}
If in the definition we take $S$ to be affine and assume $T = S$, then it's clear that gluability has to do with the vanishing of negative self-exts of $\cmplx{E}$.
This condition is necessary to avoid having to enter the realm of higher stacks.

A prestack\footnote{We are using the term \emph{prestack} in analogy with term \emph{presheaf}.} $\stack{Mum}_X$ is defined by associating with an S-scheme $T \to S$ (the associated groupoid of) the category $\text{D}_\text{pug}(\curlyO_{X_T})$ of perfect and universally gluable complexes (relatively over $T$).
The restriction functors are defined by derived pullback.
\begin{teo}[Lieblich]
	The prestack $\stack{Mum}_X$ is an Artin stack, locally of finite presentation over $S.$
\end{teo}
From now on we fix $\pi: X \to S$ flat and projective with $S$ a noetherian scheme.
We assume all rings and schemes to be locally of finite type over $S$.\footnote{For what follows, this assumption isn't substantial (as $\mum_X$ is locally of finite type over $S$) but it enables us to use the local criterion of flatness directly. This is essentially a consequence of \cite[Corollaire (10.11) (ii)]{laumon}.}

We want to construct various open substacks of $\stack{Mum}_X,$ namely stacks of complexes satisfying additional properties.
For example we would like to construct the stack of complexes with cohomology concentrated in degrees less or equal than a fixed integer $n.$
The correct way to proceed is by imposing conditions fibrewise on restrictions to geometric points.
Let us illustrate a general recipe first.
The following diagram comes in handy.
\begin{center}
	\begin{tikzpicture}
		\matrix (m) [matrix of math nodes, row sep=3em, column sep=3em, text height=1.5ex, text depth=0.25ex]
		{
		X_t & X_T & X \\
		\Spec k & T & S \\
		};
		\path[->,font=\scriptsize]
		(m-1-1) edge node[auto]{$\pi_t$} (m-2-1)
				edge node[auto]{$t_X$} (m-1-2)
		(m-1-2) edge node[auto]{$\pi_T$} (m-2-2)
				edge node[auto]{$u_X$} (m-1-3)
		(m-1-3) edge node[auto]{$\pi$} (m-2-3)
		(m-2-1) edge node[auto]{$t$} (m-2-2)
		(m-2-2) edge node[auto]{$u$} (m-2-3)
		;
	\end{tikzpicture}
\end{center}
Here $T$ is the base space for our family of complexes, together with its structure map to $S,$ and $t \in T$ is a geometric point.
Given a property P, we might define the stack of complexes satisfying $P$ as follows.
\begin{align*}
	\stack{Mum}_X^P(T) = \left\{ \cmplx{E} \in \stack{Mum}_X(T) \st \forall \text{ geometric } t\in T, \cmplx{E}\vert^L_{X_t} \text{ satisfies } P \right\}
\end{align*}
We recall that by $\cmplx{E}\vert^L_{X_t}$ we mean $Lt_X^* \cmplx{E}.$

To construct the substacks of $\mum_X$ we are interested in we make use of the following lemma.
\begin{lem}\label{superlemma}
	Let $T \to S$ be an $S$-scheme, let $t: \Spec k \to T$ be a point of $T$ and let $\cmplx{E} \in \text{D}^\text{b}(\curlyO_{X_T})$ be a bounded complex of $\curlyO_{X_T}$-modules flat over $T$.
	Let $n \in \Z$ be an integer.
	The following statements hold.
	\begin{enumerate}
		\item $\cmplx{E}\vert^L_{X_t} \in \text{D}^{\leq n}(\curlyO_{X_t}) \longiff X_t \subset U_>$, where
		\begin{align*}
			U_> &= \bigcap_{q > n} X_T \setminus \supp H^q\left( \cmplx{E} \right).
		\end{align*}
		\item $\cmplx{E}\vert^L_{X_t} \in \text{D}^{[n]}(\curlyO_{X_t}) \longiff X_t \subset U$, where
		\begin{align*}
			U &= U_> \cap U_\text{f} \cap U_< \\
			U_> &= \bigcap_{q > n} X_T \setminus \supp H^q(\cmplx{E}) \\
			U_\text{f} &= \left\{ x \in X_T \st H^n(E)_x \text{ is a flat } \curlyO_{T,\pi_T(x)}\text{-module} \right\}\\
			U_< &= \bigcap_{q < n} X_T \setminus \supp H^q(\cmplx{E}).
		\end{align*}
		\item $\cmplx{E}\vert^L_{X_t} \in \text{D}^{\geq n}(\curlyO_{X_t}) \longiff \cmplx{F} \in \text{D}^{[n]}(\curlyO_{X_t})$, where $\cmplx{F} = \sigma_{\leq n}\cmplx{E}$ is the stupid truncation of $\cmplx{E}$ in degrees less or equal than $n$.
		\begin{align*}
			\sh{F}^p =
			\begin{cases}
				\sh{E}^p, \text{ if } p\leq n \\
				0, \text{ if } p > n
			\end{cases}
		\end{align*}
	\end{enumerate}
\end{lem}
\begin{prf}
	{\scshape Proof of 1.}
	Let $t_X$ be the inclusion of the fibre $X_t \to X_T$.
	As $t_X$ is an affine map we do not lose information on the cohomologies of $\cmplx{E}\vert^L_{X_t}$ after pushing forward back into $X_T$.
	We also have isomorphisms
	\begin{align*}
		t_{X,*}\cmplx{E}\vert^L_{X_t} \simeq
		\cmplx{E} \stackrel{L}{\otimes}_{\curlyO_{X_T}}t_{X,*}\curlyO_{X,t}
		\simeq \cmplx{E} \stackrel{L}{\otimes}_{\curlyO_{X_T}} \pi_T^*t_*k
	\end{align*}
	where the first follows from the projection formula and the second from base change compatibility.
	As we are interested in the vanishing of $H^q(\cmplx{E}\vert^L_{X_t})$ we may restrict to the stalk at a point $x \in X_t$.
	Taking stalks at $x$ gives us isomorphisms
	\begin{align}\label{globaltor}
		H^q\left(\cmplx{E}\vert^L_{X_t}\right)_x \simeq H^q\left( \cmplx{E}_x \stackrel{L}{\otimes}_{\curlyO_{T,t}} k \right).
	\end{align}
	We have the page two spectral sequence of the pullback
	\begin{align}\label{pullback spectral sequence}
		L^pt_X^* H^q(\cmplx{E}) \longimplica H^{p+q}(\cmplx{E}\vert_{X_t}^L).
	\end{align}
	which, at a point $x \in X_t$ and using the isomorphism \eqref{globaltor}, boils down to
	\begin{align}\label{torss}
		\Tor_{-p}^{\curlyO_{X_t}}\left(H^q\left(\cmplx{E}\right)_x,k\right)
		\longimplica H^{p+q}\left( \cmplx{E}\vert_{X_t}^L \right)_x.
	\end{align}
	Let now $q$ be the largest integer such that $H^q(\cmplx{E}) \neq 0$.
	From the spectral sequence \eqref{torss} we have
	\begin{align*}
		H^q\left(\cmplx{E}\vert^L_{X_t}\right)_x \simeq
		H^q(\cmplx{E})_x \otimes_{\curlyO_{T,t}} k.
	\end{align*}
	Hence, by Nakayama, $H^q(\cmplx{E}\vert_{X_t}^L)_x = 0$ if and only if $x \in X_T \setminus \supp H^q(\cmplx{E})$ and finally
	$$H^q(\cmplx{E}\vert_{X_t}^L)=0 \longiff X_t \subset X_T \setminus \supp H^q(\cmplx{E}).$$
	
	{\scshape Proof of 2.} Using 1. we can assume that $\cmplx{E}\vert_{X_t}^L \in \text{D}^{\leq n}(\curlyO_{X_t})$.
	By the spectral sequence \eqref{pullback spectral sequence} we have that $H^{n-1}(\cmplx{E}\vert_{X_t}^L) \simeq L_1t_X^*H^n(\cmplx{E})$.
	Again, we may pass on to the stalk at a point $x \in X_t$ and \eqref{torss} yields
	\begin{align*}
		H^{n-1}\left(\cmplx{E}\vert_{X_t}^L\right)_x \simeq 
		\Tor_1^{\curlyO_{X_t}}\left(\cmplx{E}_x,k\right)
	\end{align*}
	the vanishing of which is equivalent, by the local criterion for flatness, to $H^q(\cmplx{E})_x$ being a flat $\curlyO_{X,t}$-module.
	
	We can thus assume that $X_t \subset U_> \cap U_\text{f}$.
	Once more, from the spectral sequence \eqref{torss} we have that
	$H^{n-1}(\cmplx{E}\vert_{X_t}^L) \simeq t_X^*H^{n-1}(\cmplx{E})$ and we proceed as in the proof of 1.
	
	{\scshape Proof of 3.} Consider the page one spectral sequence
		\begin{align*}
			L^q t_X^* \sh{E}^p \longimplica H^{p+q}\left(\csh{E}\vert^L_{X_t}\right)
		\end{align*}
		from which we get isomorphisms
		\begin{align*}
			H^p\left(\csh{E}\vert^L_{X_t}\right) \simeq H^p\left(t_X^*\csh{E}\right)
		\end{align*}
		as a consequence of flatness of the $\sh{E}^q$'s.
		Thus, for $p < n$,
		$$
			H^p\left(\cmplx{E}\vert^L_{X_t}\right)=0 \longiff
			H^p\left(t_X^*\cmplx{E}\right) = 0 \longiff
			H^p\left(t_X^*\cmplx{F}\right) = 0.
		$$
\end{prf}

\begin{prop}\label{mumlessthann}
	Define the prestack $\mum_X^{\leq n} = \stack{Mum}_X^{[-\infty,n]}$ by assigning to each $S$-scheme $T$ the groupoid
	\begin{align*}
		\stack{Mum}_X^{\leq n}(T) = \left\{ \cmplx{E} \in \stack{Mum}_X(T) \st \forall \text{ geometric } t \in T, \cmplx{E}\vert_{X_t}^L \in \text{D}^{\leq n}(\curlyO_{X_t}) \right\}
	\end{align*}
	with restriction functors induced by $\stack{Mum}_X.$
	The prestack $\stack{Mum}_X^{\leq n}$ is an open substack of $\stack{Mum}_X.$
\end{prop}
\begin{prf}
	That $\stack{Mum}_X^{\leq n}$ satisfies descent is a direct consequence of descent for $\stack{Mum}_X.$
	To prove that it is indeed an open substack it is sufficient to prove that for any affine $S$-scheme $T$, together with a morphism $T \to \stack{Mum}_X$ corresponding to a complex $\cmplx{E} \in \stack{Mum}_X(T)$, the set
	\begin{align*}
		V = \left\{ t \in T \st \cmplx{E}\vert_{X_t}^L \in \text{D}^{\leq n}(X_t) \right\}
	\end{align*}
	is an open subset of $T.$
	
	By Lemma \ref{superlemma} 1.~we know that $t \in V$ if and only if $X_t \subset U_>$ (notice that by our assumptions the complex $\cmplx{E}$ is bounded).
	Thus $\pi_T(X_T \setminus U_>) = \pi_T(X_T)\setminus V$.
	The set $U_>$ is open as the sheaves $H^q(\cmplx{E})$ are quasi-coherent and of finite type.
	Finally, the sets $\pi_T(X_T)$ and $\pi_T(X_T \setminus U_>)$ are closed, being the image of closed subsets under a proper map.
	Thus, $V$ is open.
\end{prf}
Notice that the condition of being concentrated in degrees less or equal than $n$ is in fact a global condition, i.e.~we could have requested $\cmplx{E}\in \text{D}^{\leq n}(\curlyO_{X_T})$ directly.

We now impose on our complexes the further condition of being concentrated in a fixed degree $n\in \Z.$
This stack will be isomorphic to the stack of coherent sheaves shifted by $-n$.
\begin{prop}\label{stackcoh}
	Define the prestack $\stack{Mum}_X^{[n]}$ by assigning to each $S$-scheme $T$ the groupoid
	\begin{align*}
		\mum_X^{[n]}(T) = \left\{ \cmplx{E} \in \mum_X^{\leq n}(T) \st \forall t\in T, \cmplx{E}\vert_{X_t}^L \in \text{D}^{[n]}(\curlyO_{X_t}) \right\}
	\end{align*}
	with restriction functors induced by $\mum_X.$
	The prestack $\mum_X^{[n]}$ is an open substack of $\mum_X^{\leq n}$.
\end{prop}
\begin{prf}
	The proof follows along the lines as the previous one.
	It suffices to show that for any affine scheme $T$, together with a map $T \to \mum_X^{\leq n}$ corresponding to a complex $\cmplx{E} \in \mum_X^{\leq n}(T)$, the set
	\begin{align*}
		V = \left\{ t \in T \st \cmplx{E}\vert_{X_t}^L \in \text{D}^{[n]}(\curlyO_{X_t}) \right\}
	\end{align*}
	is an open subset of $T.$
	By Lemma \ref{superlemma} 2.~we know that $t \in V$ if and only if $X_t \subset U$.
	The sets $U_<, U_>$ are open as the sheaves $H^q(\cmplx{E})$ are quasi-coherent and of finite type.
	The set $U_\text{f}$ is open by the open nature of flatness \cite[Th\'eor\`eme 11.3.1]{egaiv-3}.
	Thus $U$ is open and we conclude as in the previous proof.
\end{prf}
\noindent When $n=0$ we get back the ordinary stack of coherent sheaves on $X$.

We now turn to the opposite condition: being concentrated in degrees greater or equal than a fixed $n\in \Z.$
\begin{prop}
	Define the prestack $\mum_X^{\geq n} = \mum_X^{[n,\infty]}$ by assigning to each $S$-scheme $T$ the groupoid
	\begin{align*}
		\mum_X^{\geq n}(T) =\left\{ \cmplx{E} \in \mum_X(T) \st \forall t\in T, \cmplx{E}\vert_{X_t}^L \in \text{D}^{\geq n}(\curlyO_{X_t}) \right\}
	\end{align*}
	with restriction functors induced by $\mum_X.$
	The prestack $\mum_X^{\geq n}$ is an open substack of $\mum_X.$
\end{prop}
\begin{prf}
	As in the previous proofs we consider a complex $\cmplx{E}\in \mum_X(T)$ corresponding to a morphism $T \to \mum_X$ and prove that the set
	\begin{align*}
		V = \left\{ t \in T \st \cmplx{E}\vert_{X_t}^L \in \text{D}^{\geq n}(X_t) \right\}
	\end{align*}
	is an open subset of $T.$
	By Lemma \ref{superlemma} 3.~this set is equal to
	\begin{align*}
		V = \left\{ t \in T \st \cmplx{F}\vert_{X_t}^L \in \text{D}^{[n]}(X_t) \right\}
	\end{align*}
	which is open by the previous proof.
\end{prf}

\bibliographystyle{fabio.bst}
\bibliography{biblio.bib}
\end{document}